\documentclass[12pt]{amsart} 
\usepackage{amsmath,amsfonts,amssymb,amscd}
\usepackage[all]{xy}
\usepackage[francais]{babel}

\newenvironment{pf*}{\begin{proof}}{\end{proof}} 
\newcommand{\ead}[1]{{\it E-mail~:~}{\tt #1}}
\newtheorem{theoreme}[equation]{Th\'eor\`eme}
\newtheorem{lemme}[equation]{Lemme}
\newtheorem{lemme-definition}[equation]{Lemme-D\'efinition}
\newtheorem{proposition-definition}[equation]{Proposition-D\'efinition}
\newtheorem{proposition}[equation]{Proposition}
\newtheorem{definition}[equation]{D\'efinition}
\newtheorem{corollaire}[equation]{Corollaire}
\newtheorem{conjecture}[equation]{Conjecture}
\theoremstyle{remark}
\newtheorem{remarque}[equation]{Remarque}

\numberwithin{equation}{subsection}
\newcommand{\sub}[1]{
\setcounter{subsubsection}{\value{equation}}
\subsubsection{#1}
\setcounter{equation}{\value{subsubsection}}
}

\DeclareMathOperator{\ad}{ad}
\DeclareMathOperator{\Chow}{CH}
\DeclareMathOperator{\End}{End}
\DeclareMathOperator{\Hom}{Hom}
\DeclareMathOperator{\Ext}{Ext}
\DeclareMathOperator{\Gal}{Gal}
\DeclareMathOperator{\Id}{Id}
\DeclareMathOperator{\Ind}{Ind}
\DeclareMathOperator{\Irr}{Irr}
\DeclareMathOperator{\Res}{Res}
\DeclareMathOperator{\Sh}{Sh}
\DeclareMathOperator{\Spec}{Spec}
\DeclareMathOperator{\St}{St}
\DeclareMathOperator{\Trace}{Trace}
\newcommand{\bbA}{{\mathbb A}}
\newcommand{\bbC}{{\mathbb C}}
\newcommand{\bbF}{{\mathbb F}}
\newcommand{\bbG}{{\mathbb G}}
\newcommand{\bbN}{{\mathbb N}}
\newcommand{\bbQ}{{\mathbb Q}}
\newcommand{\bbR}{{\mathbb R}}
\newcommand{\bbT}{{\mathbb T}}
\newcommand{\bbZ}{{\mathbb Z}}
\newcommand{\bB}{{\bf B}}
\newcommand{\bF}{{\bf F}}
\newcommand{\bG}{{\bf G}}
\newcommand{\bH}{{\bf H}}
\newcommand{\bI}{{\bf I}}
\newcommand{\bJ}{{\bf J}}
\newcommand{\bL}{{\bf L}}
\newcommand{\bO}{{\bf O}}
\newcommand{\bP}{{\bf P}}
\newcommand{\bS}{{\bf S}}
\newcommand{\bT}{{\bf T}}
\newcommand{\bU}{{\bf U}}
\newcommand{\bV}{{\bf V}}
\newcommand{\bW}{{\bf W}}
\newcommand{\bX}{{\bf X}}
\newcommand{\bY}{{\bf Y}}
\newcommand{\bZ}{{\bf Z}}
\newcommand{\ba}{{\bf a}}
\newcommand{\bb}{{\bf b}}
\newcommand{\bd}{{\bf d}}
\newcommand{\bs}{{\bf s}}
\newcommand{\bt}{{\bf t}}
\newcommand{\bu}{{\bf u}}
\newcommand{\bv}{{\bf v}}
\newcommand{\bw}{{\bf w}}
\newcommand{\bx}{{\bf x}}
\newcommand{\by}{{\bf y}}
\newcommand{\bz}{{\bf z}}
\newcommand{\bpi}{{\boldsymbol\pi}}
\newcommand{\CB}{{\mathcal B}}
\newcommand{\CC}{{\mathcal C}}
\newcommand{\cD}{{\mathcal D}} 
\newcommand{\CE}{{\mathcal E}}

\newcommand{\CH}{{\mathcal H}}
\newcommand{\CL}{{\mathcal L}}
\newcommand{\CO}{{\mathcal O}}
\newcommand{\CR}{{\mathcal R}}
\newcommand{\CS}{{\mathcal S}}
\newcommand{\CT}{{\mathcal T}}
\newcommand{\us}{{\underline s}}
\newcommand{\ut}{{\underline t}}
\newcommand{\uv}{{\underline v}}
\newcommand{\uw}{{\underline w}}
\newcommand{\uy}{{\underline y}}
\newcommand{\uB}{{\underline B}}
\newcommand{\uH}{{\underline H}}
\newcommand{\uS}{{\underline S}}
\newcommand{\uW}{{\underline W}}
\newcommand{\ubs}{{\underline\bs}}
\newcommand{\ubt}{{\underline\bt}}
\newcommand{\ubv}{{\underline\bv}}
\newcommand{\ubw}{{\underline\bw}}
\newcommand{\uby}{{\underline\by}}
\newcommand{\ubz}{{\underline\bz}}
\newcommand{\ubS}{{\underline\bS}}
\newcommand{\ubW}{{\underline\bW}}
\newcommand{\dz}{{\dot z}}
\newcommand{\Fq}{{\bbF_q}}
\newcommand{\Fqbar}{{\overline\bbF_q}}
\newcommand{\Qlbar}{{\overline{\mathbb Q}_\ell}}
\newcommand{\BW}{{B^+}}
\newcommand{\GF}{{\bG^F}}
\newcommand{\inv}{^{-1}}
\newcommand{\iso}{{\xrightarrow\sim}}
\newcommand{\cf}{{\it cf.}}
\newcommand{\eg}{{\it e.g.}}
\newcommand{\ie}{{\it i.e.}}
\newcommand{\etc}{{\it etc$\ldots$}}
\newcommand{\TrH}[2]{{\sum_i (-1)^i\Trace(#1|H^i_c(#2,\Qlbar))}}
\newcommand{\lexp}[2]{\kern\scriptspace\vphantom{#2}^{#1}\kern-\scriptspace#2}
\newcommand{\scal}[2]{\langle\,#1,#2\,\rangle}
\newcommand{\genby}[1]{\mathopen<#1\mathclose>}
\newcommand{\sdp}[1]{\rtimes\genby{#1}}
\newcommand{\seof}[3]{#1\xrightarrow o#2\xrightarrow f #3}
\newcommand{\casde}[1]{\par\noindent $\bullet$ Cas de $#1$}

\begin{document}
\title{Cohomologie des vari\'et\'es de Deligne-Lusztig}

\author{Fran\c cois Digne}
\address{LAMFA, Universit\'e de Picardie,
33 Rue Saint-Leu, 80039 Amiens France
\ead{digne@u-picardie.fr}}

\author{Jean Michel}
\address{LAMFA, Universit\'e de Picardie, 
33 Rue Saint-Leu, 80039 Amiens France\newline
and \newline
Institut de Math\'ematiques de Jussieu,
2 place Jussieu, 75251 Paris Cedex 05, France
\ead{jmichel@math.jussieu.fr}}

\author{Rapha\"el Rouquier}
\address{Institut de Math\'ematiques de Jussieu,
2 place Jussieu, 75251 Paris Cedex 05, France
\ead{rouquier@math.jussieu.fr}}

\renewcommand{\abstractname}{Abstract}
\begin{abstract}
  We study the cohomology of Deligne-Lusztig varieties with aim the
construction of actions of Hecke algebras on such cohomologies, as
predicted by the conjectures of Brou\'e, Malle and Michel ultimately
aimed at providing an explicit version of the abelian defect conjecture.
  We develop the theory for varieties associated to elements of the braid
monoid and partial compactifications of them. We are able to compute the
cohomology of varieties associated to (possibly twisted) rank $2$ groups and
powers of the longest element $w_0$ (some indeterminacies remain for $G_2$).
  We use this to construct Hecke algebra actions on the cohomology of
varieties associated to $w_0$ or its square, for groups of arbitrary rank.
  In the subsequent work \cite{DMc}, we construct actions associated to more
general regular elements and we study their traces on cohomology.
\end{abstract}
\maketitle 
\tableofcontents

\section{Introduction}

L'objet   de   cet   article   est  la   construction   d'actions   d'alg\`ebres
d'Iwahori-Hecke sur la cohomologie de certaines vari\'et\'es de Deligne-Lusztig.
L'existence de  ces actions  fait partie des  conjectures pr\'ecisant,  pour les
groupes r\'eductifs finis,  la conjecture de Brou\'e sur les  blocs \`a d\'efaut
ab\'elien des groupes finis.

Dans ce travail, nous \'etablissons des propri\'et\'es g\'en\'erales des
vari\'et\'es  de De\-ligne-Lusztig et de leur cohomologie et nous montrons
l'existence  des repr\'e\-sentations d'alg\`ebres de Hecke pour deux types
de vari\'et\'es de Deligne-Lusztig, celles associ\'ees \`a l'\'el\'ement
$\bpi$  du mono{\"\i}de  de tresses  et \`a  sa racine carr\'ee $\bw_0$.
Nous  y parvenons gr\^ace  \`a un calcul  de la cohomologie de certaines
vari\'et\'es  pour des groupes  de rang $2$  (celles associ\'ees \`a des
puissances quelconques de $\bw_0$).

Rappelons maintenant les conjectures auxquelles nous nous int\'eressons.

Soit $\bG$ un  groupe r\'eductif connexe sur  une cl\^oture alg\'ebrique
d'un corps fini, muni d'une isog\'enie $F$ dont une
puissance est un endomorphisme  de Frobenius.  Nous
notons  $\GF$ le  groupe (fini)  des points  fixes de  $F$. Soit  $W$ le
groupe de Weyl de $\bG$ et $\BW$ le mono\"\i de de tresses associ\'e
\`a $W$.
Pour $w\in W$, on dispose d'un relev\'e de longueur minimale $\bw$ dans
$\BW$.

Soit $\bpi=\bw_0^2$, o\`u $w_0$ est l'\'el\'ement de plus
grande longueur de $W$.
Soit $w\in W$ tel que  $\bw$ est une
\og $F$-racine $d$-i\`eme  de $\bpi$\fg, \ie,
$(\bw F)^d=\bpi  F^d$. Dans \cite{Sydney}, une action \`a droite
de $C_B(\bw F)$ sur $H^*_c(\bX(w),\Qlbar)$ est
construite (d\'eduite d'une action de $C_{B^+}(\bw F)$ sur
$\bX(w)$) et il est conjectur\'e que
\begin{itemize}
\item[(i)] l'action de $\Qlbar C_B(\bw F)$ sur $H^*_c(\bX(w),\Qlbar)$ se
factorise en une action d'une alg\`ebre de Hecke \og cyclotomique\fg\
relative \`a $C_W(w F)$.
\item[(ii)] $\Hom_{\Qlbar\bG^F}(H^i_c(\bX(w),\Qlbar),H^j_c(\bX(w),\Qlbar))
=0$ si $i\not=j$.
\end{itemize}

\medskip
On passe de la conjecture g\'en\'erale sur les $\ell$-blocs \`a d\'efaut
ab\'elien \`a cette conjecture en faisant trois restrictions:

\begin{itemize}
\item on suppose que le centralisateur d'un groupe de d\'efaut est un tore
(son type est alors un \'el\'ement r\'egulier)
\item on \'etend les scalaires de $\overline{\bbZ}_\ell$ \`a $\Qlbar$
\item on se restreint aux caract\`eres unipotents.
\end{itemize}

Les  conjectures  (i)  et  (ii)  devraient \^etre compl\'et\'ees par une
pr\'ediction des param\`etres de l'alg\`ebre de Hecke et de sa trace sur
la cohomologie.

Les  parties  \S  \ref{sectionvarietes}  et  \S \ref{sectioncohomologie}
d\'eveloppent  des  propri\'et\'es  g\'en\'erales  des  vari\'et\'es  de
Deligne-Lusztig  g\'en\'eralis\'ees et de leur cohomologie.  Nous reprenons
et compl\'etons des r\'esultats de Deligne-Lusztig et Lusztig.

Dans   \S  \ref{sectionvarietes},  nous  pr\'esentons  la  construction,
suivant  Deligne, de  vari\'et\'es associ\'ees  \`a des  \'el\'ements du
mono{\"\i}de  de tresses.  Plus g\'en\'eralement,  nous introduisons des
compactifications  partielles de  ces vari\'et\'es,  associ\'ees \`a des
\'el\'ements  d'un mono{\"\i}de \og  compl\'et\'e\fg. Nous \'etablissons
en  particulier une relation avec des vari\'et\'es pour des sous-groupes
de Levi (proposition \ref{X_w produit de varietes}).

Dans   \S  \ref{sectiondevissages},   nous  d\'eveloppons  diff\'erentes
techniques reliant la cohomologie de la vari\'et\'e $\bX(\bw)$ \`a celle
de vari\'et\'es $\bX(\bw')$ pour des \'el\'ements $\bw'$ plus courts que
$\bw$.   Les   techniques   utilis\'ees   reposent   sur   l'\'etude  de
d\'ecompositions  des vari\'et\'es, de certains morphismes propres et de
situations de lissit\'e rationnelle.

Nous   \'etudions   la   structure   comme  $(\Qlbar\GF)$-module  de  la
cohomologie    des   vari\'et\'es   $\bX(\bw)$   dans   la   partie   \S
\ref{sectionunipotentes}.   Nous  \'etudions  les   valeurs  propres  de
l'endomorphisme  de  Frobenius.
Nous  d\'ecrivons  la  partie  de  la
cohomologie  o\`u le groupe agit trivialement ou par la repr\'esentation
de  Steinberg.  Nous  \'etudions  les  \'el\'ements  $\bw$  de  longueur
minimale    tels   qu'une   repr\'esentation   irr\'eductible   donn\'ee
appara{\^\i}t dans $H^*_c(\bX(\bw),\Qlbar)$.

Nous \'etudions les propri\'et\'es de rationalit\'e des
caract\`eres des $H^*_c(\bX(\bw),\Qlbar)$ dans la partie \S \ref{sectionQ}.

Dans  la partie  \S \ref{sectionrang2},  nous d\'eterminons la structure
comme   $(\Qlbar\GF)$-module   de   la   cohomologie   des  vari\'et\'es
$\bX(\bw_0^n)$  pour  des  groupes  de  type $A_2$, $\lexp 2A_2$, $B_2$,
$\lexp  2 B_2$  et $\lexp  2 G_2$  (nous obtenons  aussi des r\'esultats
partiels  pour $G_2$). Nous proc\'edons  par r\'ecurrence et nous devons
calculer   pour   cela   la   cohomologie   de   certaines  vari\'et\'es
$\bX(\bw\bw_0^{2m})$.  \`A d\'ecalage  et twist  de Tate  pr\`es (qui ne
d\'ependent  que  de  la  famille  de la repr\'esentation irr\'eductible
concern\'ee),  le r\'esultat ne d\'epend pas de $m$ et nous conjecturons
(\cf\  \S \ref{periodicite}) que de tels ph\'enom\`enes de p\'eriodicit\'e
sont  g\'en\'eraux.  Ces  calculs  sont  rendus possibles par les outils
d\'evelopp\'es au \S \ref{sectioncohomologie}.

La  derni\`ere partie \S \ref{sectionendo} est consacr\'ee \`a l'\'etude
de   certains   endomorphismes   des   vari\'et\'es   $\bX(\bw)$.   Elle
g\'en\'eralise  des  r\'esultats  de  \cite{LuMa} et \cite{Sydney}. Nous
expliquons   comment  l'\'etude  d'endomorphismes   associ\'es  \`a  des
\'el\'ements contenus dans un sous-groupe parabolique se ram\`ene au cas
de  vari\'et\'es associ\'ees au sous-groupe  de Levi correspondant. Nous
utilisons   ceci  pour  d\'emontrer   que  certains  endomorphismes  des
vari\'et\'es  $\bX(\bpi)$  et  $\bX(\bw_0)$  v\'erifient  des  relations
quadratiques  et nous en d\'eduisons  l'action d'alg\`ebres de Hecke sur
la cohomologie (conform\'ement \`a la conjecture (i)).

Dans une suite de ce travail \cite{DMc}, la conjecture (i) est \'etablie
pour   diverses  classes  d'\'el\'ements  r\'eguliers  et  la  trace  de
l'alg\`ebre de Hecke cyclotomique sur la cohomologie est \'etudi\'ee.

Nous remercions Luc Illusie, Bruno Kahn et G\'erard Laumon pour des
discussions utiles \`a ce travail.

\section{Vari\'et\'es de Deligne-Lusztig et groupes de tresses}
\label{sectionvarietes}

\subsection{Groupes de tresses}
\sub{}
Soit  $(W,S)$ un  syst\`eme de  Coxeter fini,  c'est-\`a-dire, un groupe
fini  $W$  et  une  partie  g\'en\'eratrice  $S$  de  $W$  telle que, si
$m_{ss'}$  est l'ordre de $ss'$ pour  $s,s'\in S$, alors $W$ admet comme
pr\'esentation

$$\genby{s\in S\mid s^2=1, \underbrace{ss'\cdots}_{m_{ss'}}=
\underbrace{s's\cdots}_{m_{ss'}}}$$

Soit   $\bW$  un  ensemble  muni  d'une  bijection  $W\iso\bW,w\mapsto \bw$
et soit $l$ la longueur sur $W$ d\'efinie par l'ensemble
de  g\'en\'erateurs  $S$.  Nous  d\'efinissons  le  {\it mono{\"\i}de de
tresses} d'Artin-Tits $B^+$ associ\'e \`a $W$ par la pr\'esentation
$$\genby{\bw\in\bW \mid \bw_1\bw_2=\bw'\text{ lorsque }w_1w_2=w'
\text{ et }l(w_1)+l(w_2)=l(w')}.$$
L'application $\bw\mapsto w$ induit un morphisme de mono\"{\i}des
$\beta:B^+\to W$.
Notons  $\bS$  le  sous-ensemble  de  $\bW$  en  bijection  avec $S$. La
pr\'esentation   pr\'ec\'edente  est  \'equivalente   (\cf\  \eg\ 
\cite[proposition 1.1]{michel}) \`a la pr\'esentation \og classique\fg
$$\genby {\bs\in\bS\mid \underbrace{\bs\bs'\cdots}_{m_{ss'}}=
\underbrace{\bs'\bs\cdots}_{m_{ss'}}}.$$
Nous notons aussi $l$ la longueur sur
$B^+$ d\'efinie par $\bS$. Alors, $\bW=\{\bw\in B^+ | l(\bw)=l(\beta(\bw))\}$.

On note $B$ le groupe de m\^eme pr\'esentation que $B^+$: c'est le
{\it groupe de tresses} d'Artin-Tits. L'application canonique
$B^+\to B$ est injective (\cf\ \eg\  \cite[corollaire 3.2]{michel}).
L'identit\'e de $\bS$ s'\'etend uniquement en un
anti-automorphisme de $B^+$, appel\'e {\em retournement}.

Nous notons \og $\le$\fg\ l'ordre
de  Bruhat sur  $W$. Pour  $w\in W$,  nous appelons {\em support} de  $w$
l'ensemble $\{s\in S\mid s\le w\}$.

Pour $I\subset  S$, nous notons  $W_I$ le sous-groupe de  $W$ engendr\'e
par  $I$. On  a $w\in  W_I$ o\`u  $I$  est le  support de  $w$. On  note
$w_0^I$  l'\'el\'ement de  plus  grande  longueur de  $W_I$  et on  pose
$w_0=w_0^S$. Nous dirons que $w\in W$ est {\em $I$-r\'eduit} (resp. {\em
r\'eduit-$I$}) si pour  tout $s\in I$ on a $sw>w$  (resp. $ws>w$). Notons
que $w\in W$ est $I$-r\'eduit (resp.
r\'eduit-$I$) si et seulement si pour  tout $v\in  W_I$,
on  a $l(v)+l(w)=l(vw)$ (resp. $l(w)+l(v)=l(wv)$).

Soit  $\bpi=\bw_0^2$; c'est un  \'el\'ement  central  de $B$ (lorsque $W$ est
irr\'eductible,  c'est le g\'en\'erateur positif du centre du groupe des
tresses pures, noyau de $\beta:B\to W$).

\sub{}
Il nous sera  commode de travailler avec une version  enrichie de $B^+$.

Soit  $\uW$  un  ensemble  muni d'une  bijection  $W\iso  \uW$,
$w\mapsto \uw$. Le {\it mono{\"\i}de  de tresses compl\'et\'e} $\uB^+$ a
pour ensemble de  g\'en\'erateurs\linebreak $\{y(w),y(\uw')\}_{w\in W,\uw'\in\uW}$
et pour relations

$\bullet$ $y(1)=y(\underline{1})=1$,

$\bullet$ $y(w_1)y(w_2)=y(w_1w_2)$ pour $l(w_1w_2)=l(w_1)+l(w_2)$,

$\bullet$ $y(\uw)y(\uw')=y(\underline{ww'})$ lorsque $w$ et $w'$ ont des
supports disjoints,

$\bullet$ $y(w)y(\uw')=y(\uw')y(w)$ lorsque $wv=vw$ et $l(wv)=l(w)+l(v)$
pour tout $v\le w'$.

On  a  une  injection  $B^+\to  \uB^+$  induite  par  $\bw\mapsto  y(w)$
qui  nous  permet  d'identifier   $B^+$  avec  le  sous-mono{\"\i}de  de
$\uB^+$  engendr\'e par  les $y(w)$.  On note  $\ubw=y(\uw)$ et  on pose
$\ubW=\{\ubw\}_{w\in  W}$. On  munit  $\uB^+$  d'une fonction  longueur:
c'est  le  morphisme  de   mono{\"\i}des  $l:\uB^+\to\bbN$  donn\'e  par
$l(\bw)=l(\ubw)=l(w)$.

On dispose d'un morphisme $\rho:\uB^+\to B^+$ donn\'e par
$y(w)\mapsto \bw$ et $y(\uw)\mapsto \bw$. C'est une section
de l'injection canonique.

\begin{remarque}
Une propri\'et\'e remarquable de $\uB^+$ est l'existence d'un
morphisme de mono{\"\i}des $\uB^+\to \bbZ B^+$ donn\'e par
$y(w)\mapsto \bw$ et $y(\uw)\mapsto \sum_{v\le w}\bv$. Si ce morphisme
s'av\'erait ne pas \^etre injectif, il y aurait lieu d'\'etudier si les 
relations suppl\'ementaires donnent
lieu \`a des relations entre vari\'et\'es correspondantes (\cf\ \S
\ref{varietesuB}).
\end{remarque}
\sub{}
Nous rappelons maintenant la construction et les propri\'et\'es des
formes normales des  \'el\'ements de $B^+$.

Nous  notons $\preccurlyeq$  la  divisibilit\'e  \`a gauche  dans  $B^+$,
\ie,  $\bx\preccurlyeq\by$ s'il existe $\bz\in  B^+$ tel que $\bx\bz=\by$.
Le mono\"\i de $B^+$ \'etant simplifiable, la relation $\preccurlyeq$ est un
ordre partiel.
Nous  notons  $B^+_I$  (resp.\ $B_I$)  le
sous-mono{\"\i}de de $B^+$ (resp.\ le
sous-groupe de $B$)  engendr\'e  par $\bI$, o\`u
$I\subset  S$  ($B^+_I$ s'identifie  au  mono{\"\i}de  des
tresses du sous-groupe $W_I$ de $W$ engendr\'e par $I$).
De m\^eme, nous notons $\uB_I^+$   le
sous-mono{\"\i}de de $\uB^+$
engendr\'e  par $\bI\cup \ubW_I$.

\begin{lemme-definition}\label{alpha}  Soit $\bw\in B^+$. Alors
\begin{enumerate}
\item  Il  existe  un  unique  diviseur  \`a  gauche  maximal  de  $\bw$
dans  $\bW$.   Nous  le  notons  $\alpha(\bw)$,   et  nous  posons
$\omega(\bw)=\alpha(\bw)\inv\bw$.
\item  Il  existe  un  unique  diviseur  \`a  gauche  maximal  de  $\bw$
dans  $B^+_I$. Nous  le  notons $\alpha_I(\bw)$  et nous  posons
$\omega_I(\bw)=\alpha_I(\bw)\inv\bw$.
\end{enumerate}
\end{lemme-definition}
\begin{pf*}{Preuve} (i) est, par exemple, \cite[proposition 2.1]{michel}.
(ii) r\'esulte  imm\'e\-dia\-tement de  \cite[lemme 1.4]{michel}  car l'ensemble
des $\bx\in B^+_I$ tel que $\bx\preccurlyeq\bw$ v\'erifie les hypoth\`eses de
{\it loc. cit.}
\end{pf*}
Nous dirons  qu'une d\'ecomposition $\bw=\bw_1\cdots \bw_k$  est la {\it
forme  normale} de  $\bw$  si  pour tout  $i$  on a  $\bw_i=\alpha(\bw_i
\bw_{i+1}\cdots  \bw_k)$. Le  mono{\"\i}de  $B^+$ \'etant  simplifiable,
tout  $\bw$  poss\`ede  une  unique  forme normale,  qui  en  donne  une
d\'ecomposition canonique en produit d'\'el\'ements de $\bW$.

Le r\'esultat suivant est classique (\cf\ par exemple \cite[Corollaire 4.4]{michel}).
\begin{proposition}
\label{monoidepointsfixes}
Soit $\sigma$ un automorphisme du syst\`eme de Coxeter $(W,S)$.
\begin{enumerate}
\item $(W^\sigma,\{t_I\}_{I\in S/\sigma})$
est un syst\`eme de Coxeter, o\`u $t_I=w_0^I$ et
$I$ d\'ecrit l'ensemble des orbites de $\sigma$ dans $S$.
De plus, les longueurs dans ce syst\`eme de Coxeter s'ajoutent si et seulement
si elles s'ajoutent dans $W$.
\item
On a un isomorphisme
$B^+_{W^\sigma}\iso C_{B^+}(\sigma),\ \bt_I\mapsto \bw_0^I$.
\end{enumerate}
\end{proposition}

\subsection{Vari\'et\'es \og classiques\fg}
\sub{}
Nous rappelons ici quelques r\'esultats g\'eom\'etriques qui nous seront utiles par
la suite.

Soit $k$ un corps alg\'ebriquement clos. On appelle vari\'et\'e sur $k$
un sch\'ema quasi-projectif sur $k$.

Soit $\{Z_i \}$ une famille finie de sous-vari\'et\'es localement ferm\'ees
d'une vari\'et\'e $X$. Nous notons $\bigcup Z_i$ l'union de ces
sous-vari\'et\'es --- ce sera une sous-vari\'et\'e localement ferm\'ee lorsque nous
utiliserons cette notation. Nous notons aussi
$\coprod Z_i$ cette union lorsque $Z_i\cap Z_j=\emptyset$ pour $i\not=j$.

On utilisera le r\'esultat suivant, pour montrer que certains morphismes
de vari\'et\'es sont des isomorphismes \cite[proposition II.6.6 et corollaire
de 6.1]{Borel}.

\begin{theoreme}
\label{isonormale}
Soit $f:X\to Y$ un morphisme bijectif. On suppose que les composantes
irr\'eductibles de $X$ sont ses composantes connexes et que $Y$ est
normal.
Si $f$ est s\'eparable (par exemple birationnel), alors c'est un
isomorphisme.
\end{theoreme}

\begin{proposition}
\label{fibration}
Soit $f:X\to Y$ une fibration localement triviale de fibre normale
(resp. lisse).
Alors, $X$ est normale (resp. lisse) si et seulement si $Y$ est normale
(resp. lisse).
\end{proposition}

\begin{pf*}{Preuve}
En effet, localement pour la topologie de Zariski, on a un produit
cart\'esien. Dans le cas d'un produit cart\'esien, l'assertion sur la lissit\'e
est \'etablie dans \cite[VI, III th\'eor\`eme 2]{Chevalley} et celle sur la
normalit\'e dans \cite[V, I proposition 3]{Chevalley}.
\end{pf*}

Nous aurons aussi besoin de comparer des propri\'et\'es \`a travers un
rev\^etement \'etale \cite[I, remarque 2.24 et proposition 3.17]{Milne}.

\begin{proposition}
\label{revetementetale}
Soit $f:X\to Y$ un rev\^etement \'etale fini. Alors, $X$ est lisse (resp.
normale) si et seulement si $Y$ est lisse (resp. normale).
\end{proposition}

\begin{definition}\label{formule des traces}
Nous dirons qu'un endomorphisme fini $\phi$ d'une vari\'et\'e $\bX$ 
\og v\'erifie la formule des traces\fg\ si
$\TrH{\phi}{\bX}=|\bX^\phi|$.
\end{definition}

Rappelons la formule des traces de Lefschetz \og classique\fg\
\cite[Theorem 12.3]{Milne}.

\begin{theoreme}
\label{Lefschetz}
Si $\bX$ est une vari\'et\'e projective
lisse et si le graphe de $\phi$ est transverse \`a la diagonale,
alors, $\phi$ v\'erifie la formule des traces.
\end{theoreme}

Le th\'eor\`eme suivant (conjecture de Deligne) est d\^u \`a Fujiwara
\cite{Fu}:
\begin{theoreme}\label{fuji}
Soit  $\bX$ une  vari\'et\'e sur la cl\^oture alg\'ebrique d'un corps fini,
munie d'un  endomorphisme de Frobenius $F$ et  soit $\phi$ un
endomorphisme fini de $\bX$. Alors pour $n$ entier suffisamment grand et
tel que $F^n$  commute \`a $\phi$, l'endomorphisme  $\phi F^n$ v\'erifie
la formule des traces. En particulier, pour un tel $n$ on a
$$\TrH\phi\bX=-\lim_{t\to\infty}\sum_{k=1}^\infty|\bX^{\phi F^{nk}}|t^k.$$
\end{theoreme}

\sub{}
Dans cet  article nous utiliserons  les notations suivantes: $p$  est un
nombre premier, $\bar{\bbF}_p$ est une  cl\^oture
alg\'ebrique du corps fini $\bbF_p$, 
$\bG$ est un groupe r\'eductif connexe sur $\bar{\bbF}_p$
muni d'une isog\'enie $F$ dont une
puissance est  l'endomorphisme de Frobenius attach\'e  \`a une structure
rationnelle de  $\bG$ sur  un corps  fini. Nous  notons $\GF$  le groupe
(fini) des  points fixes de  $\bG$ sous  $F$. Nous notons  $\delta$ le
plus petit entier tel que  $F^\delta$ soit un endomorphisme de Frobenius
munissant $\bG$  d'une structure  rationnelle d\'eploy\'ee sur  le sous-corps
$\bbF_{q^\delta}$ \`a $q^\delta$ \'el\'ements de $\bar{\bbF}_p$,
o\`u  $q$ est un nombre  r\'eel positif d\'efini
par cette condition ($q^\delta$ est une puissance enti\`ere de $p$).

Nous notons $\CB$ la vari\'et\'e des sous-groupes de Borel de $\bG$ et
$W$ le groupe de Weyl de $\bG$: ses \'el\'ements correspondent
aux orbites
de $\bG$ dans son action diagonale sur $\CB\times\CB$. Muni
de l'ensemble $S$ de g\'en\'erateurs correspondant aux
orbites de dimension $1$, $W$ est un groupe de Coxeter.
Pour $w\in W$, nous notons $\CO(w)$ l'orbite correspondante.
On a un isomorphisme canonique
$\CO(w)\times_\CB \CO(w')\iso \CO(ww')$ donn\'e par la premi\`ere
et la derni\`ere projection lorsque $l(ww')=l(w)+l(w')$.

On dit que $\bB_1$ est en {\it position relative} $w$ avec
$\bB_2$ (et on note $\bB_1\xrightarrow w\bB_2$) lorsque $(\bB_1,\bB_2)\in\CO(w)$.

L'action  de $F$  sur  $W$ stabilise  $S$. L'application  correspondante
$\bS\to\bS$ s'\'etend de mani\`ere unique  en un automorphisme $B\to B$.
On a  obtenu ainsi une  action de $F$  sur $B$ compatible  (via $\beta$)
avec  l'action  sur $W$.  On  d\'efinit  de  m\^eme  une action  de  $F$
sur  $\uB^+$ par  $F(y(w))=y(F(w))$ et  $F(y(\uw))=y(\underline{F(w)})$.
L'ordre de cette action est $\delta$.

On  fixe dans  toute  la  suite un  couple  $F$-stable $\bT\subset  \bB$
form\'e  d'un  tore  maximal  de  $\bG$ et  d'un  sous-groupe  de  Borel
le  contenant.  La  vari\'et\'e  $\CB$  s'identifie  \`a  $\bG/\bB$  par
l'isomorphisme  $\lexp  g\bB\mapsto g\bB$.  On  identifie  le groupe  de
Weyl  \`a   $N_\bG(\bT)/\bT$,  en  associant  l'orbite   de  $(\bB,\lexp
w\bB)$  \`a $w\in  N_\bG(\bT)/\bT$  ;  avec l'identification  ci-dessus,
on a 
$$\CO(w)\iso \{(g_1\bB,g_2\bB)\in\bG/\bB\times\bG/\bB\mid g_1\inv
g_2\in\bB w\bB\}.$$
 La g\'eom\'etrie des  orbites $\CO(w)$ est li\'ee \`a
celle des {\it cellules de Bruhat} $\CB(w)=\{\bB'\in\CB\mid \bB\xrightarrow
w\bB'\}$, qui sont des espaces affines de dimension $l(w)$.
L'adh\'erence de  $\CB(w)$ dans $\CB$ est  la {\it vari\'et\'e
de Schubert} $\overline{\CB(w)}=\coprod_{w'\le w}\CB(w')$.

\begin{lemme}\label{Owbar}
Soit $\overline{\CO(w)}$ l'adh\'erence de $\CO(w)$ dans $\CB\times\CB$.
\begin{enumerate}
\item On a $\overline{\CO(w)}=\coprod_{w'\le w}\CO(w')$.
\item La  vari\'et\'e $\overline{\CO(w)}$ est  lisse si et  seulement si
$\overline{\CB(w)}$ est lisse.
\end{enumerate}
\end{lemme}
\begin{pf*}{Preuve}
Consid\'erons le diagramme:
$$\xymatrix{
\bG\times\bG\ar[r]^q\ar[d]_p & \CB\times\CB\\
\bG\ar[r]_r & \CB
}$$
o\`u  les fl\`eches sont donn\'ees par
$q(g_1,g_2)=(\lexp{g_1}\bB,\lexp{g_2}\bB)$, $p(g_1,g_2)=g_1\inv g_2$
et   $r(g)=\lexp   g\bB$.   Les   morphismes   $p,q,r$   sont 
ouverts.
On  a $r\inv(\CB(w))=\bB  w\bB$.
On  en  d\'eduit  que  l'adh\'erence  de  $\bB  w\bB$  est  \'egale  \`a
$r\inv(\overline{\CB(w)})=  \coprod_{w'\le   w}\bB  w'\bB$,   car  cette
derni\`ere  vari\'et\'e est  une  union  de fibres  de  $r$. De  m\^eme,
l'adh\'erence dans $\bG\times\bG$ de
$p\inv(\bB w \bB)=\{(g_1,g_2)\in\bG\times\bG\mid g_1\inv g_2\in \bB w\bB\}$
est $\{(g_1,g_2)\in\bG\times\bG\mid
g_1\inv g_2\in \coprod_{w'\le w}\bB w'\bB\}$; enfin ces deux derni\`eres
vari\'et\'es \'etant  unions de fibres de  $q$ et ayant pour  images par
$q$ respectivement $\CO(w)$ et $\coprod_{w'\le w}\CO(w')$, on en d\'eduit
le (i) du lemme.

Le  (ii) se  d\'eduit (\cf\ proposition \ref{fibration})
du  fait que  dans le  diagramme ci-dessus  $q$ et
$r$  sont des  fibrations localement  triviales de  fibre $\bB\times\bB$
et  $\bB$  respectivement et  que  $p$  est  une fibration  triviale  de
fibre $\bG$. 
\end{pf*}

\begin{corollaire}
\label{O(W_I)lisse}
Soit $I\subset  S$. Alors,
les vari\'et\'es $\overline{\CB(w_0^I)}$ et $\overline{\CO(w_0^I)}$ sont lisses.
\end{corollaire}

\begin{pf*}{Preuve}
On a $\overline{\CB(w_0^I)}\simeq \bB W_I\bB/\bB$ o\`u $\bB W_I\bB$ est
un sous-groupe parabolique de $\bG$, donc $\overline{\CB(w_0^I)}$
est une vari\'et\'e lisse. Par cons\'equent, $\overline{\CO(w_0^I)}$ est
lisse d'apr\`es le lemme \ref{Owbar}.
\end{pf*}

\sub{}
Pour $w\in W$, nous posons
$\CO(\uw)=\overline{\CO(w)}$,
$\CB(\uw)=\overline{\CB(w)}$,
$\bB\uw\bB=\overline{\bB w\bB}$ et nous \'ecrirons
$\bB_1\xrightarrow\uw\bB_2$ la propri\'et\'e \og $(\bB_1,\bB_2)\in\CO(\uw)$\fg.

Avec ces notations, nous avons la

\begin{definition}\label{O(t)}
Pour $t_1,\ldots,t_k\in W\cup\uW$,
nous posons
\begin{align*}\CO(t_1,\ldots,t_k)&=
\CO(t_1)\times_\CB\CO(t_2)\times_\CB\cdots\times_\CB\CO(t_k)\\
&=\{(\bB_1,\ldots,\bB_{k+1})\in\CB^{k+1}\mid\bB_i\xrightarrow{t_i}\bB_{i+1}\}\\
&\iso\{(g_1\bB,\ldots,g_{k+1}\bB)\in(\bG/\bB)^{k+1}\mid \hfill
g_i\inv g_{i+1}\in\bB t_i\bB\}.
\end{align*}
Nous notons $p':\CO(t_1,\ldots,t_k)\to\CB,\ 
(\bB_1,\ldots,\bB_{k+1})\mapsto\bB_1$ la premi\`ere  projection
et $p'':(\bB_1,\ldots,\bB_{k+1})\mapsto\bB_{k+1}$ la derni\`ere projection.
\end{definition}

\begin{lemme}\label{lissiteO}
Soient $t_1,\ldots,t_k\in W\cup\uW$.
Alors,  la  vari\'et\'e  $\CO(t_1,\ldots,t_k)$  est  
normale. Si  la
vari\'et\'e $\CO(t_i)$  est lisse pour  tout $i$, alors  la vari\'et\'e
$\CO(t_1,\ldots,t_k)$ est lisse.
\end{lemme}
\begin{pf*}{Preuve}
Notons  $q'_k$  (resp.  $q_k''$)  l'oubli  du  premier  (resp.  du  dernier)
sous-groupe   de   Borel;   nous   allons   d\'emontrer   par   r\'ecurrence
sur   $k$    que   ces    applications   d\'efinissent    des   fibrations
localement triviales $\CO(t_1,\ldots,t_k)\to \CO(t_2,\ldots,t_k)$ (resp.
$\CO(t_1,\ldots,t_k)\to \CO(t_1,\ldots,t_{k-1})$),  de fibres isomorphes
\`a $\CB(t_1)$ (resp. $\CB(t_k)$).

D\'emontrons tout  d'abord le r\'esultat  pour $k=1$.
On a un carr\'e cart\'esien
$$\xymatrix{
Z\ar[rr]^s \ar[d]_p && \CO(t_1) \ar[d]^{p'=q_1''}\\
\bG \ar[rr]_r && \CB
}$$
o\`u 
$Z=\{(g,\bB')\in \bG\times\CB | (\lexp g\bB,\bB')\in \CO(t_1)\}$,
$r(g)=\lexp g\bB$, $s(g,\bB')=(\lexp g\bB,\bB')$ et $p$ est la
premi\`ere projection. On a un isomorphisme
$\bG\times\CB(t_1)\iso Z,\ (g,\bB')\mapsto (g,\lexp g\bB')$, donc
$p$ est  une fibration triviale de  fibre $\CB(t_1)$. Puisque $r$  est une
fibration localement  triviale de fibre  $\bB$, on en d\'eduit  que $p'$
est une fibration localement triviale de fibre $\CB(t_1)$ (localement, on
a un produit direct).

Le cas de la seconde projection $p''=q'_1$ se traite de mani\`ere
sym\'etrique.

Nous traitons maintenant le cas g\'en\'eral par r\'ecurrence sur $k$.
On consid\`ere le carr\'e cart\'esien
$$\xymatrix{
\CO(t_1,\ldots,t_k) \ar[rr]^{q'_k}\ar[d]_{q_k''}
 &&\CO(t_2,\ldots,t_k)\ar[d]^{q_{k-1}''}\\
\CO(t_1,\ldots,t_{k-1}) \ar[rr]_{q'_{k-1}} &&\CO(t_2,\ldots,t_{k-1})
}$$
Par r\'ecurrence, $q'_{k-1}$ est une fibration localement triviale de fibre
$\CB(t_1)$ et
$q_{k-1}''$ est une fibration localement triviale de fibre $\CB(t_k)$.
Par cons\'equent, $q_k''$ et $q'_k$ sont des
fibrations localement triviales de fibres $\CB(t_k)$ et $\CB(t_1)$
respectivement
(localement on a des produits cart\'esiens au-dessus d'un ouvert de
$\CO(t_2,\ldots,t_{k-1})$).

On  d\'eduit alors par r\'ecurrence que $\CO(t_1,\ldots,t_k)$ est normale
(proposition \ref{fibration}) car
les  vari\'et\'es de  Schubert $\CB(t_i)$ sont  normales 
\cite[th\'eor\`eme 3]{Rama}.

Si les vari\'et\'es $\CO(t_i)$ sont lisses, les vari\'et\'es $\CB(t_i)$ le
sont aussi d'apr\`es le lemme \ref{Owbar} et on d\'eduit
par r\'ecurrence que $\CO(t_1,\ldots,t_k)$ est lisse
(proposition \ref{fibration}).
%
\end{pf*}

\sub{}
\label{varietesuB}
Nous allons maintenant g\'en\'eraliser l'id\'ee de \cite{Sydney,Deligne}
consistant
\`a associer une vari\'et\'e \`a un \'el\'ement de $B^+$ au cas
d'\'el\'ements de $\uB^+$. Nous allons voir que la construction de
\cite{Deligne} se recolle convenablement.

Pour $t_1,\ldots,t_k\in W\cup\uW$ et
$t'_1,\ldots,t'_k\in W$, nous
\'ecrirons $(t'_1,\ldots,t'_k)\subset(t_1,\ldots,t_k)$ si pour tout $i$,
on a $t'_i=t_i$ ou $t_i=\uw\in \uW$ et $t'_i\le w$.
On a une d\'ecomposition en union de sous-vari\'et\'es localement ferm\'ees
\begin{equation}
\label{lt ferme}
\CO(t_1,\ldots,t_k)=\coprod_{(t'_1,\ldots,t'_k)\subset(t_1,\ldots,t_k)}
\CO(t'_1,\ldots,t'_k).
\end{equation}

\begin{proposition}\label{omnibus} Soient $w,w'\in W$.
\begin{enumerate}
\item On a une immersion ouverte canonique $\CO(w)\to\CO(\uw)$.
\item Si $w\le w'$, alors on a une immersion ferm\'ee canonique
$\CO(\uw)\to\CO(\uw')$.
\item  Si  $l(ww')=l(w)+l(w')$,  alors $(p',p'')$  induit  un  morphisme
canonique $\CO(\uw,\uw')\to\CO(\underline{ww'})$ qui se restreint
en un isomorphisme canonique \allowbreak$\CO(w,w')\iso\CO(ww')$.
\item  Si  $w$   et  $w'$  ont  des  supports   disjoints, alors  le  morphisme
canonique
$$(p',p''):\CO(\uw,\uw')\iso\CO(\underline{ww'})$$ est  un
isomorphisme.
\item Si $l(wv)=l(w)+l(v)$ pour tout $v\le w'$, alors
$(p',p'')$ induit un isomorphisme canonique
$\CO(w,\uw')\iso \coprod_{v\le w'}\CO(wv)$. De m\^eme, 
si $l(vw')=l(v)+l(w')$ pour tout $v\le w$, alors
$(p',p'')$ induit un isomorphisme canonique
$\CO(\uw,w')\iso \coprod_{v\le w}\CO(vw')$.
\item  Si $wv=vw\ge w$ pour  tout $v\le w'$, alors  on a un isomorphisme
canonique    $$\CO(w,\uw')\iso\CO(\uw',w)$$   caract\'eris\'e   par   la
commutativit\'e du diagramme
$$\xymatrix{
\CO(w,\uw')\ar[rr]^\sim\ar[dr]_{(p',p'')} && \CO(\uw',w)\ar[dl]^{(p',p'')} \\
& \CO(\uw\uw')
}.$$
\end{enumerate}
\end{proposition}
\begin{pf*}{Preuve}
(i),  (ii),  (iii)  sont  cons\'equences imm\'ediates des d\'efinitions
ou  sont  des propri\'et\'es classiques de l'ordre de  Bruhat.

D\'emontrons (iv).  La caract\'erisation  de l'ordre  de Bruhat  par les
suites extraites  d'une d\'ecomposition  r\'eduite montre  que $\{v\mid
v\le  ww'\}=\{v_1v_2\mid (v_1,v_2)\subset  (\uw,\uw')\}$. Comme  de tels
$v_1,v_2$ v\'erifient  les hypoth\`eses  de (iii),  le morphisme  est un
isomorphisme entre les termes de la d\'ecomposition (\ref{lt ferme}), il
est donc bijectif. Il se restreint  en un isomorphisme entre les ouverts
denses  $\CO(w,w')\iso\CO(ww')$, donc  il est  birationnel. Puisque  les
vari\'et\'es consid\'er\'ees sont normales (lemme \ref{lissiteO}), c'est
un isomorphisme (th\'eor\`eme \ref{isonormale}).

Pour d\'emontrer (v), nous commen\c cons par un lemme.

\begin{lemme}\label{sousmots}
Soit  $I$  un  sous-ensemble  de  $S$,  soit  $w\in  W$  un  \'el\'ement
r\'eduit-$I$, soit $w'\in  W_I$ et soit $v\le ww'$.

Si $v\ge w$, alors
$v=ww'_1$ o\`u $w'_1\le w'$. Sinon, $v=w_1 w'_1$  o\`u $w_1<w$,
$w'_1\le w'$ et $l(w_1w'_1)=l(w_1)+l(w'_1)$.
\end{lemme}

\begin{pf*}{Preuve}
Puisque  $v\le ww'$,  il  existe $v_1\le  w$ et  $v'_1\le  w'$ tels  que
$v=v_1v'_1$.  D'apr\`es le  lemme  d'\'echange, il  existe $w_1\le  v_1$
et  $w'_1\le v'_1$  tels  que $v=w_1w'_1$  et $l(v)=l(w_1)+l(w'_1)$. 

Si $v\not\ge w$,  alors $w_1\ne w$,  donc $w_1<w$,  d'o\`u la
conclusion de l'\'enonc\'e  dans ce cas.

Supposons  maintenant $v\ge w$.
Comme pr\'ec\'edemment, il existe $w_2\le w_1$ et  $w'_2\le w'_1$  tels que
$w=w_2w'_2$  et $l(w)=l(w_2)+l(w'_2)$.
Puisque $w'_2\in W_I$ et $w$ est r\'eduit-$I$, on a $w'_2=1$,
donc $w\le w_1$ et finalement $w=w_1$.
\end{pf*}

Supposons $l(wv)=l(w)+l(v)$ pour tout $v\le w'$. Alors, $w$ est
r\'eduit-$I$, pour $I$ le support de $w'$.
Il   r\'esulte   du   lemme  \ref{sousmots} que le compl\'ementaire   de
$\coprod_{w'_1\le w'}\CO(ww'_1)$ dans $\CO(\underline{ww'})$ est
l'union   des  $\CO(w_1  w'_1)$   o\`u  $w_1<w$,  $w'_1\le   w'$  et
$l(w_1w'_1)=l(w_1)+l(w'_1)$. Comme les  longueurs s'ajoutent, si $w_2\le
w_1$ et $w'_2\le  w'_1$ alors $w_2w'_2\le w_1w'_1$ donc  cette union est
une sous-vari\'et\'e ferm\'ee.
Par cons\'equent, $\coprod_{w'_1\le w'}\CO(ww'_1)$ est une sous-vari\'et\'e
ouverte  de  la  vari\'et\'e  normale  $\CO(\uw\uw')$,
donc  c'est  une vari\'et\'e  normale. On conclut
comme pr\'ec\'edemment que
$(p',p''):\CO(w,\uw')\to\coprod_{w'_1\le w'}\CO(ww'_1)$ est un isomorphisme.

La seconde partie de (v) se d\'emontre de la m\^eme mani\`ere.

L'hypoth\`ese  de (vi) montre que  $w$ est r\'eduit-$I$ et $I$-r\'eduit,
o\`u   $I$  est  le   support  de  $w'$.   L'assertion  (vi)  est  alors
cons\'equence imm\'ediate de (v).
\end{pf*}

La proposition \ref{omnibus} fournit, par produits fibr\'es, des
morphismes canoniques pour les vari\'et\'es associ\'ees \`a des suites.
Par exemple, (vi) fournit un isomorphisme canonique
$$\CO(t_1,\ldots,t_r,w,\uw',t_{r+1},\ldots,t_k)\iso
\CO(t_1,\ldots,t_r,\uw',w,t_{r+1},\ldots,t_k)$$
pour $t_1,\ldots,t_k\in W\cup\uW$.
Les  isomorphismes   donn\'es  par  les   (iii),  (iv)  et  (vi)   de  la
proposition \ref{omnibus}  sont associ\'es  \`a chacun des  trois types
de  relations  de $\uB^+$:  en  composant  ces isomorphismes  on  dispose
donc d'isomorphismes entre  deux vari\'et\'es $\CO(t_1,\ldots,t_k)$ et
$\CO(t'_1,\ldots,t'_k)$  lorsque $\bt=\bt_1\cdots\bt_k$  et
$\bt'=\bt'_1\cdots\bt'_{k'}$ sont \'egaux dans $\uB^+$.

\begin{proposition-definition}\label{transitif}
Soit $\bt\in \uB^+$.
Le syst\`eme d'isomorphismes entre vari\'et\'es $\CO(t'_1,\ldots,t'_k)$ telles
que $\bt'_1\cdots\bt'_k=\bt$ est transitif et on note $\CO(\bt)$ la limite
projective du syst\`eme. Elle est munie d'une application (premi\`ere
et derni\`ere projections, que nous notons $(p',p'')$) vers $\CB\times\CB$.
\end{proposition-definition}

\begin{pf*}{Preuve}
Lorsque $\bt\in B^+$, le r\'esultat est d\^u \`a Deligne
\cite[application 2]{Deligne}.

Dans  le  cas  g\'en\'eral,  les  items  (iv)  et (vi) de la proposition
\ref{omnibus}   montrent  que   les  isomorphismes   entre  vari\'et\'es
$\CO(t_1,\ldots,t_k)$  correspondant  \`a  une  relation  de  $\uB^+$ sont
d\'etermin\'es  par  des  isomorphismes  correspondants  donn\'es par la
proposition  \ref{omnibus} (iii) entre les  termes de la d\'ecomposition
(\ref{lt  ferme}) de ces  vari\'et\'es, c'est-\`a-dire des isomorphismes
dont  la  transitivit\'e  a  \'et\'e  d\'emontr\'ee  par  Deligne. On en
d\'eduit   donc  la  transitivit\'e   du  syst\`eme  d'isomorphismes  en
g\'en\'eral.
\end{pf*}

Chaque d\'ecomposition $\bt=\bt_1\cdots\bt_k$ avec $\bt_i\in \bW\cup \ubW$
fournit un isomorphisme canonique
$\CO(\bt)\iso \CO(t_1,\ldots,t_k)$.

L'action de $F$ sur $\CB$ induit un
morphisme $$\CO(t_1,\ldots,t_k)\to \CO(F(t_1),\ldots,F(t_k)).$$
Ces morphismes sont compatibles aux isomorphismes canoniques de la
proposition \ref{omnibus}. Par cons\'equent, $F$ induit un morphisme
$\CO(\bt)\to\CO(F(\bt))$, pour tout $\bt\in \uB^+$.

\subsection{Vari\'et\'es de Deligne-Lusztig g\'en\'eralis\'ees}
\sub{}
Nous allons maintenant introduire les vari\'et\'es de Deligne-Lusztig
attach\'ees aux \'el\'ements du mono{\"\i}de de tresses
et des compactifications partielles de ces vari\'et\'es
associ\'ees aux \'el\'ements du mono{\"\i}de de tresses compl\'et\'e.
Ces vari\'et\'es sont l'objet principal d'\'etude de cet article.

\begin{definition}\label{DL}
Soit $\bt\in\uB^+$  et soit $\Gamma\subseteq \CB\times\CB$  le graphe de
$F$.  Nous appelons  \og vari\'et\'e de  Deligne-Lusztig\fg\ attach\'ee  \`a
$\bt$  la  vari\'et\'e  
$$\bX(\bt)=\CO(\bt)\times_{(\CB\times\CB)}\Gamma=
\{x\in\CO(\bt)\mid p''(x)=F(p'(x))\}.$$
 De  m\^eme,
pour $t_1,\ldots,t_k\in W\cup\uW$,
nous     posons     
$$\bX(t_1,\ldots,t_k)=     \{(\bB_1,\ldots,\bB_k)\mid
(\bB_1,\ldots,\bB_k,F(\bB_1))\in\CO(t_1,\ldots,t_k)\}.$$
Nous  noterons
parfois  $\bX_\bG(\bt)$ (resp.  $\bX(\bt,F)$, resp.  $\bX_\bG(\bt,F)$)
pour  pr\'eciser  le groupe  (resp.  l'isog\'enie,  resp. le  groupe  et
l'isog\'enie) utilis\'e dans la d\'efinition de la vari\'et\'e.
\end{definition}

Les  vari\'et\'es  de  Deligne-Lusztig   sont  munies  d'une  action  de
$\GF$: l'action diagonale de $\bG$ sur $\CB^{k+1}$
se  restreint en  une action  de $\bG^F$  sur $\bX(t_1,\ldots,t_k)$.

L'endomorphisme $F$ induit un morphisme
$\bX(\bt)\to\bX(F(\bt))$ pour tout $\bt\in\uB^+$. En particulier,
les vari\'et\'es de Deligne-Lusztig sont munies d'une action de $F^\delta$.
Si $\bt\in (\uB^+)^F$, alors $\bX(\bt)$ est munie d'une action de $F$.

Notons  que,  de m\^eme  que  pour  les vari\'et\'es  $\CO(\bt)$,  chaque
d\'ecomposition  $\bt=\bt_1\cdots\bt_k$  avec  $\bt_i\in  \bW\cup  \ubW$
fournit   un  isomorphisme  canonique
$$\bX(\bt)\iso\bX(t_1,\ldots,t_k).$$
Par composition, on obtient un isomorphisme canonique
$$\bX(t_1,\ldots,t_k)\iso \bX(t'_1,\ldots,t'_{k'})$$
pour $t'_1,\ldots,t'_{k'}\in W\cup  \uW$ tels que 
$\bt_1\cdots \bt_k=\bt'_1\cdots\bt'_{k'}$.

Pour $\bw\in\bW$,
la vari\'et\'e $\bX(\bw)$ est isomorphe \`a la vari\'et\'e de
Deligne-Lusztig  \og ordinaire\fg\  \cite[d\'efinition 1.4]{DL}
$$\bX(\bw)
\xrightarrow[p']{\sim}\bX(w)=\{\bB'\in\CB\mid\bB'\xrightarrow w F(\bB')\}
\iso\{g\bB\in \bG/\bB\mid g\inv F(g)\in \bB w\bB\}.$$

Plus g\'en\'eralement, nous pouvons identifier les vari\'et\'es de
Deligne-Lusztig associ\'ees \`a des \'el\'ements du mono{\"\i}de de tresses
\`a des vari\'et\'es de Deligne-Lusztig ordinaires pour un autre groupe
et une autre isog\'enie en appliquant la proposition suivante dans le
cas o\`u les $\bt_i$ sont dans $\bW$  (\cf\ \cite[\S 1.18]{LuCox}).

\begin{proposition}\label{produit}
Soit     $F_1$     l'isog\'enie     de     $\bG^k$     d\'efinie     par
$$F_1(g_1,\ldots,g_k)=(g_2,\ldots,g_k,F(g_1)).$$   Soit   $\bt=\bt_1\cdots
\bt_k\in  \uB^+$.  Alors,  on  a un  isomorphisme  entre  $\bX(\bt)$  et
la  vari\'et\'e  $\bX_{\bG^k}(\bt',F_1)$  attach\'ee  \`a  l'\'el\'ement
$\bt'=(\bt_1,\ldots,\bt_k)$  du  mono\"\i  de  de  tresses  compl\'et\'e
$(\uB^+)^k$  du  groupe  de  Weyl de  $\bG^k$.  L'action  de  $F^\delta$
correspond  par  cet isomorphisme  \`a  l'action  de $F_1^{k\delta}$  et
l'action de $\GF$ \`a celle de $(\bG^k)^{F_1}$.
\end{proposition}

\begin{pf*}{Preuve} Notons $\CO_{\bG^k}(\bt')$ la vari\'et\'e associ\'ee \`a
$\bt'$ dans le groupe $\bG^k$.
Soit $\bt'=\bt'_1\cdots\bt'_l$ une d\'ecomposition en \'el\'ements de
$\bW^k\cup\ubW^k$ avec
$\bt'_i=(\bt'_{i,1},\ldots,\bt'_{i,k})$ et $\bt'_{i,j}\in \bW\cup\ubW$.
On a un isomorphisme de la vari\'et\'e
\begin{multline*}\{(x_1,\ldots,x_{l+1})\in \CO_{\bG^k}(t'_1,\ldots,t'_l)\mid\\
x_1=(\bB_1,\bB_2,\ldots,\bB_k),x_{l+1}=(\bB_2,\bB_3,\ldots,\bB_k,\bB_{k+1})\};\\
\end{multline*}
vers $\CO_{\bG}(t'_{1,1},\ldots,t'_{l,k})$ donn\'e par
\begin{multline*}(x_1,\ldots,x_{l+1})\mapsto\\
(p_1(x_1),\ldots,p_1(x_{l+1})=
p_2(x_1),\ldots,p_2(x_{l+1}),\ldots,p_k(x_1),\ldots,p_k(x_{l+1}))\\
\end{multline*}
o\`u $p_i:\bG^k\to\bG$ est la $i$-\`eme projection.

Comme $\bt=\bt'_{1,1}\cdots\bt'_{l,k}$, on a un isomorphisme
canonique $$\CO(t'_{1,1},\ldots,t'_{l,k})\iso\CO(\bt).$$
La vari\'et\'e $X(\bt)$ \'etant la sous-vari\'et\'e de $\CO(\bt)$ d\'efinie par
$p''(x)=F(p'(x))$, on obtient l'isomorphisme annonc\'e.

Par cette identification,
l'action de $F_1^{k\delta}$ correspond bien \`a celle de $F^\delta$ et celle
de $(\bG^k)^{F_1}$ \`a celle de $\GF$.
\end{pf*}

\sub{}
La proposition suivante montre que pour $\bb\in B^+$,
les vari\'et\'es $\bX(\bb)$ sont
lisses et montre aussi que si $\bb=\bw_1\cdots\bw_k$ est une
d\'ecomposition en produit d'\'el\'ements de $\bW$ d'un \'el\'ement de
$\BW$ telle que $\CO(\uw_i)$ est lisse pour chaque $i$, alors
$\bX(\ubw_1\cdots\ubw_k)$ est une compactification lisse de $\bX(\bb)$.

\begin{proposition}
\label{desingularisation}
Soit $\bt\in  \uB^+$; alors, la  vari\'et\'e $\bX(\bt)$ est  normale. De
plus, si pour tout \'el\'ement $\ubw$  de $\ubW$ qui intervient dans une
d\'ecomposition de $\bt$, la vari\'et\'e $\CO(\uw)$ est lisse, alors, la
vari\'et\'e $\bX(\bt)$ est lisse.
\end{proposition}

\begin{pf*}{Preuve}
En  consid\'erant   les  \'el\'ements  $(1,\ldots,1,\bt_i,1,\ldots,1)\in
(\uB^+)^k$, la proposition \ref{produit} ram\`ene  la preuve au cas o\`u
on a $\bt=\bt_1\cdots\bt_k$ avec $t_i\in W\bigcup \uW$ tel que pour tout
$(t'_1,\ldots,t'_k)\subset(t_1,\ldots,t_k)$ avec  $t'_1,\ldots,t'_k\in W$,
on a $l(t'_1\cdots t'_k)=\sum_i l(t'_i)$.

On  note $\bB\bt\bB=\coprod_{(t'_1,\ldots,t'_k)}\bB  t'_1\cdots t'_k\bB$
et   $\CB(\bt)=\coprod_{(t'_1,\ldots,t'_k)}\CB(t'_1,\ldots,t'_k)$,  o\`u
$(t'_1,\ldots,t'_k)\subset(t_1,\ldots,t_k)$   et    $t'_1,\ldots,t'_k\in
W$.   On   suit   la   preuve   de   \cite[lemme   4.3]{Lusympl}.   Soit
$\CL:\bG\to\bG,g\mapsto  g^{-1}F(g)$   l'application  de   Lang:  c'est
l'application quotient  par $\bG^F$, pour son  action par multiplication
\`a  gauche.   C'est  donc   un  rev\^etement   \'etale  de   groupe  de
Galois  $\bG^F$.  Soit  $r:\bG\to\CB,  g\mapsto   \lexp  g  \bB$.  On  a
$r^{-1}(\CB(\bt))=\bB\bt\bB$  et $r^{-1}(\bX(\bt))=\CL^{-1}(\bB\bt\bB)$,
o\`u on a  identifi\'e $\bX(\bt)$ \`a une sous-vari\'et\'e  de $\CB$ via
$p'$.

On a un diagramme commutatif
$$\xymatrix{
\CL^{-1}(\bB\bt\bB) \ar[d]_r\ar[rr]^{\CL} && \bB\bt\bB
\ar[d]_r \\
\bX(\bt) && \CB(\bt)
}$$

Rappelons que $r$  est une fibration localement triviale  de fibre $\bB$
(\cf\ preuve du lemme \ref{Owbar}).  Par cons\'equent, la normalit\'e ou
la lissit\'e  de $\CB(\bt)$ est  \'equivalente \`a celle  de $\bB\bt\bB$
(proposition \ref{fibration}).  Il en est  de m\^eme pour  $\bX(\bt)$ et
$\CL^{-1}(\bB\bt\bB)$ (proposition \ref{revetementetale}).
La proposition r\'esulte alors de
\cite[th\'eor\`eme 3]{Rama} et du lemme \ref{Owbar}.
\end{pf*}

Notons que, puisque $\CO(\us)$ est lisse lorsque $s\in S$
(\cf\ corollaire \ref{O(W_I)lisse}), l'hypoth\`ese de la proposition
\ref{desingularisation} sera v\'erifi\'ee si les seuls \'el\'ements de $\ubW$
qui interviennent dans
$\bt$ sont de la forme $\ubs$ o\`u $s\in S$.

\begin{proposition}\label{quasiaffine}
Soit $\bt\in\uB^+$.
\begin{enumerate}
\item La vari\'et\'e $\bX(\bt)$ est de dimension $l(\bt)$ et
ses composantes connexes sont irr\'eductibles.
\item Si $\bt\in B^+$, alors $\bX(\bt)$ est quasi-affine.
\item  Si $\bt=\bt_1\cdots\bt_r$ avec $\bt_i\in\bW$ et $q\ge h$ o\`u $h$
est le nombre de Coxeter de $\bG$, alors $\bX(\bt)$ est affine.

\item Si $\bt$ est un produit d'\'el\'ements de $\ubW$, alors
$\bX(\bt)$ est projective.
\end{enumerate}
\end{proposition}

\begin{pf*}{Preuve}
Soit $\bb$ l'\'el\'ement du mono{\"\i}de de tresses produit des
\'el\'ements obtenus en rempla\c cant dans $\bt$ chaque \'el\'ement
$\ubw$ par $\bw$. Alors,
$\bX(\bb)$ est un ouvert dense de $\bX(\bt)$;
d'autre part, la dimension de $\CO(\bb)$ est $l(\bb)+\dim\CB$ et on
conclut que $\dim\bX(\bt)=l(\bt)$ en utilisant que l'intersection de $\Gamma$
(qui est de dimension $\dim\CB$) avec $\CO(\bb)$ est transverse.

La proposition \ref{produit} ram\`ene (ii) au cas d'une vari\'et\'e de
Deligne-Lusztig ordinaire et le r\'esultat est alors
\cite[th\'eor\`eme 2.3]{Haa}.
De m\^eme, (iii) se d\'eduit du cas des vari\'et\'es ordinaires \cite[th\'eor\`eme 9.7]{DL},
le nombre de Coxeter de $\bG^r$ \'etant le m\^eme que celui de $\bG$.

Le point (iv) est imm\'ediat, puisque $\CO(\bt)$ est projective dans
ce cas.
\end{pf*}

\sub{}
La proposition suivante (\cf\ \cite[Lemma 3]{LuFi})
d\'ecrit les vari\'et\'es associ\'ees \`a 
une suite  d'\'el\'ements contenus  dans un sous-groupe  parabolique $F$-stable
du groupe de  Weyl comme induites  \`a partir  de vari\'et\'es 
associ\'ees  \`a un  sous-groupe de  Levi. Pour  $I\subset S$,  on note 
$\bP_I$  le  sous-groupe  parabolique  $\bB W_I\bB$,  on  note  $\bU_I$ 
son  radical unipotent,  on  note  $\bL_I$ le compl\'ement de Levi
contenant $\bT$  et on  note $\bB_I$  le sous-groupe  de Borel 
$\bB\cap\bL_I$ de $\bL_I$. Rappelons que si $\bB_1$ est un sous-groupe 
de Borel de $\bL_I$, alors $\bB_1\bU_I$ est un sous-groupe  de Borel de $\bG$.
On note $p_{\bL_I}$ la projection de $\bP_I$ sur $\bL_I\iso \bP_I/\bU_I$.

\begin{proposition}\label{RLG(Xs)}
Soit $I\subset S$ une partie $F$-stable et
soient $t_1,\ldots,t_k\in W_I\cup\uW_I$.
Alors, l'application
\begin{align*}
\GF/\bU_I^F\times_{\bL_I^F}
\bX_{\bL_I}(t_1,\ldots,t_k)&\rightarrow
\bX_\bG(t_1,\ldots,t_k)\\
x\bU_I^F\times_{\bL_I^F}(\bB_1,\ldots,\bB_k)&\mapsto (\lexp x(\bB_1\bU_I),
\ldots,\lexp x(\bB_k\bU_I))\\
\end{align*}
d\'efinit un isomorphisme de vari\'et\'es compatible avec l'action de
$\GF\times F^m$ pour tout $m$ tel que $F^m$ fixe $(t_1,\ldots,t_k)$.
\end{proposition}

\begin{pf*}{Preuve} Le morphisme $\gamma$ de l'\'enonc\'e est compatible \`a
l'action de $\GF\times F^m$ pour $m$ comme ci-dessus.
Soit $l_i\in \bL_I$ tel que $\bB_i=\lexp{l_i}\bB_I$. Le morphisme peut alors se
r\'e\'ecrire:
$$ x\bU_I^F\times_{\bL_I^F}(\lexp{l_1}\bB_I,\ldots,\lexp{l_k}\bB_I)\mapsto
(\lexp{xl_1}\bB,\ldots,\lexp{xl_k}\bB).$$

Nous allons montrer que $\gamma$ est bijectif. 
Consid\'erons l'application
\begin{multline*}\{(g_1,\ldots,g_k)\in\bG^k\mid g_i\inv g_{i+1}\in\bB
t_i\bB,\ i\ne k; \ g\inv_k F(g_1)\in\bB
t_k\bB\}\rightarrow\\
\GF/\bU_I^F\mathop\times\limits_{\bL_I^F}
\{(l_1,\ldots,l_k)\in\bL^k_I\mid
l_i\inv l_{i+1}\in\bB_I t_i\bB_I,\
i\ne k; \ l\inv_k F(l_1)\in\bB_I t_k\bB_I\},
\end{multline*}
donn\'ee
par $(g_1,\ldots,g_k)\mapsto x\bU_I^F\times_{\bL_I^F}(l_1,\ldots,l_k)$, o\`u
$l_1\inv F(l_1)=p_{\bL_I}(g_1\inv F(g_1))$,
$l_i=l_{i-1}p_{\bL_I}(g_{i-1}\inv g_i)$ pour $i=2,\ldots,k$ et
$x\bU_I^F=\GF\cap g_1l_1\inv\bU_I$.

On a $g_1\inv F(g_1)=g_1\inv g_2g_2\inv\cdots
g_k\inv F(g_1)\in\bL_I\bU_I$ et $g_1l_1\inv\in \GF\bU_I$, donc
l'application est bien d\'efinie:
$l_1$ est d\'efini \`a multiplication par $l\in \bL_I^F$ \`a gauche pr\`es.
Mais alors, 
tous les $l_i$ sont multipli\'es par $l$
\`a gauche et donc $x$ par $l\inv$ \`a droite.

L'action de $\bB$ par multiplication \`a droite sur les $g_i$ est
transport\'ee par cette application sur l'action de $\bB_I=\bB/\bU_I$
par multiplication \`a droite sur les $l_i$. On a donc une application induite
$f:\bX_\bG(t_1,\ldots,t_k)\to
\GF/\bU_I^F\mathop\times\limits_{\bL_I^F} \bX_{\bL_I}(t_1,\ldots,t_k)$
qui est un inverse \`a gauche de $\gamma$.
Puisque $xl_1\in g_1\bB$ et
$g_il_i\inv \bU_I=g_{i-1}l_{i-1}\inv \bU_I$ pour $i\ge 2$,
on a $xl_i\in g_i \bB$ pour tout $i$. Par cons\'equent, $f$ est inverse
\`a droite de $\gamma$.

Le morphisme produit $(x,l_1,\ldots,l_k)\mapsto (xl_1,\ldots,xl_k):
\bG^F\times \bG^k\to \bG^k$ est s\'eparable. Il le reste apr\`es passage au
quotient par $\bU_I^F\times \bB_I^k$ puis par $\bL_I^F$ et enfin par
$\bU_I^k$ \`a l'arriv\'ee.
Par restriction \`a la sous-vari\'et\'e ferm\'ee
$\bX_\bG(t_1,\ldots,t_k)$ et \`a son image inverse
$\GF/\bU_I^F\mathop\times_{\bL_I^F} \bX_{\bL_I}(t_1,\ldots,t_k)$,
le morphisme reste s\'eparable.

Le morphisme $\gamma$ \'etant bijectif,
s\'eparable et la vari\'et\'e d'arriv\'ee \'etant normale (proposition
\ref{desingularisation}),
c'est bien un isomorphisme (th\'eor\`eme \ref{isonormale}). \end{pf*}

Nous  allons maintenant  introduire une  technique (proposition \ref{X_w  produit de
varietes}) qui  nous permettra d'\'etudier par
r\'ecurrence les vari\'et\'es de Deligne-Lusztig (g\'en\'eralis\'ees) et qui
g\'en\'eralise le cas particulier de l'\'enonc\'e pr\'ec\'edent o\`u tous les $t_i$ sont
dans $\bW_I$.

Pour cela nous avons besoin de deux lemmes sur le mono{\"\i}de de tresses.

\begin{lemme}\label{ws in B+} Soient $\bw\in B^+$ et $\bs\in \bS$
tels que $\lexp\bw\bs\in B^+$ et soit $\bw_1\cdots\bw_k$ la forme
normale de $\bw$. Alors, pour tout $i$, on a
$\lexp{\bw_i\cdots\bw_k}\bs\in\bS$.
\end{lemme}
\begin{pf*}{Preuve}
Nous allons d\'emontrer le lemme par r\'ecurrence sur $k$.
L'\'el\'ement $\lexp\bw\bs\in B^+$ est  de longueur $1$,
donc $\lexp\bw\bs=\bt\in\bS$. On a $\bw_1\cdots\bw_k
\bs=\bt\bw_1\cdots\bw_k$. Par \cite[lemme 4.6]{michel}, il existe
 $\bw'_1\preccurlyeq \bw_1$ tel que $\alpha(\bt\bw_1\cdots\bw_k)=\bt\bw'_1$.
L'\'egalit\'e pr\'ec\'edente montre que $\bw_1\preccurlyeq\alpha(\bt\bw_1\cdots\bw_k)=
\bt\bw'_1\preccurlyeq\bt\bw_1$.  Il y a donc deux possibilit\'es:
\begin{enumerate}
\item $\bw'_1=\bw_1$ et alors $\bt\bw_1=\bw_1\bs_1$ pour un
certain  $\bs_1\in\bS$. Alors $\bw_2\cdots\bw_k\bs=\bs_1\bw_2\cdots\bw_k$
et on  conclut par r\'ecurrence.
\item Il existe $\bs_1\in\bS$ tel que
$\bw_1=\bw'_1\bs_1=\bt\bw'_1$. Ici encore on en d\'eduit que 
$\bw_2\cdots\bw_k\bs=\bs_1\bw_2\cdots\bw_k$ et on conclut de m\^eme.
\end{enumerate}
\end{pf*}

\begin{lemme}\label{w in N_B(BI)} Soient $\bw\in B^+$ et $I,
J\subset S$ tels que $\lexp \bw B_I\subset B_J$.
Alors, $\lexp{\omega_J(\bw)}\bI\subset\bJ$ et l'image $\beta
(\omega_J(\bw))$ dans $W$ est $J$-r\'eduite.
\end{lemme}
\begin{pf*}{Preuve}
De  l'hypoth\`ese  il r\'esulte   que  si $\bs\in\bI$,  alors
$\lexp{\omega_J(\bw)}\bs\in B_J$.  Posons $\lexp{\omega_J(\bw)}\bs
=\ba\inv\bb$  o\`u $\ba,\bb\in  B_J^+$  n'ont pas  de diviseur \`a gauche
commun (non trivial)
(ce  qui  est   toujours  possible  par  \cite[corollaire 3.2]{michel}).  On  a
$\alpha_J(\bb\omega_J(\bw))=\bb$, car,  si on  \'ecrit $\bb\omega_J
(\bw)=\bb\bx\by$  o\`u  $\bb\bx=  \alpha_J  (\bb\omega_J(\bw))$,  en
simplifiant   on   a   $\bx\preccurlyeq  \omega_J(\bw)$   ce   qui   implique
$\bx=1$  puisque   par  d\'efinition  aucun  \'el\'ement   de  $B^+_J$
ne   divise  $\omega_J(\bw)$.   Donc  tout \'el\'ement $\bt$ de $\bJ$  qui   divise
$\bb\omega_J(\bw)$  divise  $\bb$.  Mais  comme  on  a  $\ba\omega_J
(\bw)\bs=\bb\omega_J(\bw)$, tout  $\bt\in\bJ$ qui divise  $\ba$ divise
$\bb\omega_J(\bw)$; le fait que $\ba$ et $\bb$ n'aient pas de diviseur
commun implique alors que $\ba=1$. Alors, $l(\bb)=1$, donc
$\bb\in\bJ$ et $\lexp{\omega_J(\bw)}\bs\in\bJ$,  ce
qui d\'emontre la premi\`ere  assertion.

Nous d\'emontrons la deuxi\`eme
par  r\'ecurrence  sur le  nombre  de  termes  de  la forme  normale  de
$\omega_J(\bw)$,  l'assertion  \'etant  claire   si  cette  forme  n'a
qu'un  terme.  Notons $\bw_1\cdots\bw_k$  cette  forme  normale. Par  le
lemme \ref{ws  in B+},
$\bJ_1=\lexp{\bw_2\cdots \bw_k}\bI$ est inclus dans $\bS$.
La premi\`ere partie du lemme montre que
$\lexp{\omega_{J_1}(\bw_2\cdots\bw_k)}\bI\subset \bJ_1$, donc
$\omega_{J_1}(\bw_2\cdots\bw_k)=\bw_2\cdots\bw_k$. On en d\'eduit que
$\beta(\bw_2\cdots\bw_k)$ est  $J_1$-r\'eduit par l'hypoth\`ese
de r\'ecurrence.
On est ramen\'e \`a prouver  le r\'esultat dans le cas o\`u la
forme normale a  deux termes, c'est-\`a-dire \`a prouver  que si $v_1\in
W$ est un \'el\'ement $J$-r\'eduit qui conjugue $J_1$ sur $J$ et $v_2\in
W$  un  \'el\'ement $J_1$-r\'eduit  qui  conjugue  $I$ sur  $J_1$  alors
$v_1v_2$ est $J$-r\'eduit. Un \'el\'ement  $s$ de $J$ ne divise pas $v_1$,
donc divise $v_1v_2$  si et seulement s'il existe un  \'el\'ement de $S$
qui divise $v_2$ et dont le  conjugu\'e par $v_1$ vaut $s$. Ceci est
impossible car $\lexp{v_1\inv}s\in J_1$ et aucun \'el\'ement de $J_1$ ne
divise $v_2$.
\end{pf*}

\begin{definition}\label{RLG}
Soit $\bL$  (resp. $\bU$)  un  compl\'ement de  Levi (resp.  le radical
unipotent) d'un sous-groupe parabolique de $\bG$ et $n\in\bG$ tel
que $\bL$ est $nF$-stable. On pose
$$\tilde\bX^{\bL,\bU}(n)=
\{g\bU\in\bG/\bU\mid g\inv F(g)\in\bU n F(\bU)\}$$
(nous noterons  $\tilde\bX^{\bL,\bU}_\bG(n,F)$ cette vari\'et\'e  quand nous
voudrons pr\'eciser  $\bG$ et  $F$). Sur  cette vari\'et\'e  les groupes
$\GF$  et  $\bL^{nF}$  agissent respectivement  par  multiplication  \`a
gauche  et  \`a  droite.
\end{definition}

\begin{definition}\label{Xtilde_n}
Soient $\bw\in B^+$ et $I\subset S$ tels que $\lexp{\bw }F(B_I)=B_I$.
On pose $\bz=\omega_I(\bw)$, on note
$\bz_1\cdots \bz_k$ la forme normale de $\bz$ et on pose
$\bI_j=\lexp{\bz_j\cdots\bz_k}F(\bI)$ (on a $\bI_1=\bI$ et
$\bI_j\subseteq\bS$ par les lemmes
\ref{ws in B+} et \ref{w in N_B(BI)}).
Soient $\dz_1$, $\ldots$, $\dz_k$ des relev\'es dans $N_\bG(\bT)$ de 
$z_1$, $\ldots$, $z_k$;
on note $\tilde\bX_{(I)}(\dz_1,\ldots,\dz_k)$ la vari\'et\'e
$$\tilde\bX^{\bL_{I_1\times\cdots\times I_k},
\bU_{I_1\times\cdots\times I_k}}_{\bG^k}((\dz_1,\ldots,\dz_k),F_1)$$
o\`u $F_1(g_1,\ldots,g_k)=(g_2,\ldots,g_k,F(g_1))$.
\end{definition}
Autrement dit,
\begin{multline*}\tilde\bX_{(I)}(\dz_1,\ldots,\dz_k)=
\{(g_1\bU_{I_1},\ldots,g_k\bU_{I_k}) \mid \\
g_i\inv g_{i+1}\in \bU_{I_i}\dz_i\bU_{I_{i+1}}\text{ pour }i<k
\text{ et }g_k\inv F(g_1)\in\bU_{I_k}\dz_k F(\bU_{I_1})\}.\\
\end{multline*}

Le   groupe    $\GF\simeq(\bG^k)^{F_1}$   (\cf\  proposition \ref{produit})   agit
par   multiplication   \`a   gauche    sur   cette   vari\'et\'e.   Soit
$\dz=\dz_1\cdots\dz_k$. La   premi\`ere   projection induit un
isomorphisme
$\bL_{I_1\times\cdots\times    I_k}^{(\dz_1,\ldots,\dz_k)F_1}
\iso\bL_I^{\dz F}$. Ceci fournit une action \`a droite de 
$\bL_I^{\dz F}$ sur $\tilde\bX_{(I)}(\dz_1,\ldots,\dz_k)$.

\begin{proposition}\label{X_w produit de varietes} 
Sous  les hypoth\`eses et avec les notations de la d\'efinition
\ref{Xtilde_n}, soit
$\by=\alpha_I(\bw)$
et soit  $\by=\by_1\cdots\by_h$ une  d\'ecomposition de $\by$
en  \'el\'ements de $\bW$.
Alors, l'application
\begin{multline*}
((g_1\bU_{I_1},\ldots,g_k\bU_{I_k}),(\bB_1,\ldots,\bB_h))\mapsto 
(\lexp{g_1}(\bB_1\bU_I),\lexp{g_1}(\bB_2\bU_I),\ldots,\lexp{g_1}(\bB_h\bU_I),\\
\lexp{g_1}(\lexp z F(\bB_1)\bU_{I_1}),
\lexp{g_2}(\lexp{ z_2\cdots z_k}F(\bB_1)\bU_{I_2}),
\ldots,\lexp{g_k}(\lexp{z_k}F(\bB_1)\bU_{I_k}))\\
\end{multline*}
d\'efinit un isomorphisme
$$\tilde\bX_{(I)}(\dz_1,\ldots,\dz_k)\times_{\bL_I^{\dz F}}
\bX_{\bL_I}(y_1,\ldots,y_h,\dz F)\iso 
\bX(y_1,\ldots, y_h,z_1,\ldots,z_k)$$
Via les isomorphismes canoniques, on en d\'eduit un isomorphisme
$$\tilde\bX_{(I)}(\dz_1,\ldots,\dz_k)\times_{\bL_I^{\dz F}}
\bX_{\bL_I}(\by,\dz F)\iso \bX(\bw)$$
qui ne  d\'epend pas de la  d\'ecomposition choisie de
$\by$ et qui est compatible avec les  actions de $\GF$ et de $F^n$ pour
tout $n$ tel que $I$, $\by$ et $(\dz_1,\ldots,\dz_k)$ soient $F^n$-stables.
\end{proposition}
\begin{pf*}{Preuve}
Nous pouvons simplifier la formule du morphisme en introduisant $l\in\bL_I$ tel
que $\bB_1=\lexp l\bB_I$ et en posant
$l_i=\lexp{\dz_i\cdots\dz_k}F(l)\in\bL_{I_i}$.
On a alors $\lexp{z_i\cdots z_k}F(\bB_1)\bU_{I_i}=\lexp{l_i}\bB$
car $\bB_{I_i}\bU_{I_i}=\bB$ et $\lexp{z_i\cdots z_k}F(\bB_1)=
\lexp{l_i}(\lexp{z_i\cdots z_k}F(\bB_I))=\lexp{l_i}\bB_{I_i}$.
Le morphisme s'\'ecrit donc 
\begin{multline*}
((g_1\bU_{I_1},\ldots,g_k\bU_{I_k}),(\bB_1,\ldots,\bB_h))\mapsto\\
(\lexp{g_1}(\bB_1\bU_I),\lexp{g_1}(\bB_2\bU_I),\ldots,\lexp{g_1}(\bB_h\bU_I),
\lexp{g_1l_1}\bB,\lexp{g_2l_2}\bB,\ldots,\lexp{g_kl_k}\bB)\quad(*)\\
\end{multline*}
Sous cette forme on v\'erifie qu'il est \`a valeurs dans
$\bX(y_1,\ldots,y_h,z_1,\ldots,z_k)$; par exemple, en \'ecrivant que
$g_i\inv g_{i+1}\in \bU_{I_i}\dz_i \bU_{I_{i+1}}$ on trouve que
$$(g_il_i)\inv(g_{i+1}l_{i+1})\in \bU_{I_i}l_i\inv \dz_i l_{i+1}\bU_{I_{i+1}}
\subset\bB z_i\bB \text{ (car $l_i\inv \dz_i l_{i+1}=\dz_i$),}$$
donc 
$(\lexp{g_il_i}\bB,\lexp{g_{i+1}l_{i+1}}\bB)$ est bien dans $\CO(z_i)$.

Remarquons maintenant qu'on peut se ramener au cas $k=1$.
En effet par la proposition \ref{produit}, on a un isomorphisme
$$\bX(\bw)\iso\bX_{\bG^k}((\by\bz_1,\bz_2,\ldots,\bz_k),F_1)$$
o\`u $(\by\bz_1,\bz_2,\ldots,\bz_k)\in (\bB^+)^k$ est un \'el\'ement dont le
$\omega_{I_1\times\cdots\times I_k}$ est \linebreak
$(\bz_1,\bz_2,\ldots,\bz_k)\in \bW^k$;
l'image de $(\bB_1,\ldots,\bB_{h+k})\in
\bX(y_1,\ldots,y_h,z_1,\ldots,z_k)$ par cet isomorphisme est l'\'el\'ement
de $$\bX_{\bG^k}((y_1,1,\ldots,1),\ldots,(y_h,1,\ldots,1),(z_1,\ldots,z_k),F_1)$$
donn\'e par
$$((\bB_1,\bB_{h+2},\ldots,\bB_{h+k}),(\bB_2,\bB_{h+2},\ldots,\bB_{h+k}),\ldots
  (\bB_{h+1},\bB_{h+2},\ldots,\bB_{h+k})),$$
et on v\'erifie que cet isomorphisme est compatible avec la formule $(*)$.

Nous sommes donc ramen\'es au cas o\`u $\bw=\by\bz$ avec $\bz\in\bW$, $\lexp\bz
\bI=\bI$; on a $\tilde{\bX}_{(I)}(\dz)=\{g\bU_I\mid g\inv F(g)\in\bU_I\dz
F(\bU_I)\}$ et le morphisme qui va de
$\tilde{\bX}_{(I)}(\dz)\times_{\bL_I^{\dz F}}\bX_{\bL_I}(y_1,\ldots,y_h,\dz F)$
vers $\bX(y_1, \ldots,y_h,z)$
est donn\'e par
$$ (g\bU_I,(\bB_1,\ldots,\bB_h))\mapsto
(\lexp g(\bB_1\bU_I),\lexp g(\bB_2\bU_I),\ldots,\lexp g(\bB_h\bU_I),
\lexp g(\lexp z F(\bB_1)\bU_I)).\eqno(**)$$
On se ram\`ene alors au cas $h=1$ en utilisant l'isomorphisme
donn\'e par la proposition \ref{produit},
$$\bX(\bw)\iso\bX_{\bG^h}((\by_1,\ldots,\by_{h-1},\by_h\bz),F_1)$$
(notons que $\by_h\bz\in\bW$ d'apr\`es le lemme \ref{w in N_B(BI)});
on v\'erifie que cet isomorphisme est compatible avec $(**)$.

Nous sommes donc ramen\'es au cas $\bw\in\bW$, o\`u comme
expliqu\'e apr\`es la d\'efinition \ref{DL}, on a $\bX(w)\iso
\{h\bB\mid h\inv F(h)\in\bB w\bB\}$
et $\bX_{\bL_I}(y,\dz F)\iso\{l\bB_I\mid l\inv\lexp{\dz}F(l)\in\bB_I y\bB_I\}$.
Il s'agit alors de montrer que le morphisme $j_w$ donn\'e par 
$(g\bU_I,l\bB_I)\mapsto gl\bB$
r\'ealise un isomorphisme $\tilde\bX_{(I)}(\dz)\times_{\bL_I^{\dz F}}
\bX_{\bL_I}(y,\dz F)\xrightarrow{j_w}\bX(w)$.

Or le diagramme suivant est commutatif:
$$\xymatrix{
\bG/\bU_I\times_{\bL_I}\bL_I/\bB_I\ar[r]^-j & \bG/\bB\\
\tilde\bX_{(I)}(\dz)\times_{\bL_I^{\dz F}}\bX_{\bL_I}(y,\dz F)
\ar[r]_-{j_w}\ar[u]^i & \bX(w)\ar[u]_{i_w}
}$$
o\`u $i_w$ est l'inclusion canonique, $i$ l'application induite par les
deux inclusions canoniques $\tilde\bX_{(I)}(\dz)\hookrightarrow\bG/\bU_I$
et  $\bX_{\bL_I}(y,\dz  F)\hookrightarrow\bL_I/\bB_I$   et  o\`u  $j$,
donn\'e par  la m\^eme  formule que  $j_w$, est  un isomorphisme (il
admet comme inverse
$g\bB\mapsto(g\bU_I,\bB_I)$). Montrons que  $i$ est injectif:
si $(h\bU_I,l\bB_I)$ et $(hl'\bU_I,l'l\bB_I)$ sont tous deux dans
$\tilde\bX_{(I)}(\dz)\times_{\bL_I^{\dz F}}\bX_{\bL_I}(\by,\dz F)$
avec $l'\in\bL_I$ alors $h\inv F(h)$ et $(hl')\inv F(hl')$
sont dans $\bU_I\dz F(\bU_I)$ et
l'intersection $(\bU_I\lexp{\dz}F(\bU_I))\cap(\bU_I l'\lexp{\dz}F(l^{\prime-1})
\lexp{\dz}F(\bU_I))$ n'est pas vide, puisqu'elle contient
$h\inv F(h)\dz\inv$.
Or $\bL_I\iso \bU_I\setminus \bP_I\lexp zF(\bP_I)/\lexp{\dz}F(\bU_I)$,
donc $l'\in\bL_I^{\dz F}$. Par cons\'equent, $i$ est une immersion
ferm\'ee.

Montrons maintenant que $j_w$ est surjectif.
Soit $h\in\bG$ tel que $h\inv F(h)\in \bB w\bB$.
On a
$\bB w\bB=\bU_I\bB_I yz F(\bB_I)F(\bU_I)=\bU_I\bB_I y\bB_I z F(\bU_I)$
car $\lexp zF(\bB_I)=\bB_I$.
Soit $l'\in\bB_I y\bB_I $ tel que
$h\inv F(h)\in \bU_I l' \dz F(\bU_I)$. Puisque
$\lexp{\dz}F(\bL_I)=\bL_I$, il existe $l\in\bL_I$ tel que
$l\inv \lexp{\dz}F(l)=l'$. On a alors
$h\inv F(h)\in l\inv \bU_I \dz F(\bU_I)F(l)$.
Soit $g=h l\inv$. On a 
$g\inv  F(g)\in \bU_I \dz F(\bU_I)$, donc
$(g,l)\in \tilde\bX_{(I)}(\dz)\times_{\bL_I^{\dz F}}\bX_{\bL_I}(\by,\dz F)$.
On a donc montr\'e que $j_w$ est surjectif.

On en d\'eduit que $j_w$ s'obtient \`a partir de $j$ par changement de base, donc
c'est un isomorphisme.
\end{pf*}
\section{Cohomologie}
\label{sectioncohomologie}

On fixe un nombre premier $\ell$ diff\'erent de $p$.
Pour la commodit\'e de l'exposition, on fixe un isomorphisme $\Qlbar\iso\bbC$.
Pour tout $\bt\in\uB^+$, la vari\'et\'e $\bX(\bt)$ est munie d'actions
de $\GF$ et de $F^\delta$ qui commutent. Ces actions sont donn\'ees
par des morphismes propres. Elles
induisent donc des actions qui commutent sur les groupes de cohomologie
$\ell$-adique \`a support propre $H^i_c(\bX(\bt),\Qlbar)$.
Puisque l'endomorphisme $F^\delta$ induit une \'equivalence de sites \'etales,
il agit par un endomorphisme inversible sur la cohomologie. Nous avons donc
une repr\'esentation du groupe $\GF\times\genby{F^\delta}$ sur la
cohomologie.  Nous nous int\'eressons dans cette partie
au $\Qlbar(\GF\times\genby {F^\delta})$-module gradu\'e
$H^*_c(\bX(\bt),\Qlbar)$.
Nous omettrons
dans la suite les coefficients quand ils sont \'egaux \`a $\Qlbar$ et
\'ecrirons simplement $H^i_c(\bX(\bt))$.

\subsection{Constructions}
\sub{}
Pour $\sigma:A\to B$ un morphisme d'alg\`ebres, nous notons
$\sigma^*$ le foncteur de restriction \`a travers $\sigma$
de la cat\'egorie des $B$-modules vers celle des $A$-modules.

Pour $C$ un complexe et $i\in\bbZ$, nous notons $C[i]$ le complexe
d\'ecal\'e de $i$ crans vers la gauche.
Pour $k$ un anneau commutatif et
$M$ un $k\langle F^\delta\rangle$-module,  nous notons
$M(i)$ le $k\langle F^\delta\rangle$-module 
$\sigma^* M$, o\`u $\sigma:k\langle F^\delta\rangle\to k\langle
F^\delta\rangle$ est l'automorphisme qui envoie $F^\delta$ sur
$q^{-i\delta}F^\delta$.

\sub{}
Nous allons passer en revue un certain nombre de propri\'et\'es des
$(\GF\times\genby{F^\delta})$-modules $H^i_c(\bX(\bt))$.

La proposition \ref{RLG(Xs)} admet le corollaire
suivant (la cohomologie \og s'induit\fg\ d'un sous-groupe de Levi au groupe).

\begin{corollaire}\label{RLG(H(Xs))}
Soit $I$ une partie $F$-stable de $S$.
Alors,
pour toute suite $t_1,\ldots,t_k$ d'\'el\'ements de $W_I\cup\uW_I$,
on a isomorphisme de $(\GF\times F^\delta)$-modules gradu\'es
$$R_{\bL_I}^\bG(H^*_c(\bX_{\bL_I}(t_1,\ldots,t_k)))\iso
H^*_c(\bX_\bG(t_1,\ldots,t_k)).$$
\end{corollaire}
Nous avons not\'e $R_{\bL_I}^\bG$ l'induction de Harish-Chandra d\'efinie
sur les $\bL_I^F$-modules par $\Ind_{\bP_I^F}^\GF\circ p_{\bL_I^F}^*$
o\`u $p_{\bL_I^F}$ est la restriction de $p_{\bL_I}$ en un morphisme de
groupes $\bP_I^F\to \bL_I^F$.

\sub{}
Nous g\'en\'eralisons maintenant une partie de \cite[th\'eor\`eme 1.6]{DL}.
Pour $\bt',\bt''\in\uB^+$, on a un isomorphisme canonique
$$\bX(\bt'\bt'')\iso\{(a,b)\in\CO(\bt')\times\CO(\bt'')\mid p''(b)=F(p'(a))
\text{ et }p''(a)=p'(b)\}.$$

\begin{definition}\label{Dt} Pour $\bt',\bt''\in\uB^+$,
nous d\'efinissons le morphisme $D_{\bt'}:\bX(\bt'\bt'')\to \bX(\bt'' F(\bt'))$
comme la restriction du morphisme 
$\CO(\bt')\times\CO(\bt'')\to\CO(\bt'')\times\CO(F(\bt')),(a,b)\mapsto(b,F(a))$.
\end{definition}

Choisissons des d\'ecompositions $\bt'=\bt_1\cdots\bt_k$ et
$\bt''=\bt_{k+1}\cdots\bt_n$ avec $t_i\in W\cup\uW$.
Alors, le morphisme $$D_{\bt'}:\bX(t_1,\ldots,t_n)
\to\bX(t_{k+1},\ldots,t_n,F(t_1),\ldots,F(t_k))$$
est donn\'e par
$D_{\bt'}(\bB_1,\ldots,\bB_n)=(\bB_{k+1},\ldots,\bB_n,
F(\bB_1),\ldots,F(\bB_k))$.

\begin{proposition}\label{xy=yFx} Soient $\bt',\bt''\in\uB^+$.
Alors, le morphisme $D_{\bt'}$ induit une \'equivalence de sites \'etales
$\bX(\bt'\bt'')\xrightarrow{D_{\bt'}}\bX(\bt'' F(\bt'))$.
Il induit donc un isomorphisme de
$(\GF\times\genby{F^\delta})$-modules gradu\'es
$$H^*_c(\bX(\bt'\bt''))\iso H^*_c(\bX(\bt'' F(\bt'))).$$
\end{proposition}

\begin{pf*}{Preuve} La d\'emonstration de \cite[theorem 1.6 (case 1)]{DL}
s'applique ici. On a un diagramme commutatif
$$\xymatrix{
  \bX(\bt'\bt'')\ar[rr]^{D_{\bt'}}\ar[d]_F && \bX(\bt'' F(\bt'))\ar[d]^F
  \ar[dll]_{D_{\bt''}}\\
  \bX(F(\bt')F(\bt''))\ar[rr]_{D_{F(\bt')}} &&
  \bX(F(\bt'')F^2(\bt'))\\}
$$
Puisque $F^\delta$ est une \'equivalence de sites \'etales, $F$ aussi en est
une; les fl\`eches verticales sont donc des \'equivalences de sites \'etales,
il en est alors de m\^eme des fl\`eches horizontales.
\end{pf*}

Deux \'el\'ements $\bt$ et $\bt'$ de $B$ sont {\em $F$-conjugu\'es}
s'il existe $\bw\in B$ tel que $\bt'=\bw^{-1}\bt F(\bw)$.
La proposition \ref{xy=yFx}
montre l'invariance de la cohomologie $H^*_c(\bX(\bt))$
pour une forme particuli\`ere de $F$-conjugaison. 
\begin{conjecture}\label{groupe de tresses}
La classe d'isomorphisme du 
$(\GF\times\genby{F^\delta})$-module gradu\'e
$H^*_c(\bX(\bt))$
ne d\'epend que de la $F$-classe de conjugaison de $\bt\in B^+$.
\end{conjecture}
Nous le prouverons pour $\bG$ de type $A_2$, \cf\ corollaire \ref{H fclasse}.

\begin{proposition}\label{autG}
Soit $\sigma$ une isog\'enie sur $\bG$ commutant \`a $F$. Notons encore
$\sigma$ l'automorphisme correspondant de $S$ et celui qui s'en d\'eduit
sur $\uB^+$.  Alors, pour tout $\bt\in\uB^+$, l'isog\'enie $\sigma$
induit
un isomorphisme de $(\GF\times\genby{F^\delta})$-modules gradu\'es
$H^*_c(\bX(\bt))\iso\sigma^*H^*_c(\bX(\sigma(\bt)))$.
\end{proposition}

\begin{pf*}{Preuve} La proposition se d\'eduit du fait qu'une isog\'enie
induit une \'equivalence de sites \'etales.
\end{pf*}

\subsection{D\'evissages}
\label{sectiondevissages}
\sub{}

Pour $j:\bU\to\bX$ l'inclusion d'une sous-vari\'et\'e localement ferm\'ee,
on notera $\Lambda_\bU$ le faisceau sur $\bX$ donn\'e par $j_!(\Qlbar)$.
Si $j$ est une immersion
ouverte et $i:\bZ\to\bX$ est l'immersion ferm\'ee compl\'ementaire, on
a une suite exacte de faisceaux sur $\bX$
$$0 \to \Lambda_\bU \to \Lambda_\bX \to \Lambda_\bZ \to 0,$$
donc une suite exacte longue
$$\cdots\to H^i_c(\bU,\Qlbar)\to H^i_c(\bX,\Qlbar)\to H^i_c(\bZ,\Qlbar)\to
H^{i+1}_c(\bU,\Qlbar)\to\cdots.$$

Par faisceau sur une vari\'et\'e $\bX(\bt)$ nous entendrons 
un faisceau $\bG^F$-\'equivariant et il en sera
de m\^eme pour les morphismes.

\begin{proposition}\label{ouvertferme}
Soient   $w\in   W$   et  $\bt,\bt'\in\uB^+$.   Alors   la   vari\'et\'e
$\bX(\bt\ubw\bt')$   admet    une   filtration    par des   sous-vari\'et\'es
ferm\'ees      $\bX_0\subset\bX_1\subset\cdots\subset\bX_{l(w)}=
\bX(\bt\ubw\bt')$    o\`u
$\bX_i=\coprod_{v\le   w,   l(v)\le   i}\bX(\bt   \bv\bt')$.   De   plus
$\bX_i-\bX_{i-1}=\coprod_{v\le    w,l(v)=i}\bX(\bt\bv\bt')$   est    une
d\'ecomposition   en union de composantes  connexes.   On en d\'eduit une
filtration du faisceau   constant
$\Lambda_{\bX(\bt\ubw\bt')}$  dont    les
quotients  successifs  sont  $$\Lambda_{\bX_i-\bX_{i-1}}=\bigoplus_{v\le
w,l(v)=i}\Lambda_{\bX(\bt\bv\bt')}.$$
\end{proposition}

\begin{pf*}{Preuve}
On  obtient  la  filtration  de $\Lambda_{\bX(\bt\ubw\bt')}$
\`a   partir  des  suites  exactes 
associ\'ees  aux  d\'ecompositions  \og ouvert-ferm\'e\fg\ provenant  de  la
stratification en sous-vari\'et\'es ferm\'ees.
\end{pf*}

Soient $w\in W$ et $s\in S$ tels que $ws<w$. Alors, le morphisme
$(p',p''):\CO(\ubw\ \ubs)\to \CB\times\CB$ se factorise par l'immersion
ferm\'ee $\CO(\ubw)\to\CB\times\CB$. Par produit fibr\'e et restriction,
on en d\'eduit un morphisme canonique
$\bX(\bt\ubw\ \ubs\bt')\to \bX(\bt\ubw\bt')$ pour $\bt, \bt'\in\uB^+$.

\begin{proposition}\label{wss}
Soient $\bt, \bt'\in\uB^+$, $w\in W$ et $s\in S$ tels que $ws<w$.
\begin{enumerate}
\item Soit
$\pi:\bX(\bt\ubw \bs\bt')\to \bX(\bt\ubw\bt')$ l'application canonique.
Alors,
$$R\pi_!\Lambda_{\bX(\bt\ubw \bs\bt')}\simeq\Lambda_{\bX(\bt\ubw\bt')}[-2](-1).$$
En particulier, pour tout $i$, on a un isomorphisme de
$(\GF\times\genby{F^\delta})$-modules:
$$H^i_c(\bX(\bt\ubw \bs\bt'))\simeq H^{i-2}_c(\bX(\bt\ubw\bt'))(-1).$$
\item Soit
$\pi':\bX(\bt\ubw\ \ubs\bt')\to \bX(\bt\ubw\bt')$ l'application canonique.
Alors,
$$R\pi'_!\Lambda_{\bX(\bt\ubw\ \ubs\bt')}\simeq
 \Lambda_{\bX(\bt\ubw\bt')}[-2](-1)\oplus \Lambda_{\bX(\bt\ubw\bt')}.$$
En particulier, pour tout $i$, on a un isomorphisme de
$(\GF\times\genby{F^\delta})$-modules:
$$H^i_c(\bX(\bt\ubw\ \ubs\bt'))\simeq H^{i-2}_c(\bX(\bt\ubw\bt'))(-1)\oplus
H^i_c(\bX(\bt\ubw\bt')).$$
\end{enumerate}
\end{proposition}

\begin{pf*}{Preuve} L'application $\pi$ provient par changement de base
par $\Gamma\to\CB\times\CB$ (\cf\ d\'efinition \ref{DL})
de l'application canonique analogue 
$\tilde\pi:\CO(\bt\ubw \bs\bt')\to \CO(\bt\ubw\bt')$.
La deuxi\`eme  projection $\CO(\bs)\xrightarrow{p''}\CB$ est  une fibration
en droites affines. On a un carr\'e cart\'esien
$$\xymatrix{
\CO(\bt\ubw \bs\bt')\ar[r]\ar[d]_{\tilde\pi} & \CO(\bs) \ar[d]^{p''}\\
\CO(\bt\ubw \bt')\ar[r]          & \CB
}$$
Par cons\'equent, le morphisme $\tilde\pi$ est une fibration en droites
affines. Par changement de base, on d\'eduit (i).

De m\^eme l'application $\pi'$
provient par changement de base 
d'une application analogue $\tilde\pi'$. 
La restriction de $\tilde\pi'$ au ferm\'e $\CO(\bt\ubw\bt')$ est un isomorphisme
et la restriction \`a $\CO(\bt\ubw\bs\bt')$ est $\tilde\pi$.
La suite exacte \og ouvert-ferm\'e\fg
$$0\to \Lambda_{\CO(\bt\ubw \bs\bt')}\to \Lambda_{\CO(\bt\ubw\ \ubs\bt')} \to
\Lambda_{\CO(\bt\ubw\bt')} \to 0$$
devient donc, par application de $R\tilde\pi'_!$, un triangle distingu\'e
$$\Lambda_{\CO(\bt\ubw\bt')}[-2](-1) \to
R\tilde\pi'_! \Lambda_{\CO(\bt\ubw\ \ubs\bt')}\to
\Lambda_{\CO(\bt\ubw\bt')}\rightsquigarrow$$
Le morphisme associ\'e $\Lambda_{\CO(\bt\ubw\bt')}\to
\Lambda_{\CO(\bt\ubw\bt')}[-1](-1)$ est nul (il correspond \`a un
\'el\'ement de $\Ext^{-1}(\Lambda_{\CO(\bt\ubw\bt')},
\Lambda_{\CO(\bt\ubw\bt')}(-1))$), d'o\`u
$$R\tilde\pi'_!\Lambda_{\CO(\bt\ubw\ \ubs\bt')}\simeq
 \Lambda_{\CO(\bt\ubw\bt')}[-2](-1)\oplus \Lambda_{\CO(\bt\ubw\bt')}$$
et (ii) s'en d\'eduit par changement de base.
\end{pf*}

Lorsque $sw<w$,
on a  une proposition  analogue pour  les
applications  canoniques $\bX(\bt  \bs\ubw\bt')\to \bX(\bt\ubw\bt')$  et
$\bX(\bt \ubs\ \ubw\bt')\to \bX(\bt\ubw\bt')$.

\sub{}
Soit $a:\bX\to\Spec \Fq$ une vari\'et\'e sur $\Fq$ de dimension $n$.
Suivant Deligne
\cite[\S 3.3.11]{Weil2}, nous dirons que la vari\'et\'e $\bX$ est
rationnellement
lisse  si  le  morphisme $(\Qlbar)_X(n)[2n]\to  Ra^!\Qlbar$  adjoint  du
morphisme  trace est  un isomorphisme.  En particulier,  une vari\'et\'e
lisse est rationnellement lisse. Les conjectures de Weil sur les valeurs
propres de  l'endomorphisme de  Frobenius s'appliquent  aux vari\'et\'es
rationnellement lisses.

La vari\'et\'e $\bX$ est rationnellement
lisse si pour tout point ferm\'e $i_x:\{x\}\to \bX$, le morphisme
canonique
$\Qlbar(n)[2n]\to Ri_x^! \Qlbar$ est un isomorphisme, \ie, si
$H^i_{<x>}(\bX,\Qlbar)=0$ pour $i\not=n$ et
$H^{2n}_{<x>}(\bX,\Qlbar)\simeq\Qlbar(-n)$
\cite[d\'efinition A1.a]{KaLu}.
Cela  implique que  le  faisceau  constant  d\'ecal\'e
$\Qlbar[n]$ est pervers \cite[\S 4.0]{BBD}.

\begin{proposition}\label{X rationnellement lisse}
Soit $w\in W$. Les propri\'et\'es suivantes sont \'equivalentes:
\begin{enumerate}
\item $\CO(\uw)$ est rationnellement lisse.
\item $\bX(\uw)$ est rationnellement lisse.
\item Pour tout $y\leq w$ le polyn\^ome de Kazhdan-Lusztig
$P_{y,w}$ vaut
$1$.
\item Le polyn\^ome de Kazhdan-Lusztig $P_{1,w}$ vaut $1$.
\item La vari\'et\'e de Schubert $\overline{\CB(w)}$ est rationnellement
lisse.
\end{enumerate}
\end{proposition}
\begin{pf*}{Preuve}
Montrons d'abord que (i) et (v) sont  \'equivalents.

Restreignons le diagramme de la preuve du lemme \ref{Owbar}. On obtient
$$\begin{CD}q\inv(\CO(\uw))@>q >>\CO(\uw)\\
@V p VV\\
r\inv(\overline{\CB(w)})@> r>>\overline{\CB(w)}\\
\end{CD}$$
o\`u les morphismes $p,q,r$ sont des fibrations localement
triviales, de fibres respectives
$\bG$, $\bB\times\bB$ et $\bB$. 
Par   cons\'equent   la   lissit\'e  rationnelle   de   $\CO(\uw)$   est
\'equivalente \`a  celle de $\overline{\CB(w)}$ (localement,  on compare
la lissit\'e rationnelle d'une vari\'et\'e  et celle du produit de cette
vari\'et\'e par une vari\'et\'e lisse).

D'apr\`es \cite[th\'eor\`eme A2]{KaLu},
la vari\'et\'e $\overline{\CB(w)}$ est rationnellement lisse au
point $b\in\CB(y)\subset\overline{\CB(w)}$
si et seulement si $P_{y,w}=1$. La
sous-vari\'et\'e ferm\'ee $\CS$ des points non rationnellement lisses
de $\overline{\CB(w)}$ est donc l'union des $\CB(y)$ tels que
$P_{y,w}\ne 1$.
L'\'equivalence de (v), (iii) et (iv)
en r\'esulte, car $\CS$ \'etant ferm\'ee, on a
$\CB(1)\subset\CS$ si et seulement si $\CS$ est non vide.

Les fibrations de la preuve de la proposition \ref{desingularisation}
montrent que
(ii) et (v) sont \'equivalents (\cf\ \cite[preuve du lemme 4.3]{Lusympl}).
\end{pf*}

Si $x$  est une ind\'etermin\'ee  ou un \'el\'ement  de $\Qlbar^\times$,
nous notons  $\CH_x(W)$ l'alg\`ebre de  Hecke de $W$  avec param\`etre
$x$,  \ie, le  quotient  de l'alg\`ebre  du groupe  de  tresses $B$  sur
$\bbZ[x^{1/2},x^{-1/2}]$  par l'id\'eal  bilat\`ere  engendr\'e par  les
$(\bs+1)(\bs-x)$  pour $\bs\in\bS$.  On  note  $T_{\bw}$ (ou  simplement
$T_w$) l'image de $\bw\in\bW$ dans $\CH_x(W)$. On pose
$T_{\uw}=\sum_{w'\le w}T_{w'}$ pour $w\in W$.

\begin{lemme}\label{Twbar}
Le  morphisme canonique
$\bbZ[x^{1/2},x^{-1/2}]B^+\to  \CH_x(W)$  s'\'etend en  un
morphisme d'alg\`ebres  $\bbZ[x^{1/2},x^{-1/2}]\uB^+\to\CH_x(W)$ donn\'e
par $\ubw\mapsto T_{\uw}$ pour $w\in W$.
\end{lemme}
\begin{pf*}{Preuve}
La v\'erification que les  relations qui  d\'efinissent $\uB^+$  sont
v\'erifi\'ees  dans  $\CH_x(W)$ se fait par un calcul imm\'ediat dans
l'alg\`ebre de Hecke.
\end{pf*}

Nous notons $T_\bt$ l'image de $\bt\in\uB^+$ par ce morphisme.

Nous \'etudions maintenant la lissit\'e rationnelle de $\CO(\uw\us)$ sous
l'hypoth\`ese que $\CO(\uw)$ est rationnellement lisse.

\begin{lemme}\label{lemme T}
Soit $w\in W$ tel que $\CO(\uw)$ est rationnellement lisse et soit
$s\in S$.  Alors, une des assertions suivantes est vraie:
\begin{enumerate}
\item $ws<w$. Alors, $T_\uw T_\us=(x+1) T_\uw$.
\item $s\not\le w$. Alors, $T_\uw T_\us=T_{\underline{ws}}$ et
$\CO(\underline{ws})$ est rationnellement lisse.
\item $s<w,ws>w$ et il existe $y\in W$ tel que $ys<y<w$ et $l(y)=l(w)-1$.
Alors, $y$ est l'unique \'el\'ement maximal de $\{v\in W|vs<v<w\}$, les vari\'et\'es
$\CO(\underline{ws})$ et $\CO(\uy)$ sont
rationnellement lisses et $T_\uw T_\us=T_{\underline{ws}}+xT_\uy$.
\item $s<w,ws>w$ et il n'existe aucun $y\in W$ tel que $ys<y<w$ et
$l(y)=l(w)-1$.
Alors, $\CO(\underline{ws})$ n'est pas rationnellement lisse.
\end{enumerate}
\end{lemme}
\begin{pf*}{Preuve}
Remarquons  tout  d'abord  que   $\CO(\uw)$  est  rationnellement  lisse
si  et   seulement  si   $T_\uw=D_w$  (o\`u  $C'_w$   est  la
base   de   Kazhdan-Lusztig   \cite[\S 5.1]{Lubook} et $D_w=x^{l(w)/2}C'_w$).
En   effet,   on
a   $D_w=\sum_{y\le  w}   P_{y,w}T_y$,   et  $\CO(\uw)$   est
rationnellement lisse  si et  seulement si  $P_{y,w}=1$ pour  tout $y\le
w$.  Selon la proposition \ref{X  rationnellement  lisse},
la vari\'et\'e  $\CO(\uw)$
est  rationnellement  lisse  si  et   seulement  si  le  coefficient  de
$D_w$ sur $T_1$ vaut $1$ (si  ce coefficient ne vaut pas $1$,
c'est un polyn\^ome  \`a coefficients positifs, de terme  constant 1 et
de degr\'e strictement  plus grand que 1). La  formule de multiplication
pour les $D_w$ s'\'ecrit (\cf\ \cite[\S 5.1.15]{Lubook}):
$$D_w D_s=
\begin{cases}
(x+1) D_w,&\text{si $ws<w$}\\
 D_{ws}+\sum\limits_{\genfrac{}{}{0pt}{}{y<w, ys<y }{l(y)\not\equiv l(w)\pmod2}}
 \mu(y,w) x^{\frac{l(ws)-l(y)}2}\cdot D_y &\text{sinon.}\\
\end{cases}\\
$$
Le cas $ws<w$ donne l'assertion (i) du lemme.

L'entier $\mu(y,w)$
est le  coefficient de  $x^{(l(w)-l(y)-1)/2}$ dans  $P_{y,w}$. Puisque
$\CO(\uw)$  est   rationnellement  lisse,  on  a   $\mu(y,w)=0$  sauf  si
$l(y)=l(w)-1$ et on a alors $\mu(y,w)=1$.

Supposons $ws>w$ et $\CO(\uw)$ rationnellement lisse.
Alors la formule ci-dessus donne:
$$T_\uw T_\us=D_{ws}+\sum_{\genfrac{}{}{0pt}{}{y<w, ys<y}{l(y)=l(w)-1}}
xD_y.\eqno(1)$$
D'autre part, un calcul direct donne:
$$T_\uw T_\us=T_{\underline{ws}}+x(\sum_{v<w, vs<v}
(T_v+T_{vs})).\eqno(2)$$
Si $s\not\le w$ on a donc $T_\uw T_\us=T_{\underline{ws}}=D_{ws}$
d'o\`u le (ii) du lemme.

Supposons maintenant $s<w$. L'\'egalit\'e (2) montre que le
coefficient  de $T_\uw  T_\us$ sur  $T_1$  vaut $x+1$.  Il en  r\'esulte
que  dans  (1),  il  y   a  au  plus  un  $y$.  S'il  n'existe
aucun  $y$,   alors  $D_{ws}=T_{\underline{ws}}+x(\sum_{v<w,
vs<v}   (T_v+T_{vs}))$   a  un   coefficient   $x+1$   sur  $T_1$   donc
$\CO(\uw\us)$  n'est  pas  rationnellement   lisse  d'o\`u  l'assertion (iv).

S'il  existe   un  $y$,   alors  les   coefficients  de   $T_1$
dans  $D_y$  et   dans  $D_{ws}$  doivent  \^etre
\'egaux  \`a  $1$,   sinon  le  coefficient  de  $T_1$   dans  le  membre
de   droite   de  (1)   exc\'ederait   $x+1$.  Par cons\'equent,
$\CO(\uw\us)$   et
$\CO(\uy)$ sont  rationnellement lisses, donc  $T_\uy=D_y$ et
$T_{\underline{ws}}=D_{\underline{ws}}$. La  comparaison des
\'egalit\'es  (1)  et  (2)  montre  alors  que
$T_{\uy}=\sum_{v<w, vs<v} (T_v+T_{vs})$ et on d\'eduit l'assertion (iii).
\end{pf*}

\sub{}Nous \'etudions maintenant le lien entre vari\'et\'es $\bX$ apr\`es oubli
d'un sous-groupe de Borel.

\begin{proposition}\label{fibrationws}
Soient   $\bt,\bt'\in\uB^+$,   $w\in   W$   et   $s\in   S$   tels   que
$s<w$    et   $ws>w$.    Alors,   on    a   un    triangle   distingu\'e
de     complexes     de    faisceaux     sur     $\bX(\bt\ubw\ubs\bt')$:
$$\Lambda_{\bX(\bt\ubw\ubs\bt')}\to          R\pi_!\Lambda_{\bX(\bt\ubw\
\ubs\bt')}\to        \Lambda_{\bY}[-2](-1)\rightsquigarrow$$        o\`u
$\pi:\bX(\bt\ubw\  \ubs\bt')\to  \bX(\bt\underline{\bw\bs}\bt')$ est  la
projection canonique et $\bY$ est la sous-vari\'et\'e ferm\'ee
$\bigcup_{v\in W,vs<v<w} \bX(\bt\ubv\bt')$
de $\bX(\bt\underline{\bw\bs}\bt')$.

En  outre,  le  faisceau  $\Lambda_\bY$  est  dans  la  sous-cat\'egorie
triangul\'ee pleine  de la cat\'egorie $D^b(\bX(\bt\ubw\ubs\bt'))$ 
engendr\'ee  par les
$\Lambda_{\bX(\bt\ubv\bt')}$, $v\in W$, $vs<v<w$.

On suppose maintenant
que $\CO(\uw)$ et $\CO(\uw\us)$ sont rationnellement lisses.
Alors, il existe un unique \'el\'ement $y\in W$ tel que
$ys<y<w$ et $l(y)=l(w)-1$. On a $\bY=\bX(\bt\uby\bt')$ et le triangle
distingu\'e est scind\'e, \ie,
$$R\pi_!\Lambda_{\bX(\bt\ubw\ \ubs\bt')}\simeq
\Lambda_{\bX(\bt\ubw\ubs\bt')}\oplus
\Lambda_{\bX(\bt\uby\bt')}[-2](-1).$$
\end{proposition}

\begin{pf*}{Preuve}
L'ensemble des \'el\'ements $v\in W$ tels que $v<w$
et $vs<w$ est clos inf\'erieurement pour l'ordre de Bruhat. Soit
$\tilde{\bV}=\coprod_{y\in W, y\le w,ys\not< w}\CO(y,\us)$
l'ouvert de $\CO(\uw,\us)$ compl\'ementaire de
$$\tilde{\bZ}=\coprod_{x\in W,x<w,xs<w}\CO(x,\us)=
\bigcup_{v\in W,vs<v<w}\CO(\uv,\us).$$
Soit maintenant $\bV=\coprod_{y\in W,y\le w,ys\not<w}
\left (\CO(ys)\coprod\CO(y)\right)$
l'ouvert de $\CO(\uw\us)$ comp\-l\'ementaire de
$$\bZ=\coprod_{x\in  W,x<w,xs<w}\CO(x)= \bigcup_{v\in  W, vs<v<w}  \CO
(\uv).$$
 Alors, on prouve par la proposition \ref{omnibus}(v) appliqu\'ee
avec $w'=s$ que
le morphisme $p=(p',p''):\CO(\uw,\us)\to\CO(\uw\us)$ se  restreint en  un
isomorphisme $r:\tilde{\bV}\iso\bV$ et en un morphisme
$q:\tilde{\bZ}\to\bZ$. On a alors
un diagramme commutatif,
o\`u les deux carr\'es sont cart\'esiens et les fl\`eches verticales propres:
$$\xymatrix{
\tilde{\bV}  \ar@{^{(}->}[r]^-{j'}\ar[d]_{r}& \CO(\uw,\us)
 \ar[d]^p& \tilde{\bZ} \ar@{_{(}->}[l]_-{i'}\ar[d]^q \\
\bV\ar@{^{(}->}[r]^-j      & \CO(\uw\us)
          & \bZ \ar@{_{(}->}[l]_-i
}$$
Par changement de base, on a $j^*Rp_*\Lambda_{\CO(\uw,\us)}\simeq
Rr_*j^{\prime *}\Lambda_{\CO(\uw,\us)}\simeq
\Lambda_{\bV}$. Par cons\'equent, $j^*$ \'etant exact, on a
$0=\tau_{\ge 1}j^*Rp_*\Lambda_{\CO(\uw,\us)}=
j^*\tau_{\ge 1}Rp_*\Lambda_{\CO(\uw,\us)}$, \ie,
$\tau_{\ge 1}Rp_*\Lambda_{\CO(\uw,\us)}$ est support\'e par
$\bZ$. Puisque $i^*$ est exact, on a
\begin{align*}
\tau_{\ge 1}Rp_*\Lambda_{\CO(\uw,\us)}&\simeq
i_*i^*\tau_{\ge 1}Rp_*\Lambda_{\CO(\uw,\us)}\simeq
i_*\tau_{\ge 1}i^*Rp_*\Lambda_{\CO(\uw,\us)}\\
&\simeq i_*\tau_{\ge 1}Rq_*i^{\prime *}\Lambda_{\CO(\uw,\us)}\simeq
i_*\tau_{\ge 1}Rq_*\Lambda_{\tilde{\bZ}}.\\
\end{align*}

On v\'erifie, comme dans la preuve de la proposition \ref{wss} (ii),
que $$Rq_*\Lambda_{\tilde{\bZ}}\simeq
\Lambda_{\bZ}\oplus\Lambda_{\bZ}[-2](-1).$$
Par cons\'equent, on a $\tau_{\ge 1}Rq_*\Lambda_{\tilde{\bZ}}\simeq
\Lambda_{\bZ}[-2](-1)$.
Puisque les fibres de $p$ sont connexes, on a
$p_*\Lambda_{\CO(\uw,\us)}\simeq \Lambda_{\CO(\uw\us)}$
et le triangle distingu\'e
$$p_*\Lambda_{\CO(\uw,\us)}\to
Rp_*\Lambda_{\CO(\uw,\us)}\to
\tau_{\ge 1}Rp_*\Lambda_{\CO(\uw,\us)}\rightsquigarrow$$
fournit le triangle distingu\'e
$$\Lambda_{\CO(\uw\us)}\to Rp_*\Lambda_{\CO(\uw,\us)}\to
\Lambda_{\bZ}[-2](-1)\rightsquigarrow$$
Comme
$\bY=(\CO(\bt)\times_\CB\bZ\times_\CB\CO(\bt'))\times_{\CB\times\CB}\Gamma$,
nous d\'eduisons par changement de base le triangle distingu\'e
voulu pour les vari\'et\'es de Deligne-Lusztig.

Montrons maintenant l'affirmation sur $\Lambda_\bY$.
Soient $v,v'\in W$ distincts avec $vs<v<w$ et $v's<v'<w$.
L'ensemble $\{x\in W|x<v,x<v'\}$ est r\'eunion
des intervalles $I_z=\{x\in W, x\le z\}$, o\`u $zs<z$, $z<v$ et $z<v'$.
Par cons\'equent,
$\bX(\bt\ubv\bt')\cap\bX(\bt\ubv'\bt')=\bigcup_z \bX(\bt\ubz\bt')$ o\`u
$z$ d\'ecrit les \'el\'ements de $W$ tels que $z<v$, $z<v'$ et $zs<z$.
On en d\'eduit
le r\'esultat sur $\Lambda_\bY$ par application it\'er\'ee de la suite exacte
de Mayer-Vietoris.

Nous allons maintenant utiliser le lemme \ref{lemme T} pour \'etablir la
derni\`ere partie de la proposition.
Puisque $\CO(\uw)$ et $\CO(\underline{ws})$ sont
rationnellement lisses,
l'\'egalit\'e du lemme \ref{lemme T} (iii)  peut se r\'e-\'ecrire
$(T_\uw - T_\uy) T_\us= T_{\underline{ws}}-T_\uy$, car
$T_\uy T_\us=(x+1) T_\uy$ d'apr\`es le lemme \ref{lemme T} (i).
Dans la preuve du lemme \ref{lemme T} (iii), nous avons montr\'e que
$T_\uy=\sum_{v<w,vs<v}(T_{v}+T_{vs})$, ce \'equivaut \`a
l'\'egalit\'e: $\{v\mid vs<v<w\}=\{v\mid  vs<v\le y\}$.
On a donc $\bZ=\CO(\uy)$ et on a un triangle distingu\'e
$$\Lambda_{\CO(\uw\us)}[d]\to Rp_*\Lambda_{\CO(\uw,\us)}[d]\to
\Lambda_{\CO(\uy)}[d-2](-1)\rightsquigarrow.$$
La vari\'et\'e $\CO(\uw\us)$ est rationnellement lisse, donc
$\Lambda_{\CO(\uw\us)}[d]$ est un faisceau pervers pur de poids $d$ sur
$\CO(\uw\us)$, o\`u $d=l(w)+1$.
Puisque la vari\'et\'e $\CO(\uy)$ est rationnellement lisse
(lemme \ref{lemme T} (iii)),
le complexe $\Lambda_{\CO(\uy)}[d-2](-1)$ est un faisceau
pervers pur de poids $d$.
On a donc une suite
exacte de faisceaux pervers
$$0\to \Lambda_{\CO(\uw\us)}[d]\to
Rp_*\Lambda_{\CO(\uw,\us)}[d]
\to \Lambda_{\CO(\uy)}[d-2](-1)\to 0.$$
Puisque $\Lambda_{\CO(\uw\us)}[d]$ et
$\Lambda_{\CO(\uy)}[d-2](-1)$ sont purs de m\^eme poids,
le faisceau pervers $Rp_!\Lambda_{\CO(\uw,\us)}[d]$ est pur, donc la suite
exacte se scinde d'apr\`es le th\'eor\`eme de d\'ecomposition
\cite[th\'eor\`eme 5.3.8]{BBD}, \ie,
$$Rp_*\Lambda_{\CO(\uw,\us)}[d]\simeq
\Lambda_{\CO(\uw\us)}[d]\oplus
\Lambda_{\CO(\uy)}[d-2](-1).$$
Par changement de base, on d\'eduit une d\'ecomposition analogue au niveau
des vari\'et\'es de Deligne-Lusztig.
\end{pf*}

Nous g\'en\'eralisons maintenant \cite[theorem 1.6 (case 2)]{DL}. Noter que ce
r\'esultat, contrairement aux pr\'ec\'edents,
ne provient pas directement par changement de base d'un r\'esultat
sur les vari\'et\'es $\CO$.
Nous remercions George Lusztig de nous avoir communiqu\'e le principe de la
d\'emonstration de la remarque de \cite[\S 1.6.3]{DL} que les morphismes
\og$\partial$\fg\  sont nuls.

\begin{proposition}\label{sws} Soient $\bs\in\bS$ et
$\bb\in\BW$. Soit $\rho:\bX(\bs\bs\bb)\to\bX(\ubs \bb)$ d\'eduit
par changement de base de $(p',p''):\CO(s,s)\to\CO(\us)$.
Alors, on a un triangle distingu\'e de complexes de faisceaux
sur $\bX(\ubs \bb)$:
$$\Lambda_{\bX(\bs\bb)}[-2](-1)\oplus \Lambda_{\bX(\bs\bb)}[-1]\to
R\rho_!\Lambda_{\bX(\bs\bs\bb)}\to \Lambda_{\bX(\bb)}[-2](-1)\rightsquigarrow$$
En particulier, on a
une suite exacte longue de $(\GF\times\genby{F^\delta})$-modules:
\begin{align*}\cdots&\to H^{i-3}_c(\bX(\bb))(-1)\to
H^{i-2}_c(\bX(\bs\bb))(-1)\oplus H^{i-1}_c(\bX(\bs\bb))
\to H^i_c(\bX(\bs\bs\bb))\\
&\to H^{i-2}_c(\bX(\bb))(-1)\to\cdots\\
\end{align*}
\end{proposition}

\begin{pf*}{Preuve}
Choisissons une d\'ecomposition $\bb=\bw_1\cdots\bw_k$ avec $w_i\in W$.
Soit $\bX_1$ la sous-vari\'et\'e ferm\'ee 
$\{(\bB_1,\bB_2,\bB_3,\ldots,\bB_{k+2})\in\bX(s,s,w_1,\ldots,w_k)
\mid \bB_1=\bB_3\}$ et soit $\bX_2$ l'ouvert compl\'ementaire dans
$\bX(s,s,w_1,\ldots,w_k)$.
L'application $$(\bB_1,\bB_2,\bB_3,\bB_4,\ldots,\bB_{k+2})\mapsto
(\bB_1,\bB_4,\ldots,\bB_{k+2})$$
fait de $\bX_1$ un fibr\'e en droites affines au-dessus de
$\bX(w_1,\ldots,w_k)$,
donc $H^i_c(\bX_1)\simeq H^{i-2}_c(\bX(\bb))(-1)$.

La vari\'et\'e $\bX_2$ est un ouvert de la vari\'et\'e 
$$\bY=\{(\bB_1,\ldots,\bB_{k+2})\mid
(\bB_1,\bB_3,\ldots,\bB_{k+2})\in\bX(s,w_1,\ldots,w_k),
\bB_1\xrightarrow s\bB_2\}.
$$
La restriction $\pi:\bY\to\bX(s,w_1,\ldots,w_k)$ de $\rho$
est une fibration en droites affines et $\pi$ restreint \`a $\bY-\bX_2$
est un isomorphisme.
On d\'eduit alors de la suite exacte
$0\to \Lambda_{\bX_2}\to\Lambda_\bY\to
\Lambda_{\bY-\bX_2}\to 0$ un triangle distingu\'e
$$R\pi_!\Lambda_{\bX_2}\to\Lambda_{\bX(\bs\bb)}[-2](-1)
\xrightarrow\partial \Lambda_{\bX(\bs\bb)}\rightsquigarrow$$
o\`u
$$\partial\in\Hom(\Lambda_{\bX(\bs\bb)}[-2](-1),
\Lambda_{\bX(\bs\bb)})=\Ext^2(\Lambda_{\bX(\bs\bb)},
\Lambda_{\bX(\bs\bb)}(1))=H^2(\bX(\bs\bb))(1).$$
Toutes les fl\`eches \'etant $\GF$-\'equivariantes, on a
$\partial\in H^2(\bX(\bs\bb))^\GF$.
Puisque $\bs\bb\in B^+$, on d\'eduit de la proposition \ref{Id} \`a venir (mais
dont la preuve n'utilise que des r\'esultats d\'ej\`a prouv\'es) que
$H^{2l(\bs\bb)-2}_c(\bX(\bs\bb))^\GF=0$. La vari\'et\'e
 $\bX(\bs\bb)$ est lisse, donc on obtient par dualit\'e
$H^2(\bX(\bs\bb))^\GF=0$, donc $\partial=0$. Par cons\'equent, le triangle
distingu\'e donne un isomorphisme
$R\pi_!\Lambda_{\bX_2}\simeq\Lambda_{\bX(\bs\bb)}[-2](-1)
\oplus\Lambda_{\bX(\bs\bb)}[-1]$ d'o\`u
$H^i_c(\bX_2)\simeq H^{i-2}_c(\bX(\bs\bb))(-1)\oplus
H^{i-1}_c(\bX(\bs\bb))$.
On conclut en utilisant la suite exacte associ\'ee \`a la
d\'ecomposition ouvert-ferm\'e $\bX(\bs\bs\bb)=\bX_2\coprod\bX_1$.
\end{pf*}

\subsection{Repr\'esentations unipotentes}
\label{sectionunipotentes}
\sub{}
Nous allons g\'en\'eraliser \cite[th\'eor\`eme 3.8]{LuMa}.

On dira que $w,w'\in W$ sont {\it $F$-conjugu\'es} s'il existe $v\in W$
tel que $w'=v w F(v)^{-1}$.

\begin{proposition}\label{XtxXtp} Pour $\bb,\bb'\in\BW$,  on a
\begin{enumerate}
\item $\left(H^i_c(\bX(\bb))\otimes H^{i'}_c(\bX(\bb'))\right)^\GF=0$
pour $i+i'<l(\bb)+l(\bb')$.
\item Les valeurs propres de $F^\delta$ sur
$\left(H^i_c(\bX(\bb))\otimes H^{i'}_c(\bX(\bb'))\right)^\GF$
sont des puissances enti\`eres de $q^\delta$.
\item Si $i+i'=l(\bb)+l(\bb')$, les valeurs propres de $F^\delta$ sur
$$\left(H^i_c(\bX(\bb))\otimes H^{i'}_c(\bX(\bb'))\right)^\GF$$
sont des puissances de $q^\delta$ inf\'erieures \`a
$q^{\delta(l(\bb)+l(\bb'))/2}$ et cette valeur ne peut  \^etre atteinte
que si $\beta(\bb)$ est $F$-conjugu\'e \`a $\beta(\bb')$.  
\end{enumerate}
\end{proposition}

\begin{pf*}{Preuve} Par le th\'eor\`eme de K\"unneth,
$$\left(H^i_c(\bX(\bb))\otimes H^{i'}_c(\bX(\bb'))\right)^\GF$$ est un facteur direct de
$H^{i+i'}_c(\GF\backslash(\bX(\bb)\times\bX(\bb')))$. On est donc ramen\'e \`a prouver les propri\'et\'es
ci-dessus pour cette cohomologie.

Soient $(x_1,\ldots,x_k)$ et
$(x'_1,\ldots,x'_k)$ deux suites d'\'el\'ements de $W$ telles que
$\bb=\bx_1\cdots\bx_k$ et $\bb'=\bx'_1\cdots\bx'_k$ respectivement.
Identifions $\bX(\bb)$ et $\bX(\bb')$, comme dans la proposition
\ref{produit}, aux
vari\'et\'es de Deligne-Lusztig ordinaires du groupe $\bG^k$, pour
l'isog\'enie $F_1$, associ\'ees aux \'el\'ements $(x_1,\ldots,x_k)$ et
$(x'_1,\ldots,x'_k)$ respectivement. La vari\'et\'e
$\GF\backslash(\bX(\bb)\times\bX(\bb'))$ est alors identifi\'ee \`a la
vari\'et\'e
$(\bG^k)^{F_1}\backslash(\bX(x_1,\ldots,x_k)\times\bX(x'_1,\ldots,x'_k))$.
D'autre part si les \'el\'ements $(x_1,\ldots,x_k)$ et
$(x'_1,\ldots,x'_k)$ sont $F_1$-conjugu\'es par $(v_1,\ldots,v_k)$,
alors les \'el\'ements $\beta(\bb)=x_1\cdots x_k$ et $\beta(\bb')=x'_1\cdots
x'_k$ sont $F$-conjugu\'es par $v_1$.   On est donc ramen\'e \`a
d\'emontrer la proposition pour des vari\'et\'es de Deligne-Lusztig
ordinaires. Nous supposons donc dans la suite que $\bb$ et $\bb'$ sont
dans $\bW$.

Nous pouvons alors appliquer les r\'esultats de Lusztig
\cite[proposition 3.4 et lemme 3.5]{LuMa}: la vari\'et\'e
$\GF\backslash(\bX(\bb)\times\bX(\bb'))$ admet une stratification par
l'union, quand $w_1$ parcourt $W$, des vari\'et\'es
$$\bY(\bb,\bb',w_1)=\GF\backslash
\{(\bB_1,\bB'_1)\in\bX(\bb)\times\bX(\bb')
\mid \bB_1\xrightarrow{w_1}\bB'_1\}.$$
On note $\bZ(\bv,w_1)$ la fibre de
$(\bB,\lexp{w_1}\bB)$ par $(p',p''):\CO(\bv)\to\CB\times\CB$, pour
$\bv\in\BW$ et $w_1\in W$. Si
$\bv=\bs_1\cdots\bs_k$ est une d\'ecomposition avec $s_i\in S$, alors
on a un isomorphisme
$$\bZ(\bv,w_1)\simeq\{(\bB_0,\ldots,\bB_k)\in\CB^{k+1}\mid
\bB_{i-1}\xrightarrow{s_i}\bB_i, (\bB_0,\bB_k)=(\bB,\lexp{w_1}\bB)\}.$$

Prenons $\bv=\bb F(\bw_1)\tilde\bb'$
o\`u $\tilde{\bb'}$ est le
retourn\'e du mot $\bb'$. Alors, on a un isomorphisme
de $F^\delta$-modules (\cf\ \cite[lemme 3.5]{LuMa} o\`u la
vari\'et\'e $\bZ(\bv,w_1)$ est not\'ee $\bZ_{b,b',w_1}$)
$$H^i_c(\bY(\bb,\bb',w_1))\simeq
H^i_c(\bZ(\bb F(\bw_1)\tilde\bb',w_1)).$$
La d\'emonstration donn\'ee par Lusztig s'\'etend du cas o\`u $F$ est
un endomorphisme de Frobenius au cas o\`u $F$ est une isog\'enie dont une
puissance est un endomorphisme de Frobenius, la seule propri\'et\'e de
$F$ utilis\'ee \'etant le th\'eor\`eme de Lang.

\'Enon\c cons maintenant le lemme clef:

\begin{lemme}\label{Zvw}\begin{enumerate}
\item On a $H^i_c(\bZ(\bv,w_1))=0$ pour
$i<l(\bv)-l(w_1)$.
\item Les valeurs propres de $F^\delta$ sur
$H^k_c(\bZ(\bv,w_1))$ sont des puissances de $q^\delta$.
\item De plus si $k= l(\bv)-l(w_1)$, les valeurs propres de $F^\delta$ sont
inf\'erieures \`a $q^{\delta(l(\bv)-l(w_1))/2}$, cette valeur ne pouvant
\^etre atteinte que si $\beta(\bv)=w_1$.
\end{enumerate}
\end{lemme}

\begin{pf*}{Preuve} Nous allons d\'emontrer le
r\'esultat par r\'ecurrence sur $l(\bv)$.
Si la vari\'et\'e $\bZ(\bv,w_1)$ est non vide, alors une d\'ecomposition r\'eduite
de $w_1$ peut \^etre extraite de $\bv$.
Cette vari\'et\'e est donc vide si
$l(\bv)<l(w_1)$ ou si $l(\bv)=l(w_1)$ et $\beta(\bv)\ne w_1$. Si
$l(\bv)=l(w_1)$ et $\beta(\bv)=w_1$ la vari\'et\'e est r\'eduite \`a un
point: ceci d\'emontre le lemme pour $l(\bv)=l(w_1)$.
C'est le point de d\'epart de la r\'ecurrence.

Supposons $l(\bv)>l(w_1)$. Si $\bv\in\bW$, alors la vari\'et\'e est vide et
l'\'enonc\'e est vrai. Nous pouvons donc supposer que $\bv\notin\bW$,
\ie\ que $\bv=\ba\bt\bt\bb$ o\`u $\bt\in \bS$ et $\ba,\bb\in B^+$.
Soient $\ba=\bs_1\cdots\bs_i$ et $\bb=\bs_{i+3}\cdots\bs_k$ des
d\'ecompositions avec $s_i\in S$. On a alors la d\'ecomposition
$\bv=\bs_1\cdots\bs_i\bt\bt\bs_{i+3}\cdots\bs_k$.
Soit $\bF$ le ferm\'e de $\bZ(\bv,w_1)$ form\'e des
$$(\bB_0,\ldots,\bB_i,\bB_{i+1},\bB_{i+2},\ldots,\bB_k)$$ tels que
$\bB_i=\bB_{i+2}$ et soit $\bO$ l'ouvert compl\'ementaire.
Alors, $\bF$ est un fibr\'e en droites sur $\bZ(\ba\bb,w_1)$, donc
$H^i_c(\bF)\simeq H^{i-2}_c(\bZ(\ba\bb,w_1))(-1)$. Quant \`a $\bO$, il est
isomorphe au compl\'ement de la section nulle d'un fibr\'e  en
droites sur $\bZ(\ba\bt\bb,w_1)$. On a une suite exacte longue de
cohomologie associ\'ee \`a la fibration
\begin{multline*}\cdots \to H^{i-2}_c(\bZ(\ba\bt\bb,w_1))(-1)\to
H^i_c(\bZ(\ba\bt\bb,w_1))\to\\ \to H^{i+1}_c(\bO)
\to H^{i-1}_c(\bZ(\ba\bt\bb,w_1))(-1)\to\cdots
\end{multline*}
L'hypoth\`ese de r\'ecurrence montre alors que
$H^i_c(\bF)=0$ et $H^i_c(\bO)=0$ pour $i<l(\bv)-l(w_1)$ et que
$H^{l(\bv)-l(w_1)}_c(\bF)\simeq H^{l(\ba\bb)-l(w_1)}_c(\bZ(\ba\bb,w_1))(-1)$
et $H^{l(\bv)-l(w_1)}_c(\bO)\simeq
H^{l(\ba\bt\bb)-l(w_1)}_c(\bZ(\ba\bt\bb,w_1))$.

La suite exacte longue associ\'ee \`a la d\'ecomposition
\og ouvert-ferm\'e\fg
$$\cdots\to H^i_c(\bO)\to H^i_c(\bZ(\bv,w_1))\to H^i_c(\bF)\to
H^{i+1}_c(\bO)\to\cdots$$
montre que $H^i_c(\bZ(\bv,w_1))=0$ pour $i<l(\bv)-l(w_1)$, et
que pour tout $i$, les valeurs propres de $F^\delta$ sur $H^i_c(\bZ(\bv,w_1))$
sont des puissances de $q^\delta$. Elle montre aussi
que les valeurs propres de $F^\delta$ sur $H^{l(\bv)-l(w_1)}_c(\bZ(\bv,w_1))$
sont des puissances de $q^\delta$ inf\'erieures \`a
$q^{\delta(l(\bv)-l(w_1))/2}$.  De plus, cette valeur ne peut \^etre
atteinte que si $F^\delta$ a une valeur propre
$q^{\delta(l(\ba\bb)-l(w_1))/2}$ sur
$H^{l(\ba\bb)-l(w_1)}_c(\bZ(\ba\bb,w_1))$.  Par l'hypoth\`ese de
r\'ecurrence, il faut donc que $w_1=\beta(\ba\bb)=\beta(\ba\bt\bt\bb)=\beta(\bv)$,
d'o\`u le lemme.  \end{pf*}

De ce lemme nous d\'eduisons maintenant facilement la proposition:
aucune des strates $\bY(\bb,\bb',w_1)$ n'ayant de cohomologie en degr\'e
strictement inf\'erieur \`a
$l(\bb F(\bw_1)\tilde\bb')-l(w_1)=l(\bb)+l(\bb')$,
les suites exactes longues \og ouvert-ferm\'e\fg\ donn\'ees par la stratification
montrent qu'il en est de m\^eme pour $\GF\backslash(\bX(\bb)\times\bX(\bb'))$;
de m\^eme, les valeurs propres de $F^\delta$ sur la cohomologie de toute strate \'etant des puissances de $q^\delta$, 
les suites exactes longues \og ouvert-ferm\'e\fg\ donnent la m\^eme propri\'et\'e pour
$H^k_c(\GF\backslash(\bX(\bb)\times\bX(\bb')))$. Elles montrent aussi que
$H^{l(\bb)+l(\bb')}_c(\GF\backslash(\bX(\bb)\times\bX(\bb')))$
est isomorphe \`a un sous-espace de $\bigoplus_{w_1}
H^{l(\bb)+l(\bb')}_c(\bZ(\bb F(\bw_1)\tilde\bb',w_1))$.
Les valeurs propres de
$F^\delta$ sur $H^{l(\bb)+l(\bb')}_c(\GF\backslash(\bX(\bb)\times\bX(\bb')))$
sont donc des puissances de $q^\delta$ inf\'erieures \`a
$q^{\delta(l(\bb)+l(\bb'))/2}$ et cette valeur ne
peut \^etre atteinte que s'il existe $w_1$ tel que
$\beta(\bb F(\bw_1)\tilde\bb')=w_1$, \ie\ si $\beta(\bb)$ est $F$-conjugu\'e \`a
$\beta(\bb')$.\end{pf*}

Le  corollaire suivant  est  une g\'en\'eralisation  des r\'esultats  de
\cite[corollaire 3.9]{LuMa} et de \cite[III, th\'eor\`eme 2.3]{DM}.

\begin{corollaire}\label{valeurs propres}
\begin{enumerate}
\item
Les   valeurs   propres   de   $F^\delta$   sur   un   sous-$\GF$-module
irr\'educ\-tible   $\rho$   de  $H^i_c(\bX(\bt))$   (pour   $\bt\in\uB^+$)
sont dans $q^{\delta\bbN}\lambda_\rho\omega_\rho$, o\`u  $\lambda_\rho$  est
une  racine  de   l'unit\'e  et
$\omega_\rho\in\{1,q^{\delta/2}\}$, tous deux ind\'ependants  de   $i$  et
de  $\bt$.
\item
Pour tous $\bb$ et $\bb'$ dans $\uB^+$, les
valeurs propres de $F^\delta$ sur
$\left(H^i_c(\bX(\bb))\otimes H^{i'}_c(\bX(\bb'))\right)^\GF$
sont des puissances enti\`eres de $q^\delta$.
\end{enumerate}
\end{corollaire}

\begin{pf*}{Preuve}
Montrons d'abord la propri\'et\'e (i) pour $\bt\in B^+$.
Dans ce cas, nous pouvons suivre telle
quelle la preuve de \cite[corollaire 3.9]{LuMa} en rempla\c cant
\cite[th\'eor\`eme 3.8]{LuMa} par le (ii) de
notre proposition \ref{XtxXtp} qui en est la
g\'en\'eralisation pour les \'el\'ements du mono{\"\i}de de tresses.
On obtient alors que les valeurs propres de $F^\delta$ sur les
sous-$\GF$-modules irr\'eductibles de m\^eme
type dans les $H^i_c(\bX_\bb)$ diff\`erent par des puissances de $q^\delta$.
On voit que ces valeurs propres sont des racines de l'unit\'e \`a une
puissance de $q^{\delta/2}$ pr\`es en prenant un $\bb\in\bW$ dans la
cohomologie duquel $\rho$ intervient: dans ce cas on a une vari\'et\'e de
Deligne-Lusztig ordinaire pour laquelle on conna{\^\i}t le r\'esultat.

En g\'en\'eral, on raisonne par r\'ecurrence
sur le nombre d'\'el\'ements de $\ubW$ dans une d\'ecomposition de $\bt$.
Si $\bt=\bt\ubw\bt'$ avec $\ubw\in\ubW$ et $\bt,\bt'\in\uB^+$,
alors dans la stratification de la proposition \ref{ouvertferme},
les vari\'et\'es $\bX_i-\bX_{i-1}$
v\'erifient la conclusion du th\'eor\`eme (par r\'ecurrence). En consid\'erant les
suites exactes longues de cohomologie associ\'ees \`a cette stratification,
on voit qu'il en est de m\^eme pour $\bX(\bt)$, d'o\`u (i).

Montrons (ii).
Par (i), les valeurs propres de $F^\delta$ dans
$$\left(H^i_c(\bX(\bb))\otimes H^{i'}_c(\bX(\bb'))\right)^\GF$$
sont de la forme $\lambda_\rho\lambda_{\rho'}\omega_\rho\omega_{\rho'} q^{k\delta}$ avec
$k\in\bbN$, o\`u $\rho$ et $\rho'$ sont des repr\'esentations contragr\'edientes l'une de
l'autre. Si on prend $\bb$ et $\bb'$ dans $B^+$ on sait qu'on doit avoir une puissance
enti\`ere de $q^\delta$. Donc $\omega_\rho=\omega_{\rho'}$ et
$\lambda_\rho=\lambda_{\rho'}\inv$. On en d\'eduit le r\'esultat.
\end{pf*}

\sub{}
\label{convention}
Pour  $H$   un  groupe,  on   note  $\Irr(H)$  l'ensemble   des  classes
d'isomorphisme de  $\Qlbar H$-modules  simples de dimension  finie. Nous
notons $\CR(H)$ le groupe de  Grothendieck de la cat\'egorie des $\Qlbar
H$-modules de dimension finie et $\CR^+(H)$ le sous-mono{\"\i}de donn\'e
par les  classes des modules. On  note $[V]$ l'image dans  $\CR(H)$ d'un
$\Qlbar H$-module $V$. Par abus  de notation, on note
simplement $\rho$ l'image d'un \'el\'ement $\rho\in\Irr(H)$.
On pose $V_\rho=\scal V\rho\rho$.

Soit  $L$   le  sous-groupe  du groupe de Grothendieck
des $\Qlbar(\GF\times\genby {F^\delta})$-modules gradu\'es
de dimension finie,
engendr\'e  par  les  classes  des repr\'esentations  o\`u  les  valeurs
propres  de   $F^\delta$  sont   de  la   forme  $\zeta   q^{\delta  j}$
avec   $\zeta$   une   racine   de   l'unit\'e   et   $j\in\bbN/2$. 
On   d\'efinit   un
morphisme $L\to\bbZ[t^{1/2},h]\otimes\CR(\GF)$  comme  suit.
Soit   $\rho$  une   repr\'esentation  irr\'eductible  de   $\GF$  et
$\lambda$  une  repr\'esentation  de  $\genby  {F^\delta}$  donn\'ee  par
$\lambda(F^\delta)=\zeta  q^{\delta  j}$.
Alors,   on  envoie   la  classe   de  la
repr\'esentation  $\rho\otimes\lambda$,  plac\'ee  en degr\'e  $i$,  sur
$h^it^j\rho$.

D'apr\`es le corollaire \ref{valeurs propres}, la classe
$[H^*_c(\bX(\bt))]$
de $H^*_c(\bX(\bt))$ est dans $L$ et le sous-groupe
de $L$ engendr\'e par les classes de telles repr\'esentations s'envoie
injectivement dans $\bbZ[t^{1/2},h]\otimes\CR(\GF)$.
Dans la suite du texte nous identifierons
$[H^*_c(\bX(\bt))]$ \`a un \'el\'ement de $\bbZ[t^{1/2},h]\otimes\CR(\GF)$,
et pour $\rho\in\Irr(\GF)$, nous consid\'ererons $\scal{H^*_c(\bX(\bt))}\rho$
comme un \'el\'ement de $\bbZ[t^{1/2},h]$.
Dans de nombreux cas, nous pourrons d\'ecrire cet \'el\'ement en termes de
caract\`eres de l'alg\`ebre de Hecke.

On    sait    que    $\CH_x(W)\otimes\bbQ(x^{1/2})$    est    absolument
semi-simple  \cite[th\'eor\`eme 9.3.5]{GP2}
et   qu'apr\`es   sp\'ecialisation   des   scalaires   \`a
$\bbQ$   \`a   travers   l'homomorphisme   $f:\sqrt{x}\mapsto   1$,   on
a   $\CH_x(W)\otimes_f\bbQ\iso\bbQ[W]$.  On   a  donc   une  bijection
(\og sp\'ecialisation\fg) de l'ensemble des caract\`eres irr\'eductibles de
$\CH_x(W)\otimes\bbQ(x^{1/2})$ vers  l'ensemble des  caract\`eres
irr\'eductibles de  $W$. Si
$\chi$ est  un caract\`ere de  $W$, on  note $\chi_x$ le  caract\`ere de
$\CH_x(W)$ qui se sp\'ecialise  en $\chi$; on a $\chi(w)=f(\chi_x(T_w))$
pour tout $w\in W$.

Pour $\psi\in\CR(W\sdp F)\otimes\bbQ$, on note $R_\psi$ le
\og caract\`ere fant\^ome\fg,
\'el\'ement de $\CR(\GF)\otimes \Qlbar$ donn\'e par
$$R_\psi=|W|\inv\sum_{i\ge 0,w\in W}\psi(wF)(-1)^i [\Res_\GF H^i_c(\bX(w))]$$
(ici la notation $\Res_\GF$ signifie que nous oublions l'action de
$\genby{F^\delta}$).
On note $\CE(\GF,1)$ l'ensemble des caract\`eres unipotents 
irr\'eductibles de $\GF$.

Pour tout caract\`ere $F$-stable $\chi$ de $W$, nous choisissons
une extension $\tilde\chi$
de $\chi$ \`a $W\sdp F$ triviale sur $F^\delta$ et rationnelle.
Si $W$ est irr\'eductible, il y a une unique telle extension rationnelle,
sauf si $\delta=2$, auquel cas une telle extension est unique
au signe pr\`es \cite[proposition 3.2]{Lubook}.
Soit $m$ un entier multiple de $\delta$.
Nous notons $\CH_{q^m}(W)$ l'alg\`ebre de
Hecke sp\'ecialis\'ee en $\sqrt{x}\mapsto q^{m/2}\in\bbR$.
Nous notons
$\tilde\chi_{q^m}$ le caract\`ere de $\CH_{q^m}(W)\sdp F$ qui se
sp\'ecialise en $\tilde\chi$
(c'est une extension du caract\`ere $\chi_{q^m}$
de $\CH_{q^m}(W)$).

\sub{}Nous allons maintenant \'etendre et g\'en\'eraliser  des r\'esultats de
\cite{DM}.

\begin{proposition}\label{Nmw} Pour tous $\bt\in\uB^+$, $g\in \GF$
et pour tout entier positif $m$ multiple de $\delta$, le nombre
de points fixes sous $gF^m$ de $\bX(\bt)$ est
$$|\bX(\bt)^{gF^m}|=\sum_{\rho\in\CE(\GF,1)}\lambda_\rho^{m/\delta}\rho(g)
\sum_{\chi\in\Irr(W)^F}\scal\rho{R_{\tilde\chi}}
\cdot\tilde\chi_{q^m}(T_{\bt} F).$$
\end{proposition}

\begin{pf*}{Preuve}
Fixons une d\'ecomposition $\bt=\bt_1\cdots\bt_k$ avec $t_i\in W\cup\uW$.
Comme en (\ref{lt ferme}), on a
$$\bX(t_1,\ldots,t_k)=\coprod_{(t'_1,\ldots,t'_k)\subset(t_1,\ldots,t_k)}
\bX(t'_1,\ldots,t'_k)$$ et
$T_{\bt}=\sum_{(t'_1,\ldots,t'_k)\subset(t_1,\ldots,t_k),t'_i\in W}
T_{\bt'_1\cdots\bt'_k}$.
Les deux membres de la proposition \'etant additifs, il suffit donc de
d\'emontrer la proposition pour $\bt\in B^+$.
La propri\'et\'e r\'esulte alors de
\cite[III proposition 1.2, th\'eor\`eme 1.3 et th\'eor\`eme 2.3]{DM} appliqu\'es \`a la vari\'et\'e
$\bX_{\bG^k}(t_1,\ldots,t_k,F_1)$ isomorphe \`a $\bX(\bt)$
(\cf\ proposition \ref{produit}).
\end{pf*}

\begin{corollaire}\label{rationnellement lisse} Soit $\bt\in\uB^+$. On a
\begin{enumerate} 
\item 
\begin{multline*}[H^*_c(\bX(\bt))]_{h=-1} =
\sum_{\chi\in\Irr(W)^F}R_{\tilde\chi}\cdot\tilde\chi_{t}(T_\bt F)
\\=\sum_{\rho\in\CE(\GF,1)}\rho
\sum_{\chi\in\Irr(W)^F}\scal\rho{R_{\tilde\chi}}\cdot
\tilde\chi_{t}(T_\bt F).
\end{multline*}

\item Si de plus $\bt$ est
un produit d'\'el\'ements $\ubw$ de $\ubW$ tels que $\bX(\uw)$ est
rationnellement lisse, alors,
\begin{multline*}[H^*_c(\bX(\bt))]=
\sum_{\chi\in\Irr(W)^F}R_{\tilde\chi}\cdot\tilde\chi_{h^2t}(T_\bt F)
\\=\sum_{\rho\in\CE(\GF,1)}\rho
\sum_{\chi\in\Irr(W)^F}\scal\rho{R_{\tilde\chi}}\cdot
\tilde\chi_{h^2t}(T_\bt F)
\end{multline*}
\end{enumerate}
\end{corollaire}

\begin{pf*}{Preuve}
D'apr\`es le corollaire \ref{valeurs propres} et la convention \ref{convention},
si $m$ est tel que
$\lambda_\rho^{m/\delta}=1$  pour  tout  $\rho$, 
on a 
$$\sum_i(-1)^i\Trace(gF^m|H^i_c(\bX(\bt)))=
\Trace(g|[H^*_c(\bX(\bt))]_{h=-1,t=q^m}).$$
Par la formule des traces de Lefschetz, le premier membre est \'egal \`a
$|\bX(\bt)^{gF^m}|$.
On obtient donc pour, tout $m$ suffisamment divisible,
$$[H^*_c(\bX(\bt))]_{h=-1,t=q^m}=
\sum_{\rho\in\CE(\GF,1)}\rho
\sum_{\chi\in\Irr(W)^F}\scal\rho{R_{\tilde\chi}}
\cdot\tilde\chi_{q^m}(T_{\bt} F),$$
ce qui montre (i), les deux membres \'etant des polyn\^omes \'egaux pour une infinit\'e de
valeurs de $t$.

D\'emontrons (ii).
La proposition \ref{produit} permet de nous
ramener  au cas  o\`u la  suite $\bt$  est r\'eduite  \`a un  seul terme
$\uw$.
Puisque $\bX(\uw)$  est rationnellement  lisse, les  conjectures
de  Weil  sur  les  valeurs  propres  de  l'endomorphisme  de  Frobenius
montrent que
$$\Trace(gF^m|H^i_c(\bX(\bt)))=
\sum_\rho \lambda_\rho^{m/\delta} q^{im/2} \rho(g)\scal\rho{H^i_c(\bX(\bt))}.
$$
Jointe \`a la proposition \ref{Nmw}, la formule des traces de Lefschetz fournit
alors l'\'egalit\'e
$$\sum_i (-1)^i q^{im/2} \scal\rho{H^i_c(\bX(\bt))}=\sum_\chi
\scal\rho{R_{\tilde\chi}}\tilde\chi_{q^m}(T_{\bt} F)$$
pour tout $m$. On en d\'eduit (ii).
La proposition s'en d\'eduit.
\end{pf*}

\sub{}Nous \'etudions maintenant les propri\'et\'es de rationalit\'e des
caract\`eres des groupes de cohomologie.

\begin{proposition}\label{allunip} \begin{enumerate}
\item Pour  tout $\bt\in\uB^+$ et pour
tout $i$, le groupe  $H^i_c(\bX(\bt))$ est unipotent comme $\GF$-module,
et les valeurs propres de $F^\delta$  dans $H^i_c(\bX(\bt))$ sont de module
inf\'erieur ou \'egal \`a $q^{{\delta}i/2}$.
\item Si de plus  $\bt\in\uB^+$  est un  produit d'\'el\'ements $\ubw$
tels  que $\bX(\ubw)$ est rationnellement  lisse,
alors pour tout $i$ le groupe  $H^i_c(\bX(\bt))$
a  un caract\`ere  rationnel comme  $(\GF\times\genby{F^\delta})$-module.
\end{enumerate}
\end{proposition}
\begin{pf*}{Preuve}
Supposons pour commencer que $\bt$  est comme dans (ii).
Le corollaire \ref{rationnellement lisse} (ii)
montre  que  $H^*_c(\bX(\bt))$  est
unipotent et pour tout entier strictement positif $m$, on a:
\begin{align*}
|\bX(\bt)|^{gF^{\delta m}}&=\sum_i(-1)^i\Trace(gF^{\delta m}\mid
H^i_c(\bX(\bt))\\
&= 
\sum_i(-1)^i\sum_{\rho\in\CE(\GF,1)}(\lambda_\rho
q^{\delta i/2})^m\Trace(g\mid H^i_c(\bX(\bt))_\rho).
\end{align*}
Soit  $[\lambda]$   le  caract\`ere  de   $\genby{F^\delta}$  d\'efini
par   $F^\delta\mapsto\lambda$. On pose   $H^i_c(\bX(\bt))_\lambda=
\oplus_{\{\rho\mid\lambda_\rho=\lambda\}}H^i_c(\bX(\bt))_\rho$.
Alors,  pour    tout   $g$,   l'\'el\'ement   $$\sum_{i,\lambda}
(-1)^i[\lambda  q^{\delta   i/2}]\Trace(g\mid  H^i_c(\bX(\bt))_\lambda)$$
de  $\CR(\genby{F^\delta})$   ne  prend  que  des   valeurs  enti\`eres
sur les puissances strictement positives de $F^\delta$,
donc   est   stable    par   tout   $\sigma\in\Gal(\overline\bbQ/\bbQ)$.
Comme   $\sigma([\lambda  q^{\delta   i/2}])=  [\sigma(\lambda)q^{\delta
i/2}]$   et   que   les   caract\`eres   $[\lambda   q^{\delta   i/2}]$
sont  tous   distincts  quand  $i$   d\'ecrit $\bbZ$   et  $\lambda$
d\'ecrit les   racines    de    l'unit\'e,   on    en   d\'eduit    que
$\sigma(\Trace(g\mid    H^i_c(\bX(\bt))_\lambda)
=\Trace(g\mid H^i_c(\bX(\bt))_{\sigma(\lambda)})$. Par cons\'equent, le
caract\`ere du  $(\GF\times\genby{F^\delta})$-module $H^i_c(\bX(\bt))$ est
$\sigma$-stable, donc rationnel. On a prouv\'e (ii) ainsi que le cas
particulier de (i) o\`u $t$ est comme dans (ii).

On suppose maintenant que $\bt$ est  de la
forme  $\bt=\bt_1\cdots\bt_n$ avec  $t_1,\ldots,t_n\in  S\cup \uS$. On
va \'etablir (i), dans ce cas,  par
r\'ecurrence sur le nombre de $i$ tels que $t_i\in S$.
Soit    $i$    tel   que    $t_i\in    S$.    Alors,   la    vari\'et\'e
$\bX(\bt_1\cdots     \bt_{i-1}\ubt_i\bt_{i+1}\cdots\bt_n)$  est union de
l'ouvert  $\bX(\bt)$   et du  ferm\'e   compl\'ementaire
$\bX(\bt_1\cdots\bt_{i-1}\bt_{i+1}\cdots\bt_n)$.
Par     r\'ecurrence,   le (i) de la  proposition  est  \'etabli    pour
$\bX(\bt_1\cdots    \bt_{i-1}\ubt_i\bt_{i+1}\cdots\bt_n)$     et    pour
$\bX(\bt_1\cdots\bt_{i-1}\bt_{i+1}\cdots\bt_n)$. La  suite exacte longue
associ\'ee \`a  la d\'ecomposition  ouvert-ferm\'e montre  le r\'esultat
pour $\bX(\bt)$.

Consid\'erons  maintenant un  \'el\'ement  quelconque $\bt  =\bt_1\cdots
\bt_n$  de $\uB^+$  (avec  $t_i\in W\cup\uW$).  Nous  allons montrer (i)
par r\'ecurrence  sur  le  nombre de  $t_i$  dans $\uW$.  Le
r\'esultat est d\'ej\`a \'etabli pour $\bt\in B^+$.

Soit  $j$ tel  que  $t_j=\uw\in\uW$ et  consid\'erons la  stratification
de la proposition \ref{ouvertferme}, o\`u  $\bt=\bt_1\cdots\bt_{j-1}$ et
$\bt'=\bt_{j+1}\cdots\bt_n$.  Par  r\'ecurrence, tous  les  $\bX_i-\bX_{i-1}$
v\'erifient  (i).  La propri\'et\'e (i) pour $\bX(\bt)$ se d\'eduit
alors de la suite spectrale associ\'ee \`a la stratification. 
\end{pf*}

Nous  utiliserons  parfois  le   raffinement  suivant  de la proposition
\ref{allunip} (ii), qui \'etudie l'action de $F$.
Lorsque $F$ n'est pas un endomorphisme de Frobenius, la formule des traces
ne s'applique pas directement et nous contournons le probl\`eme
en utilisant  le  th\'eor\`eme  \ref{fuji} de  Fujiwara.

\begin{proposition}\label{rationalite de  F} Soit  $\bt\in\uB^+$
produit d'\'el\'ements
$\ubw\in\ubW^F$ tels  que $\bX(\ubw)$ est rationnellement  lisse; alors,
pour tout $i$, le caract\`ere du $(\GF\sdp F)$-module  $H^i_c(\bX(\bt))$
est rationnel.
\end{proposition}
\begin{pf*}{Preuve}
%


\'Ecrivons $\bt=\bt_1\cdots\bt_k$ avec $\bt_i\in\uW^F$ et 
$\bX(\bt_i)$ rationnellement lisse.
Soit $i\in\{0,\ldots,\delta-1\}$ et  $(\phi,F_1)=(gF^i,F^\delta)$.
Pour $m\equiv i\pmod \delta$, on a $gF^m=\phi(F^\delta)^{(m-i)/\delta}$,
donc, pour $m$  assez
grand, on a (\cf\ th\'eor\`eme \ref{fuji}):
$$|\bX(t_1,\ldots,t_k)|^{gF^m}=\sum_i(-1)^i\Trace(gF^m\mid
H^i_c(\bX(t_1,\ldots,t_k)).$$  

Puisque $\bX(t_1,\ldots,t_k)$   est   rationnellement   lisse, le
corollaire \ref{valeurs   propres}, joint aux conjectures de Weil,
montre que les   valeurs
propres   de   $F$   sur  $H^i_c(\bX(t_1,\ldots,t_k)_\rho$   sont   de
la   forme   $\mu_\rho  q^{i/2}$   o\`u   $\mu_\rho$   est  une   racine
$\delta$-i\`eme   de   $\lambda_\rho$.  Donc,   si  on   note
$H^i_c(\bX(t_1,\ldots,t_k))_\rho^{\mu_\rho}$    l'espace propre
g\'en\'eralis\'e de $F$ dans $H^i_c(\bX(t_1,\ldots,t_k))_\rho$ pour la
valeur propre $\mu_\rho q^{i/2}$, on a
\begin{multline*}
|\bX(t_1,\ldots,t_k)^{gF^m}|\\=
\sum_i(-1)^i\sum_\rho
\sum_{\{\mu_\rho\mid\mu_\rho^\delta=\lambda_\rho\}}(\mu_\rho
q^{i/2})^m\Trace(g\mid H^i_c(\bX(t_1,\ldots,t_k)^{\mu_\rho}_\rho),
\end{multline*}
et on conclut  comme dans la preuve du (ii) de la proposition
\ref{allunip}.

\end{pf*}

\begin{remarque}
Nous conjecturons que tout \'el\'ement de $C_{\uB^+}(F)$ est
produit d'\'el\'ements de  $C_{B^+}(F)$ et d'\'el\'ements de $\ubW^F$.
\end{remarque}

\sub{}
Soit $\Id$ le caract\`ere identit\'e de $\GF$. Soit $\Id_{h^2t}$
le caract\`ere de $\CH_{h^2t}(W)$ qui se
sp\'ecialise en le caract\`ere trivial de $W$ (il prend la valeur
$(h^2t)^{l(w)}$ sur $T_w$).

\begin{proposition}\label{Id} Soit $\bt\in\uB^+$.
La classe de la partie $\GF$-invariante de
$H^*_c(\bX(\bt))$ est donn\'ee par $\scal{H^*_c(\bX(\bt))}\Id=\Id_{h^2t}(T_{\bt})$.
\end{proposition}

\begin{pf*}{Preuve} Nous
proc\'edons par r\'ecurrence sur le nombre d'\'el\'ements de $\ubW$
dans une d\'ecomposition de $\bt$.

D\'emontrons d'abord le r\'esultat quand $\bt\in B^+$.
En consid\'erant une d\'ecomposi\-tion
$\bt=\bw_1\cdots\bw_k$ o\`u $\bw_i\in\bW$ et en utilisant
la proposition \ref{produit},
on se ram\`ene au cas o\`u $\bt\in\bW$.
Soit $\dot t$ un \'el\'ement de $N_\bG(\bT)$ qui rel\`eve $t$.
Soit $\bU$ le radical unipotent de $\bB$.

Soit $\bY(\dot t)=\{g\in\bG\mid g\inv F(g)\in\dot t\bU\}$. C'est une
vari\'et\'e munie d'une action de $\GF$ par multiplication \`a gauche et
et de $\bT^{tF}$ par multiplication \`a droite. On dispose aussi
d'une action de $\bU\cap \lexp {\dot t}\bU$ par multiplication \`a droite.
On a un isomorphisme de vari\'et\'es compatible \`a l'action de $\GF$
\cite[\S 1.11]{DL}
$$g\mapsto \lexp g \bB,\ (\bY(\dot t)/(\bU\cap \lexp {\dot t}\bU))/\bT^{tF}\iso \bX(t)$$
Le morphisme quotient $\pi:\bY(\dot t)\to \bY(\dot t)/
(\bU\cap \lexp {\dot t}\bU))$
est une fibration \`a fibres des espaces affines de dimension 
$l(w_0)-l(t)$. Par cons\'equent, $\pi_*\Qlbar\simeq
\Qlbar[2(l(t)-l(w_0))](l(t)-l(w_0))$, donc
on a un isomorphisme de $G^F$-modules
$H_c^*(\bX(\bt))\simeq
H_c^{*+2(l(w_0)-l(t))}(\bY(\dot t)/\bT^{tF})(l(w_0)-l(t))$.
La vari\'et\'e $\bY(\dot\bt)/\GF$ est isomorphe \`a $\bU$ et l'action correspondante
de $\bT^{tF}$ est l'action par conjugaison. Cette action s'\'etend
en l'action par conjugaison du groupe connexe $\bT$. Or, un groupe
connexe agit trivialement sur la cohomologie \cite[proposition 6.4]{DL}, donc
\begin{multline*}H_c^*(\bU)\simeq H_c^*(\bU)^{\bT^{tF}}\simeq
H_c^*(\bY(\dot\bt)/\bT^{tF})^\GF\\\simeq
H_c^{*-2(l(w_0)-l(t))}(\bX(\bt))^\GF(l(t)-l(w_0)).
\end{multline*}
Puisque $\bU$ est un espace affine de dimension $l(w_0)$, on en
d\'eduit l'\'egalit\'e $\scal{H^*_c(\bX(\bt))}\Id=(h^2t)^{l(\bt)}$.

Passons maintenant \`a l'\'etape g\'en\'erale de la r\'ecurrence.
Consid\'erons une d\'e\-composition
$\bt=\bt\ubw\bt'$ avec $\bt,\bt'\in \uB^+$ et $w\in W$.
On consid\`ere la stratification correspondante
(\cf\ proposition \ref{ouvertferme}).
L'hypoth\`ese de r\'ecurrence s'applique aux composantes des vari\'et\'es
$\bX_i-\bX_{i-1}$. En particulier, la partie $\GF$-invariante de leur
cohomologie est en degr\'e pair.
Il en r\'esulte que les morphismes de connection des suites exactes longues
ouvert-ferm\'e de cohomologie associ\'ees \`a la stratification
sont nulles lorsque l'on se restreint \`a la partie o\`u $\GF$ agit
trivialement.
Par cons\'equent,
$H^*_c(\bX(\bt))^{\GF}=\oplus_{i}H^*_c(\bX_i-\bX_{i-1})^{\GF}$,
d'o\`u la proposition.
\end{pf*}

Soit $\varepsilon$ le caract\`ere signe de $W$ et
$\tilde\varepsilon$ l'extension telle que $F$ agit trivialement.
Nous notons $\St$ le caract\`ere de Steinberg $R_{\tilde\varepsilon}$ de
$\GF$.

\begin{proposition}\label{Steinberg}
Soit $\bt\in\uB^+$.
On a  $\scal{H^*_c(\bX(\bt))}\St=\varepsilon_{h^2t}(T_\bt)$.
En d'autres termes:
\begin{enumerate}
\item
 La repr\'esentation de Steinberg n'intervient pas dans $H_c^*(\bX(\bt))$
si $\bt\not\in B^+$.
\item
Pour $\bt\in\BW$, le seul groupe de cohomologie o\`u la repr\'esentation de
Steinberg intervient est
$H_c^{l(\bt)}(\bX(\bt))$ ; dans celui-ci, elle a
multiplicit\'e 1 et est associ\'ee \`a la valeur propre 1 de $F^\delta$.
\end{enumerate}
\end{proposition}

\begin{pf*}{Preuve}
Lorsque   $\bt$   est   un   produit  d'\'el\'ements   de   $\ubS$,   le
r\'esultat   est   fourni   par  le corollaire  \ref{rationnellement
lisse} (ii),   car  $\varepsilon_{th^2}(T_\ubs)=0$   (si  $\ubs\in\ubS$)   et
$\varepsilon_{th^2}(T_1)=1$.

Soit maintenant $\bt=\bt_1\cdots\bt_n$ avec $t_1,\ldots,t_n\in S\cup\uS$. Soit
$i$ tel que $t_i\in S$. Par r\'ecurrence sur le cardinal
de $\{i|t_i\in S\}$, la repr\'esentation de Steinberg
n'intervient pas dans la cohomologie de la vari\'et\'e
$$\bX(\bt_1\cdots \bt_{i-1}\ubt_i\bt_{i+1}\cdots\bt_n).$$ La suite exacte
associ\'ee \`a la d\'ecomposition ouvert-ferm\'e
$\CO(\ubt_i)=\CO(\bt_i)\coprod\CO(1)$
montre que $H^k_c(\bX(\bt))_{\St}=
H^{k-1}_c(\bX(\bt_1\cdots\bt_{i-1}\bt_{i+1}\cdots\bt_n))_{\St}$
et la proposition en d\'ecoule par r\'ecurrence.

Soit $\bt=\bt_1\cdots\bt_n$ avec $t_1,\ldots,t_n\in W\cup\uW$.
Nous allons proc\'eder par r\'ecurrence sur $\sum_{i} (l(t_i)-1)$, o\`u
$i$ d\'ecrit les entiers tels que $t_i\notin W$.
La proposition est d\'ej\`a \'etablie lorsque cet entier est nul.
Sinon, il existe $i$ tel que $t_i=\uv\in\uW-\uS$.
 Soient $w\in W$ et $s\in S$ tels
que $v=ws$ et $w<v$. Par r\'ecurrence, la repr\'esentation de Steinberg
n'intervient pas dans la cohomologie des vari\'et\'es
$\bX(\bt_1\cdots\bt_{i-1}\ubw\ \ubs\bt_{i+1}\cdots\bt_n)$ et
$\bX(\bt_1\cdots\bt_{i-1}\ubw'\bt_{i+1}\cdots\bt_n)$ pour $w'<w$.
La proposition \ref{fibrationws}
permet alors de d\'eduire la proposition pour $\bX(\bt)$.
\end{pf*}

Pour $w\in W$, on pose $R_w=\sum_i(-1)^i[H_c^i(\bX(w))]\in\CR(G^F)$.

\begin{proposition}\label{H_s} Pour tout $\bs\in\bS^F$ et tout $n\in \bbN$,
on a
$$[H^*_c(\bX(\bs^n))]=\frac{t^n h^{2n}}{2}(R_1+R_s)+\frac{h^n}{2}(R_1-R_s).$$
\end{proposition}

\begin{pf*}{Preuve} Lorsque le rang de $\bG$ est $1$, la proposition r\'esulte des
propositions \ref{Id} et \ref{Steinberg}, car alors
$R_1+R_s=2[\Id]$ et $R_1-R_s=2[\St]$.
Le cas g\'en\'eral s'en d\'eduit par induction de Harish-Chandra
\`a partir du sous-groupe de Levi d\'efini par $s$ (\cf\ corollaire
\ref{RLG(H(Xs))}).
\end{pf*}

\begin{proposition}\label{s^mb} Soient  $\bs\in\bS$ et $\bb\in
B^+$  tels  que  $\bb  F(\bs)=\bs\bb$.   Alors,  il  existe  $H_s,  H_i\in
\bbZ[t^{1/2},h]\otimes\CR(\GF)$  tels que
$[H^*_c(\bX(\bs^m\bb))]=h^m H_s+(h^2 t)^m H_i$  pour tout  $m\in\bbN$.
\end{proposition}

\begin{pf*}{Preuve}  Nous appliquons  la proposition \ref{X_w  produit de  varietes}
avec  $I=\{s\}$   et  $\bw=\bs^m\bb$.   Soit  $\bz=\omega_I(\bb)$
et $n\in\bbN$  tel que  $\alpha_I(\bb)=\bs^n$. On a alors
$\omega_I(\bs^m\bb)=\bz$ et  $\alpha_I(\bs^m\bb)=\bs^{m+n}$.
Soit $\bz=\bz_1\cdots\bz_k$ la forme normale de $\bz$ et
$\dz_i$ un relev\'e de $z_i$.
On pose $F'=\dz_1\cdots\dz_k F$. Soit $\bL=\bL_I$   un  sous-groupe  de 
 Levi  de $\bG$  de  type  $A_1$. 
On note $\bY=\tilde\bX_{(I)}(\dz_1,\ldots,\dz_k)$, vari\'et\'e
 munie  d'une  action  \`a  droite de  $\bL^{F'}$  et  d'une  action
\`a   gauche   de   $\GF$.

On a
$\bX(\bs^m\bb)\simeq \bY\times_{\bL^{F'}}\bX_{\bL}(\bs^{m+n},F')$.
Le  groupe   $\bL^{F'}$   poss\`ede   deux
caract\`eres   unipotents,  $\Id_\bL$   et   $\St_\bL$  et la proposition
\ref{H_s}
montre que
$[H^*_c(\bX_{\bL}(\bs^{m+n},F'))]=h^{m+n}\St_\bL+(h^2t)^{m+n}\Id_\bL$.   On
d\'eduit de la  formule   de  K\"unneth que $$[H^*_c(\bX(\bs^m\bb))]
=h^{m+n}[H^*_c(\bY)_{\St_\bL}]  +(h^2t)^{m+n}[H^*_c(\bY)_{\Id_\bL}],$$
d'o\`u la  proposition  avec
$H_s=h^n[H^*_c(\bY)_{\St_\bL}]$  et $H_i=(h^2t)^n[H^*_c(\bY)_{\Id_\bL}]$.
\end{pf*}

\begin{proposition}\label{w0Ibar}
Soit $\tilde\chi\in\CR(W\sdp F)\otimes\bbQ$
tel que $R_{\tilde\chi}\in\CR^+(\GF)$ et
soit $I$ une partie $F$-stable de $S$. On suppose que
$\tilde{\chi}$ est combinaison lin\'eaire de caract\`eres
$\psi\in\Irr(W\sdp F)$ tels que $\scal{\Res^{W\sdp F}_{W_I}\psi}\Id=0$.

Alors, pour tous
$\bt,\bt'\in\uB^+$ et
tout caract\`ere unipotent $\rho$ tel que
$\scal\rho{R_{\tilde\chi}}\ne 0$, on $H^*_c(\bX(\bt\uw_0^I\bt'))_\rho=0$.
\end{proposition}

\begin{pf*}{Preuve} Nous proc\'edons comme dans la preuve de la proposition
\ref{Steinberg} pour nous ramener
au cas o\`u $\bt,\bt'$ sont dans le sous-mono\"\i de engendr\'e par $\uS$.

Puisque $T_{\uw_0^I}$ est un multiple de l'idempotent central
de $\CH_{h^2t}(W_I)$ associ\'e au caract\`ere trivial $\Id_{h^2t}$, il
agit par $0$ dans toute repr\'esentation $V$ de $\CH_{h^2t}(W)$ dont
la restriction \`a $W_I$ ne contient pas le caract\`ere trivial.
Par cons\'equent,
pour tout $\psi\in\Irr(W\sdp F)$ tel que
$\scal{\Res^{W\sdp F}_{W_I}\psi}\Id=0$,
on a $\psi_{h^2t}(T_{\bt\uw_0^I\bt'}F)=
\psi_{h^2t}(T_{\bt}T_{\uw_0^I}T_{\bt'}F)=0$.
Par hypoth\`ese, $\tilde\chi$ est combinaison lin\'eaire de caract\`eres $\psi$
ayant cette propri\'et\'e, donc $\tilde\chi_{h^2t}(T_{\bt\uw_0^I\bt'}F)=0$.

L'orthogonalit\'e des caract\`eres $R_{\tilde\phi}$ \cite[\S 3.19.2]{LuMa}
permet de d\'eduire du
corollaire \ref{rationnellement lisse} (ii) que
$\scal{H^*_c(\bX{(\bt\uw_0^I\bt')})}{R_{\tilde\chi}}=
\tilde\chi_{h^2t}(T_{\bt\uw_0^I\bt'}F)=0$. Puisque
$R_{\tilde\chi}\in\CR^+(\GF)$, on obtient finalement
$H^*_c(\bX(\bt\uw_0^I\bt'))_\rho=0$.
\end{pf*}

\begin{corollaire}\label{w0bar}
Pour  tous   $\bt,\bt'\in\uB^+$  et  pour  tout   caract\`ere  unipotent
$\rho\neq \Id$, on a $H^*_c(\bX(\bt\uw_0\bt'))_\rho=0$.
\end{corollaire}
\begin{pf*}{Preuve}
Lusztig \cite[\S 5.10 et theorem 6.17]{Lubook} a montr\'e qu'il existe une
famille d'\'el\'ements $a_{wF}\in \CR(W\sdp F)\otimes\bbQ$, pour $w\in W$,
v\'erifiant les propri\'et\'es suivantes (nous suivons les notations de
Lusztig au signe pr\`es):
\begin{itemize}
\item Pour $w\not=1$, on a $\scal{\Res_W a_{wF}}{\Id}=0$.
\item On a $R_{a_{wF}}\in\CR^+(\GF)$.
\item Le sous-groupe de $\CR(G^F)$ engendr\'e par les $R_{a_{wF}}$, pour $w\not=1$,
co{\"\i}ncide avec le sous-groupe engendr\'e par les $R_{\tilde\phi}$,
pour $\phi\not=\Id$.
\end{itemize}
Le r\'esultat se d\'eduit alors imm\'ediatement de la proposition \ref{w0Ibar},
appliqu\'ee \`a $\tilde\chi=a_{wF}$ pour $w\in W-\{1\}$.
\end{pf*}

\sub{} Nous g\'en\'eralisons maintenant des r\'esultats de \cite{LuMa}.
\begin{proposition}\label{repmin}
Soit $\rho$  un caract\`ere unipotent  de $\GF$ et
soit $\bb_0\in\BW$  un \'el\'ement de  longueur minimale tel  que
$H^*_c(\bX(\bb_0))_\rho\not=0$. Alors
\begin{enumerate}
\item $\bb_0\in \bW$ et, si $\rho^*$ est
le caract\`ere dual de $\rho$, l'\'el\'ement $\bb_0$ est aussi de longueur
minimale pour la propri\'et\'e $H^*_c(\bX(\bb_0))_{\rho^*}\not=0$. 
\item On a
$H^i_c(\bX(\bb_0))_\rho=0$ pour $i\not=l(\bb_0)$ et
$F^\delta$ agit sur $H^{l(\bb_0)}_c(\bX(\bb_0))_\rho$ avec une valeur
propre de module $q^{\delta l(\bb_0)/2}$.
\item   $\beta(\bb_0)$  est   de  longueur   minimale  dans   sa  classe   de
$F$-conjugaison et cette  classe de $F$-conjugaison ne d\'epend
que de $\rho$.
\item
Soit $\bb\in B^+$ tel que $H^{l(\bb)}_c(\bX(\bb))_\rho\not=0$.
Alors, les valeurs propres  de $F^\delta$ 
sur $H^{l(\bb)}_c(\bX(\bb))_\rho$  sont de  module inf\'erieur ou  \'egal \`a
$q^{\delta l(\bb)/2}$ et cette valeur  ne peut \^etre atteinte  que si
$\beta(\bb)$ est $F$-conjugu\'e \`a $\beta(\bb_0)$.
\end{enumerate}
\end{proposition}
\begin{pf*}{Preuve} 
L'assertion que $\bb_0$ est dans  $\bW$ est \cite[example 3.10.c]{LuMa}. On pose
$\bw=\bb_0$ et soit $w=s_1\cdots s_n$ une d\'ecomposition r\'eduite de $w$. Soit
$\bt'=\us_1\cdots\us_n$.
Le compl\'ementaire de l'ouvert
$\bX(\bw)$ de $\bX(\bt')$ est stratifi\'e
par des  vari\'et\'es  $\bX(\bb')$  avec $\bb'\in  B^+,
l(\bb')<l(\bw)$. Par minimalit\'e de $\bw$, les groupes de  cohomologie de ces
vari\'et\'es ne
font pas intervenir $\rho$. Les  suites exactes  longues de
cohomologie  montrent  alors  que  les composantes  de  type  $\rho$
de $H^i_c(\bX(\bw))$  et  de  $H^i_c(\bX(\bt'))$ sont  isomorphes.
D'apr\`es la proposition \ref{allunip} (ii), 
on a $H^*_c(\bX(\bt'))_{\rho*}\neq 0$. En \'echangeant les r\^oles de $\rho$
et $\rho^*$ on d\'eduit que
$\rho^*$ ne peut pas intervenir dans $H^*_c(\bX(\bw'))$ si $l(\bw')<l(\bw)$.
La stratification montre donc que
$H^*_c(\bX(\bw))_{\rho^*}\simeq H^*_c(\bX(\bt'))_{\rho^*}\neq 0$ et $\bw$ est aussi un \'el\'ement de longueur minimum tel que
$\rho^*$ intervient dans la cohomologie de la vari\'et\'e associ\'ee, d'o\`u (i).

On
d\'eduit  de la  proposition \ref{XtxXtp} (i)
que si  $H^i_c(\bX(\bw))_\rho\not=0$, alors, $i\geq l(\bw)$.

La vari\'et\'e  $\bX(\bt')$   est  une
vari\'et\'e projective lisse (\cf\  proposition
\ref{desingularisation}). La dualit\'e  de Poincar\'e
montre que  si $H^i_c(\bX(\bw))_\rho\not=0$, alors
$\rho$  et $\rho^*$  interviennent  aussi  dans
$H^{2l(\bw)-i}(\bX(\bt'))$, donc dans  $H^{2l(\bw)-i}_c(\bX(\bw))$. On a
donc aussi  $2l(\bw)-i\geq l(\bw)$, d'o\`u $i=l(\bw)$.  D'autre part, les
conjectures de Weil montrent que les valeurs propres de
$F^\delta$ sur $H^{l(\bw)}_c(\bX(\bt'))$ sont  de module
$q^{\delta l(\bw)/2}$. Il en  est donc de m\^eme pour 
de $H^{l(\bw)}_c(\bX(\bw))_\rho$ et l'assertion (ii) est \'etablie.

Prouvons (iii). Supposons que  $w$ n'est pas de longueur minimale  dans
sa $F$-classe de conjugaison. D'apr\`es
\cite[th\'eor\`eme 1.1]{GP} et \cite[th\'eor\`eme 2.6]{GKP}, il existe $s\in S$
tel que $w=sw'F(s)$  avec $l(w)=l(w')+2$, \ie, $\bw=\bs\bw'\bs$  dans $\BW$.
Alors, les propositions \ref{sws} et \ref{xy=yFx} montrent
que $H^*_c(\bX(\bw'))_\rho\not=0$ ou  $H^*_c(\bX(\bs\bw'))_\rho\not=0$,
ce qui  contredit la minimalit\'e de $\bw$.

Soit maintenant $\bw'\in\bW$ tel que $l(w')=l(w)$ et tel que
$H^{l(w)}_c(\bX(\bw'))_\rho\not=0$. D'apr\`es (ii), les
valeurs  propres de  $F^\delta$ sur $H^{l(w)}_c(\bX(\bw'))_\rho$
sont  de module  $q^{\delta
l(w)/2}$, donc $F^\delta$  a pour unique valeur propre $q^{\delta  l(w)}$
sur le sous-module
(non nul) $H^{l(w)}_c(\bX(\bw))_{\rho}\otimes H^{l(w)}_c(\bX(\bw'))_{\rho^*}$
de  $(H^{l(w)}_c(\bX(\bw))\otimes H^{l(w)}_c(\bX(\bw')))^\GF$.
La proposition \ref{XtxXtp} affirme alors que
$w$ et $w'$ sont $F$-conjugu\'es.

Consid\'erons enfin (iv).  On a vu dans (i) que 
$H^{l(w)}_c(\bX(\bw))_{\rho^*}\not=0$ et dans (ii) que $F^\delta$ agit sur cet espace
avec une valeur  propre de module  $q^{\delta l(w)/2}$.
Le  r\'esultat  r\'esulte  alors  de   la
proposition \ref{XtxXtp} appliqu\'ee \`a $(H^{l(\bb)}_c(\bX(\bb))\otimes
H^{l(w)}_c(\bX(\bw)))^\GF$.
\end{pf*}

Le corollaire suivant r\'esulterait de l'affinit\'e des vari\'et\'es
$\bX(\bb)$ pour $\bb\in B^+$, qui n'est connue que pour $q$ assez grand
(proposition \ref{quasiaffine} (iii)).

\begin{corollaire}\label{Haastert}
Pour tout $\bb\in\BW$ et tout $i<l(\bb)$ on a $H^i_c(\bX(\bb))=0$.
\end{corollaire}
\begin{pf*}{Preuve}
Soit   $\rho\in\CE(\GF,1)$ et
soit   $w\in W$  un \'el\'ement de longueur minimale tel que
$H^{l(w)}_c(\bX(w))_{\rho^*}\not=0$ (\cf\ proposition \ref{repmin}).
Soit   $i$   tel  que 
$H^i_c(\bX(\bb))_\rho\not=0$.
On  a donc  $(H^{l(\bw)}_c(\bX(\bw))\otimes H^i_c(\bX(\bb)))^\GF\neq0$.
On d\'eduit de la proposition \ref{XtxXtp} que
$l(\bw)+i\geq l(\bw)+l(\bb)$, donc que $i\geq l(\bb)$, d'o\`u le corollaire.
\end{pf*}

\sub{}
\label{periodicite}
Nous proposons maintenant plusieurs conjectures mettant en \'evidence
des ph\'enom\`enes de p\'eriodicit\'e. Elles sont v\'erifi\'ees dans les cas
des groupes de rang $2$ et des \'el\'ements pour lesquels nous calculons la
cohomologie (th\'eor\`emes \ref{A2}, \ref{2A2}, \ref{B2}, \ref{2B2},
\ref{G2} et \ref{2G2}).

On note $N=l(w_0)$. Pour  $\chi\in\Irr(W)^F$ nous notons $A_\chi$ (resp.
$a_\chi$) le degr\'e (resp. la  valuation) du degr\'e g\'en\'erique du 
caract\`ere $\chi_t\in\Irr(\CH_t(W))$ correspondant \`a $\chi$ 
(\cf\ \cite[4.1.1]{Lubook}). Par
ailleurs, comme  les  coefficients  $\scal{R_{\tilde\chi}}\rho$  sont
ind\'ependants  de  $q$  (\cf\   \cite[theorem  4.23]{Lubook}),  et  comme
$R_{\tilde\chi}(1)$ est  la valeur en  $q$ d'un polyn\^ome  \`a coefficients
ind\'ependants  de  $q$ (\og degr\'e fant\^ome\fg),
pour  tout $\rho\in  \CE(\GF,1)$  il  existe  un
polyn\^ome \`a coefficients ind\'ependants de $q$ dont $\rho(1)$ est la valeur
en $q$ (en effet, la partie unipotente de la repr\'esentation r\'eguli\`ere
est donn\'ee par $\sum_{\chi\in\Irr(W)^F}R_{\tilde\chi}(1)R_{\tilde\chi}$).
Nous d\'efinissons $A_\rho$  et $a_\rho$ comme le  degr\'e et la
valuation de ce polyn\^ome (voir aussi \cite[theorem 1.32]{BMM}).

\begin{conjecture}\label{conjA}
Soit $\bw\in\uB^+$. On a
$$[H_c^*(\bX(\bpi \bw))_\rho]=h^{4N-2A_\rho}t^{2N-a_\rho-A_\rho}
[H_c^*(\bX(\bw))_\rho].$$
\end{conjecture}

On note $P\mapsto \overline{P}$ l'involution de $\bbZ[t^{\pm 1/2},
h^{\pm 1}]$ donn\'ee par
$h\mapsto h\inv$ et $t^{1/2}\mapsto t^{-1/2}$.

\begin{conjecture}\label{conjB}
Soit $\bw\in\BW$ tel que $\bw\preccurlyeq\bpi^n$.
Alors on a
$$\scal{H_c^*(\bX(\bw\inv\bpi^n))}\rho=(h^{4N-2A_\rho}t^{2N-a_\rho-A_\rho})^n
\overline{\scal{H_c^*(\bX(\bw))}\rho}.$$
\end{conjecture}

Soit  $E$  l'involution  d'Ennola  de  $\CE(\GF,1)$  (\cf  \cite[theorem
3.3]{BMM};  {\it  loc.  cit.}  consid\`ere  une  permutation  avec  signes
$\sigma_\bbG$  dont  l'effet  sur  les dimensions  est  de  changer  $q$
en  $-q$;  nous  consid\'erons  ici  la  permutation  sous-jacente;  on  a
$E(\rho)=(-1)^{A_\rho}\sigma_\bbG(\rho)$).

\begin{conjecture}\label{conjC}
Supposons $w_0$ central dans $W$. 
Soit $\bw\in\uB^+$.
On a
$$\scal{H_c^*(\bX(\bw_0 \bw))}\rho=h^{2N-A_\rho}t^{N-(a_\rho+A_\rho)/2}
\scal{H_c^*(\bX(\bw))}{E(\rho)}.$$
\end{conjecture}

Un indice pour ces conjectures est la
\begin{proposition} Les conjectures \ref{conjA}, \ref{conjB} et
\ref{conjC} sont vraies pour $h=-1$.
\end{proposition}
\begin{pf*}{Preuve}
D'apr\`es le corollaire \ref{rationnellement lisse}(i), on a
$${\scal{H_c^*(\bX(\bw))}\rho}_{|h=-1}=\sum_{\chi\in\Irr(W)^F}
\scal{R_{\tilde\chi}}\rho\tilde\chi_t(T_\bw F).$$
En sommant sur
$\rho\in \CE(\GF,1)$, nous transformons la conjecture \ref{conjA}
sp\'e\-cia\-lis\'ee en $h=-1$ en l'assertion \'equivalente
\begin{multline*}\sum_{\chi\in\Irr(W)^F}\sum_{\rho\in\CE(\GF,1)}
\rho\scal{R_{\tilde\chi}}\rho\tilde\chi_t(T_\bpi T_\bw F)
\\=\sum_{\chi\in\Irr(W)^F}\sum_{\rho\in\CE(\GF,1)}
\rho t^{2N-a_\rho-A_\rho}\scal{R_{\tilde\chi}}\rho\tilde\chi_t(T_\bw F).
\end{multline*}
Comme    $A_\chi$   et    $a_\chi$    sont   les    m\^emes   pour    tous
les   $\chi$   tels   que   $\scal{R_{\tilde\chi}}\rho\ne   0$   (\cf\
\cite[theorem 4.23]{Lubook}), on a  $A_\rho=A_\chi$ et  $a_\rho=a_\chi$ lorsque
$\scal{R_{\tilde\chi}}\rho\ne  0$; nous  pouvons  donc r\'e\'ecrire  cette
\'egalit\'e
$$\sum_{\chi\in\Irr(W)^F} R_{\tilde\chi}\tilde\chi_t(T_\bpi T_\bw F)
=\sum_{\chi\in\Irr(W)^F}
t^{2N-a_\chi-A_\chi}R_{\tilde\chi}\tilde\chi_t(T_\bw F),$$
ce qui \'equivaut par ind\'ependance lin\'eaire des $R_{\tilde\chi}$
\`a l'assertion:
$\tilde\chi_t(T_\bpi T_\bw F)=t^{2N-a_\chi-A_\chi}\tilde\chi_t(T_\bw F)$.
Soit $\omega(\chi_t)$ le   caract\`ere   central  de   $\chi_t$.
Puisque $\bpi$   est   central, cela   revient  \`a   voir   que
$\omega(\chi_t)(T_\bpi)=t^{2N-a_\chi-A_\chi}$, \'egalit\'e bien connue (voir
par exemple \cite[5.12.2]{Lubook}).

Par    un    calcul    analogue, la conjecture \ref{conjB} se transforme en
$$\tilde\chi_t(T_\bw\inv     T_\bpi^n     F)=     t^{n(2N-a_\rho-A_\rho)}
\overline{\tilde\chi_t(T_\bw  F)}.$$
Vue la  formule  ci-dessus
pour   $\omega(\chi_t)(T_\bpi)$, ceci \'equivaut   \`a   $\tilde\chi_t(T_\bw\inv
F)=\overline{\tilde\chi_t(T_\bw F)}$;  ceci r\'esulte de ce  que $t\mapsto
t\inv, T_\bw\mapsto  T_\bw\inv, F\mapsto F$ est un  anti-automorphisme
semi-lin\'eaire de $\CH_t(W)\rtimes\langle F\rangle$ qui fixe les
caract\`eres irr\'eductibles (car il les fixe pour $t=1$, puisque tous
les caract\`eres irr\'eductibles de $W$ sont r\'eels et que nous avons choisi
l'extension $\tilde\chi$ rationnelle).

Enfin, par un calcul analogue, la conjecture \ref{conjC} se transforme en
$$\sigma_\bbG(R_{\tilde\chi})\tilde\chi_t(T_\bw T_{\bw_0}F)=
t^{N-(a_\rho+A_\rho)/2}R_{\tilde\chi}\tilde\chi(T_\bw F).$$
Puisque par {\it loc. cit.} on a
$\omega(\chi_t)(T_{\bw_0})=t^{N-(a_\rho+A_\rho)/2}\omega(\chi)(w_0)$, cela
revient \`a l'\'egalit\'e $\sigma_\bbG(R_{\tilde\chi})=\omega(\chi)(w_0) R_{\tilde\chi}$.
Cette \'egalit\'e est \`a  son tour \'equivalente \`a $\sigma_\bbG(R_w)=R_{ww_0}$.
Cette derni\`ere \'egalit\'e se ram\`ene  au cas des groupes quasi-simples; pour
les groupes exceptionnels c'est une v\'erification d'un nombre fini de cas
que  nous avons  effectu\'ee  par  ordinateur. Pour  les  cas des  groupes
classiques, nous pouvons extraire la preuve de celle de \cite[theorem 3.3]{BMM}.
En effet, avec les notations de {\it  loc. cit.}, $\bbT^-$ est un tore de
type $ww_0$ si  $\bbT$ est de type $w$;  l'\'enonc\'e \cite[theorem 3.3]{BMM}
suppose
le  tore  $d$-d\'eploy\'e  pour  un  certain  $d$;  mais  la  preuve  de  la
commutation de  l'induction de  Lusztig avec $\sigma_\bbG$  donn\'ee pages
52--53 de {\it loc. cit.} n'utilise  pas cette hypoth\`ese dans le cas des
tores.
\end{pf*}

\subsection{Ind\'ependance de $\ell$}
\label{sectionQ}

\sub{}
Soit $X$ une vari\'et\'e quasi-projective \'equidimensionelle sur
$\bar{\bbF}_p$.
On note $\Chow^r(X,n)$ les groupes de Chow sup\'erieurs de Bloch.
On pose $H^r_{mot}(X,\bZ(n))=\Chow^n(X,2n-r)$ (\og cohomologie motivique\fg). Nous
renvoyons \`a \cite[II, \S II.2]{Le} et
\cite{Ge} pour leurs principales propri\'et\'es.

Lorsque $X$ est lisse, on dispose de morphismes \og cycle \fg:
$H^r_{mot}(X,\bZ(n))\otimes_\bZ \bbQ_\ell\to H^r(X,\bbQ_\ell(n))$ et 
en sommant, on obtient
$$c^r(X):\bigoplus_{n}H^r_{mot}(X,\bZ(n))\otimes_\bZ \bbQ_\ell(-n)\to
H^r(X,\bbQ_\ell).$$
Si $\gamma$ est un endomomorphisme propre de $X$, alors il induit un
endomorphisme de $H^r_{mot}(X,\bZ(n))$ et $c^r(X)$ est \'equivariant pour
cette action.


\begin{proposition}
\label{Chow}
Soient $\bb,\bb'\in\uB^+$ tels que $\bX(\bb)$ et $\bX(\bb')$ sont lisses. On
consid\`ere l'action diagonale de $\bG^F$ sur $\bX(\bb)\times\bX(\bb')$.
Alors, pour tout $r$, le morphisme cycle $c^r$
induit un isomorphisme
$$\bigoplus_n H^r_{mot}(\bX(\bb)\times\bX(\bb'),\bZ(-n))^{\bG^F}
\otimes_\bZ\bbQ_{\ell}(n)\iso H^r(\bX(\bb)\times\bX(\bb'),\bbQ_\ell)^{\bG^F}.$$
\end{proposition}

\begin{pf*}{Preuve}
Commen\c{c}ons par introduire deux d\'efinitions.
Un morphisme $f:Y\to X$ est un un {\em quasi fibr\'e vectoriel
de rang $r$} si c'est
une fibration localement triviale (pour la topologie de Zariski),
de fibres isomorphes \`a l'espace affine $\bbA^r$.
Le morphisme $f^*:H^r_{mot}(X,\bZ(n))\to H^r_{mot}(Y,\bZ(n))$ est alors
un isomorphisme (\cf\ par exemple \cite[\S 1.2.7]{Ge}).
On introduit une relation d'\'equivalence sur les vari\'et\'es comme la
cl\^oture de la relation d'\^etre un quasi-fibr\'e vectoriel et on dit
qu'une vari\'et\'e est un {\em quasi espace affine} si elle est dans
la classe d'un point. Pour un quasi espace affine $X$, les morphismes $c^r$
sont donc des isomorphismes.

Soient $t_1,\ldots,t_n,t'_1,\ldots,t'_{n'}\in\uS$.
Lusztig \cite[\S 2.6]{LuHom} introduit une stratification par des
sous-vari\'et\'es ferm\'ees $(\bG^F\times\langle F\rangle)$-\'equivariantes
$$0= Z_{L_0}\subset Z_{L_1}\subset\cdots\subset Z_{L_m}=
\bX(t_1,\ldots,t_n)\times \bX(t'_1,\ldots,t_{n'}).$$
On fixe $i$ et on pose $Z=Z_{L_i}-Z_{L_{i-1}}$.
Lusztig \cite[\S 3.3 et 3.5]{LuHom} introduit des vari\'et\'es lisses
$Z_0$ et $Z_1$ munies de l'action libre
d'un groupe fini $\CT$ commutant aux actions de $\bG^F$ et de $F$.
Il construit un isomorphisme
$(\bG^F\times\langle F\rangle)$-\'equivariant $Z_0/\CT\iso Z$
\cite[Lemma 3.4]{LuHom} et un quasi fibr\'e vectoriel
$(\bG^F\times\langle F\rangle\times \CT)$-\'equivariant
$f:Z_1\to Z_0$. Il montre \cite[\S 3.26]{LuHom} que 
$\bG^F\setminus Z_1$ est un quasi espace affine.

On en d\'eduit que $c(Z_1)^{\bG^F}$ est un isomorphisme, donc que
$c(Z_0)^{\bG^F}$ est un isomorphisme et finalement que
$c(Z)^{\bG^F}$ est un isomorphisme.

Les suites spectrales qui calculent la comohologie $\ell$-adique et
la cohomologie motivique
de $\bX(t_1,\ldots,t_n)\times \bX(t'_1,\ldots,t_{n'})$ \`a partir de
celle des $Z_{L_i}-Z_{L_{i-1}}$ sont compatibles avec les morphismes
\og cycle\fg. On d\'eduit que
$c(\bX(t_1,\ldots,t_n)\times \bX(t'_1,\ldots,t_{n'}))^{\bG^F}$
est un isomorphisme.

\smallskip
Soient maintenant $s_1,\ldots,s_n,s'_1,\ldots,s'_{n'}\in S$. Soient
$t_i=\us_i$ et $t'_i=\us'_i$. D'apr\`es ce qui pr\'ec\`ede,
les morphismes $c^{\bG^F}$ sont des isomorphismes pour les
vari\'et\'es $\bX(t_{i_1},\ldots,t_{i_r})\times
\bX(t'_{i'_1},\ldots,t'_{i'_{r'}})$, o\`u $1\le i_1<\cdots i_r\le n$ et
$1\le i'_1<\cdots i'_{r'}\le n'$. On d\'eduit de la suite spectrale
de Mayer-Vietoris qui calcule la cohomologie de
$\bX(s_1,\ldots,s_n)\times\bX(s'_1,\ldots,s'_{n'})$ \`a partir de celle des
$\bX(t_{i_1},\ldots,t_{i_r})\times \bX(t'_{i'_1},\ldots,t'_{i'_{r'}})$ que
$c(\bX(s_1,\ldots,s_n)\times\bX(s'_1,\ldots,s'_{n'}))^{\bG^F}$ est
un isomorphisme.

\smallskip
Dans le cas g\'en\'eral, on utilise la suite spectrale qui provient
de la stratification (\ref{lt ferme}) de $\bX(\bb)$ dont les pi\`eces sont des
$\bX(s_1,\ldots,s_n)$ avec $s_i\in S$.
\end{pf*}

\begin{remarque}
Lorsque $\bb$ et $\bb'$ sont produit d'\'elements de $\uS$, alors, 
$\Chow^r(\bX(\bb)\times\bX(\bb'),n)^{\bG^F}=0$ pour $n\not=0$: parmi les
groupes de Chow sup\'erieurs de $\bX(\bb)\times\bX(\bb')$, seuls les groupes de
Chow ordinaires ont des invariants sous $\bG^F$.
\end{remarque}

Par dualit\'e de Poincar\'e, on d\'eduit de la proposition \ref{Chow}
le corollaire suivant.

\begin{corollaire}
\label{corChow}
Soient $\bb,\bb'\in\uB^+$ tels que $\bX(\bb)$ et $\bX(\bb')$ sont lisses.
Alors, pour tout $i$, il existe un $\bbQ\langle F\rangle$-module
$L^i$ tel que pour tout $\ell{\not|}q$, on a un isomorphisme de
$\bbQ_\ell\langle F\rangle$-modules
$$\bbQ_\ell\otimes_\bbQ L^i\iso 
H^i_c(\bX(\bb)\times\bX(\bb'),\bbQ_\ell)^{\bG^F}.$$
\end{corollaire}
\sub{}
\label{indep}
\newcommand{\CF}{{\mathcal F}}
Pour $\bb\in\uB^+$, soit $H^i_c(\bX(\bb),\bbQ_\ell)_j$ la partie de
$H^i_c(\bX(\bb),\bbQ_\ell)$ o\`u les valeurs propres de $F^\delta$ sont de
module $q^{\delta j/2}$. 
D'apr\`es le corollaire \ref{corChow}, on  a, pour $\bt,\bb\in\uB^+$, en 
 consid\'erant la   partie   de
$H^k_c(\GF\backslash  \bX(\bt)\times\bX(\bb))$ o\`u  $F^\delta$ agit par
$q^{\delta       k'/2}$,      que
$$S_{k,k',\bt,\bb}:=\sum_{i,j}\scal{H^{k-i}_c(\bX(\bt),\bbQ_\ell)_{k'-j}}
{H^i_c(\bX(\bb),\bbQ_\ell)^*_{j}}$$  est ind\'ependant  de $\ell$ (o\`u pour une
repr\'esentation $V$ de $\GF$ on note $V^*$
la repr\'esentation contragr\'ediente).

\begin{proposition}\label{Rchi indepl} Pour tout $\bb\in\uB^+$,
tout $\chi\in\Irr(W)^F$ et tous $i, j$, le produit scalaire
$\scal{R_{\tilde\chi}}{H^i_c(\bX(\bb),\bbQ_\ell)_j}$ est ind\'ependant de 
$\ell$ et est nul si $a_\chi\not\equiv j\pmod 2$.
\end{proposition}
\begin{pf*}{Preuve}
Si $\bt$ est dans le sous-mono\"\i de de $\uB^+$ engendr\'e par $\uS$, on a
par le corollaire \ref{rationnellement lisse}(ii):
$$H^*_c(\bX(\bt),\bbQ_\ell)=\sum_{\chi\in\Irr(W)^F}R_{\tilde\chi}
\tilde\chi_{h^2t}(T_\bt F).$$
En   tenant   compte   de   cette   formule   on  a  $$S_{k,k',\bt,\bb}=
\sum_i\scal{H^{k-i}_c(\bX(\bt),\bbQ_\ell)_{k-i}}
{H^i_c(\bX(\bb),\bbQ_\ell)^*_{i+k'-k}};$$  et comme  les valeurs propres
de $F^\delta$ sur $H^k_c(\GF\backslash \bX(\bt)\times\bX(\bb))$ sont des
puissances  enti\`eres de $q^\delta$
(voir corollaire \ref{valeurs propres} (ii)),
cette somme  est nulle si $k'$ est
impair.  
Soit $\tilde\chi_{h^2m}(T_\bt F)_{(j)}$
le  coefficient de $(h^2t)^{j/2}$ dans $\tilde\chi_{h^2m}(T_\bt F)$ et soit
$$f(k,k',m):=\sum_{i,\chi\in\Irr(W)^F}\tilde\chi_{h^2t}(m
F)_{(k-i)}   \scal{R_{\tilde\chi}}{H^i_c(\bX(\bb),\bbQ_\ell)_{i+k'-k}}$$
pour $m\in \CH_{h^2t}(W)$.
On trouve donc que
$f(k,k',T_\bt)$
est  ind\'ependant de $\ell$ pour tout $\bt$ dans le sous-mono\"\i de de
$\uB^+$ engendr\'e par $\uS$ et pour tous $k,k'$ et est nul si
$k'$  est  impair. Nous  avons  retir\'e  l'\'etoile sur
$H^i_c(\bX(\bb),\bbQ_\ell)_{i+k'-k}$  en utilisant  que $R_{\tilde\chi}$
est rationnel, car $\tilde{\chi}$ est rationnel.

Pour $w\in W$, on fixe une d\'ecomposition r\'eduite
$w=s_1\cdots s_n$ et on pose $B_w=T_{\us_1}\cdots T_{\us_n}$. Alors,
$\{B_w\}_{w\in W}$ est une base de $\CH_{h^2t}(W)$.
Soit $\{C_w\}_{w\in W}$ la base de Kazhdan-Lusztig de $\CH_{h^2t}(W)$.
Il existe des entiers $p_{v,w,i}$ tels que
$C_w=\sum_{i,v}p_{v,w,i} (h^2t)^i B_v$.
On a $f(k,k',C_w)=\sum_{i,v}p_{v,w,i}f(k-2i,k'-2i,B_v)$.
On en d\'eduit que  $f(k,k',C_w)$ est ind\'ependant de
$\ell$ et est nul si $k'$ est impair.

Nous fixons maintenant $\bb$ et nous notons $g(\chi,i,j)$ l'assertion de
l'\'enonc\'e,                         \`a                        savoir:
$\scal{R_{\tilde\chi}}{H^i_c(\bX(\bb),\bbQ_\ell)_j}$  est  ind\'ependant
de $\ell$ et nul si $a_\chi\not\equiv j\pmod 2$. Nous allons d\'emontrer
$g(\chi,i,j)$ par r\'ecurrence d\'ecroissante sur $a_\chi$ et croissante
sur  $i$. Nous fixons  donc une famille  $\CF$, nous posons $a=a_\CF$ et
nous  fixons un \'el\'ement $w$  de la cellule bilat\`ere d\'etermin\'ee
par  $\CF$. Nous supposons d\'emontr\'ee $g(\chi,i',j)$ pour tout $j$ et
tout  $i'$ si $a_\chi>a$, ou pour tout $j$ et $i'<i$ si $a_\chi=a$. Nous
voulons    en   d\'eduire    $g(\chi,i,j)$   pour    $\chi\in\CF$.   Par
\cite[5.2.1]{Lubook}   on   a   $\tilde\chi_{h^2t}(C_w  F)_{(i')}=0$  si
$i'<-a_\chi$;  et en lisant la preuve  de \cite[5.2]{Lubook} on voit que
$\tilde\chi_{h^2t}(C_wF)=0$  si  $a_\chi<a$,  ou  si $a_\chi=a$ et $\chi
\notin\CF$.  Comme dans {\sl loc. cit.} nous notons $c_{wF,\tilde\chi}:=
(-1)^{l(w)}\tilde\chi_{h^2t}(C_w F)_{(-a_\chi)}$.

On obtient
\begin{multline*}
f(i-a,j-a,C_w)=
(-1)^{l(w)}\sum_{\chi\in\CF}c_{wF,\tilde\chi}\scal{R_{\tilde\chi}}
{H^i_c(\bX(\bb),\bbQ_\ell)_j}
+\\
\sum_{i',\{\chi\in\Irr(W)^F\mid a_\chi>a\}}
\tilde\chi_{h^2t}(C_w F)_{(-a+i-i')}
\scal{R_{\tilde\chi}}{H^{i'}_c(\bX(\bb),\bbQ_\ell)_{i'+j-i}}+\\
\sum_{i'<i,\chi\in\CF}
\tilde\chi_{h^2t}(C_w F)_{(-a+i-i')}
\scal{R_{\tilde\chi}}{H^{i'}_c(\bX(\bb),\bbQ_\ell)_{i'+j-i}}.
\end{multline*}

Les  autres termes \'etant nuls par les remarques qui pr\'ec\`edent. Par
hypoth\`ese   de   r\'ecurrence   les   deux   derni\`eres  sommes  sont
ind\'ependantes   de  $\ell$  (en   utilisant  $g(\chi,i',i'+j-i)$  pour
$a_\chi>a$  et    pour  $i'<i$  , respectivement). On obtient donc que
$\sum_{\chi\in\CF}c_{wF,\tilde\chi}\scal{R_{\tilde\chi}}
{H^i_c(\bX(\bb),\bbQ_\ell)_j}$, qui   s'\'ecrit
$\scal{R_{a_{wF}}}{H^i_c(\bX(\bb),\bbQ_\ell)_j}$  avec les  notations de
\cite[6.17.2]{Lubook}, est ind\'ependant de $\ell$ pour tout $w$ et est
nul si $j\not\equiv a\pmod 2$.

Or d'apr\`es \cite[(b) above 6.17.1]{Lubook} l'espace engendr\'e par les
$R_{a_{wF}}$    co\"\i   ncide    avec   celui    engendr\'e   par   les
$R_{\tilde\chi}$, d'o\`u le r\'esultat.
\end{pf*}

Rappelons une version du th\'eor\`eme de densit\'e de Chebotarev.

\begin{theoreme}
\label{Chebotarev}
Soit $L$ une extension galoisienne de $\bbQ$. 
Il existe une infinit\'e de $\ell$ tels que
$L\cap \bbQ_\ell=\bbQ$ (pour tout plongement de $\bbQ_\ell$ dans $\bbC$).
\end{theoreme}

\begin{proposition}\label{indep conj}
Soit  $\rho$ un caract\`ere unipotent et soit  $b\in  \uB ^+$.
Supposons que pour tous $i$ et $j$,   le   produit   scalaire
$\scal\theta{H^i_c(\bX(\bb),\bbQ_\ell)_j}$  est ind\'ependant de $\ell$,
pour  $\theta=\rho+\rho^*$ et pour $\theta$ un caract\`ere unipotent
irr\'eductible  diff\'erent de $\rho$  et $\rho^*$. Alors  pour tous $i$,
$j$,         on        a        $\scal\rho{H^i_c(\bX(\bb),\bbQ_\ell)_j}=
\scal{\rho^*}{H^i_c(\bX(\bb),\bbQ_\ell)_j}$  et  ce  produit scalaire est
ind\'ependant de $\ell$.
\end{proposition}
\begin{pf*}{Preuve}
L'hypoth\`ese   implique  que   la  partie   $\rho,\rho^*$-isotypique  de
$H^i_c(\bX(\bb),\bbQ_\ell)_j$  que nous noterons  simplement $H^i_j$ est
de  la  forme  $a_{i,j}\rho+a'_{i,j}\rho^*$  o\`u  $a_{i,j}+a'_{i,j}$ est
ind\'ependant  de $\ell$.  Nous prouvons  la conclusion par r\'ecurrence
croissante  sur $i$, et pour $i$ donn\'e par r\'ecurrence croissante sur
$j$.  Par r\'ecurrence, on sait donc que pour $k<i$ et tout $l$, ou pour
$k=i$   et  $l<j$, on a $a'_{k,l}=a_{k,l}$.
En ne retenant que les termes correspondant \`a la
partie $\rho,\rho^*$-isotypique dans $S_{2i,2j,\bb,\bb}$ on obtient que
\begin{multline*}
\sum_{k,l}\scal{H^k_{l}}{(H^{2i-k}_{2j-l})^*}=
2a_{i,j} a'_{i,j}
+\sum_{l<j}2a_{i,l}(a_{i,2j-l}+a'_{i,2j-l})\\
+\sum_l\sum_{k<i}2a_{k,l}(a_{2i-k,2j-l}+a'_{2i-k,2j-l})
\end{multline*}
est ind\'ependant
de  $\ell$. On obtient  donc que $a_{i,j}a'_{i,j}$  est ind\'ependant de
$\ell$. D'apr\`es le th\'eor\`eme \ref{Chebotarev},
il existe une infinit\'e de $\ell$ tels que ni $\rho$
ni  $\rho^*$ ne soient \`a valeur  dans $\bbQ_\ell$. Pour de tels $\ell$,
on      a     donc      n\'ecessairement     $a_{i,j}=a'_{i,j}$     donc
$a_{i,j}a'_{i,j}=(a_{i,j}+a'_{i,j})^2/4$.  Pour  tout  $\ell$ on a
donc    $a_{i,j}a'_{i,j}=(a_{i,j}+a'_{i,j})^2/4$    ce    qui   implique
$a_{i,j}=a'_{i,j}=(a_{i,j}+a'_{i,j})/2$. \end{pf*}

\section{Cohomologies dans les groupes r\'eductifs de rang 2}
\label{sectionrang2}

\subsection{G\'en\'eralit\'es}
Nous  allons  d\'eterminer  la   cohomologie  d'un  certain  nombre  de
vari\'et\'es de Deligne-Lusztig  dans des groupes $\bG$ de rang $2$
(d\'eploy\'es  ou  non).  Dans  tous  les  cas  notre  r\'esultat
principal   sera   une   \og quasi-p\'eriodicit\'e\fg\   des   cohomologies
calcul\'ees par rapport \`a la  multiplication par un \'el\'ement central
de $B^+$.

Les  parties $\Id$-isotypiques  et $\St$-isotypiques  de la  cohomologie
d'une  vari\'et\'e  $\bX(\by)$  sont  connues  par  les propositions
\ref{Id}  et \ref{Steinberg}.   Nous   nous   int\'eressons   donc
seulement   \`a
la   partie   isotypique   de   ces  cohomologies   pour   les   autres
caract\`eres   unipotents.    Si   $\rho_1,\ldots,\rho_n$    sont   ces
caract\`eres,   nous   \'ecrirons   alors   $H(\by)$   pour   la   partie
$(\rho_1,\ldots,\rho_n)$-isotypique  de  $\sum_i  h^i  H_c^i(\bX(\by))$,
repr\'esent\'ee  suivant  notre   convention  de \S \ref{convention}  par  un
\'el\'ement de $\bbZ[t^{1/2},h][\rho_1,\ldots,\rho_n]$.

Nous commen\c cons  par donner un certain nombre  de cons\'equences des
propositions  du \S \ref{sectioncohomologie}.
Nous notons  $S=\{s,t\}$ les  g\'en\'erateurs de $W$.

Pour  raccourcir les d\'emonstrations, si  $x,y,z\in\uB^+$ nous  \'ecrirons
dans le texte  qui suit qu'un r\'esultat est obtenu  \og par $\seof xyz$\fg\ 
pour  signifier que  la  vari\'et\'e $\bX(x)$  est une  sous-vari\'et\'e
ouverte  de  $\bX(y)$, que  $\bX(z)$  est  la sous-vari\'et\'e  ferm\'ee
compl\'ementaire  et  que  le  r\'esultat est  obtenu  en  consid\'erant
la suite  exacte longue $\cdots   H^i_c(\bX(x))\to  H^i_c(\bX(y))\to
H^i_c(\bX(z))\to H^{i+1}_c(\bX(x))\to\cdots$ qui en r\'esulte.

\begin{lemme}\label{rappel}
Soit $y\in\uB^+$. Alors,
\begin{enumerate}
\item pour tout $x\in\uB^+$ on a $H(xy)=H(yF(x))$
\item $H(\us\,\ut y)=H(\us\ut y)$
\item $H(s \us y)=H(\us s y)=h^2 t H(\us y)$
\item $H(\us\,\us y)= (h^2 t+1) H(\us y)$
\item $H(\uw_0y)=0$.
\end{enumerate}
\end{lemme}
\begin{pf*}{Preuve} (i) r\'esulte de la proposition \ref{xy=yFx}.
(ii)  est  une   relation  du  mono\"\i  de  $\uB^+$.   (iii)
et (iv) proviennent de la proposition \ref{wss} appliqu\'ee  \`a  $w=s$. 
(iv)  provient  du corollaire \ref{w0bar}.
\end{pf*}

\subsection{Type $A_2$}
\sub{}
Nous consid\'erons maintenant le cas  d'un groupe (d\'eploy\'e ou non) de
type $A_2$. Dans  ce cas il existe un unique  caract\`ere unipotent
$\rho$ diff\'erent de $\St$  et $\Id$. Nous
omettrons la mention  de $\rho$ dans la notation pour  $H(\by)$ qui sera
repr\'esent\'e par un \'el\'ement de $\bbZ[t^{1/2},h]$. 

On note $R_x$ la repr\'esentation de 
$\CH_x(W)\sdp   F$  donn\'ee   par  $T_s\mapsto\begin{pmatrix}-1&0\cr
\sqrt{x}&x\cr\end{pmatrix}$,
$T_t\mapsto\begin{pmatrix}x&\sqrt{x}\cr  0&-1\cr\end{pmatrix}$
et   $F\mapsto\begin{pmatrix}1&0\cr0&1\cr\end{pmatrix}$   si    $\bG$   est 
d\'eploy\'e, $F\mapsto\begin{pmatrix}0&1\cr1&0\cr\end{pmatrix}$ sinon.

\begin{lemme}\label{rappelA2}
Soit $y\in\uB^+$. Alors,
\begin{enumerate}
\item $H(\us\,\ut\,\us y)= H(\underline{st}\,\us y)=h^2t H(\us y)$.
\item $H(\ut\,\bs\,\ut y)= h H(\ut y)$.
\item $H(\us\bt\bs y)=h H(\us\ut y)$.
\item
Soit  $\sigma$  l'automorphisme  de  $\uB^+$ donn\'e par l'\'echange
de $\bs$ (resp.  $\us$) et $\bt$ (resp.  $\ut$)
(c'est  l'action  de   $F$  si  $(\bG,F)$  est   non  d\'eploy\'e,  \ie\
de  type   $\lexp  2A_2$).  Alors  $H(y)=H(\sigma(y))$.
\item
Si  $y$ est  dans le  sous-mono{\"\i}de engendr\'e  par $\underline  W$,
alors  on a $H(y)=\Trace(T_y  F\mid  R_{h^2t})$.
\end{enumerate}
\end{lemme}
\begin{pf*}{Preuve}
(i) vient des relations de $\uB^+$,
de la proposition \ref{fibrationws} appliqu\'ee \`a $w=st$ et $y=s$, et de \ref{rappel}(v).
(ii) est obtenu par \ref{rappel}(v) et par
$\seof{\ut\bs\ut y}{\underline{sts} y}{\ut y}$.
(iii) est obtenu par \ref{rappel}(v) et par
$\seof{\us\bt\bs y}{\underline{sts} y}{\us\ut y}$.
(iv) est une cons\'equence de la proposition \ref{autG} en prenant comme 
isog\'enie l'automorphisme d'opposition de $\bG$.
Enfin  (v) est  une cons\'equence  imm\'ediate du corollaire \ref{rationnellement
lisse} (ii)  et  du  fait que  dans  un  groupe  de  type $A_2$  toutes  les
vari\'et\'es  $\bX(\uw)$ pour  $w\in  W$ sont  rationnellement
lisses (\cf\ proposition \ref{X rationnellement lisse}).
\end{pf*}

\sub{}
Nous  d\'ecrivons  maintenant  pour  un certain  nombre  de  vari\'et\'es  une
p\'eriodicit\'e   par  rapport   \`a  la   multiplication  par   $\bpi$;  nous
conjecturons qu'une  telle p\'eriodicit\'e a  lieu pour toutes  les vari\'et\'es
associ\'ees \`a des \'el\'ements de $\uB^+$.
\begin{theoreme}\label{A2}
Soit $(\bG,F)$ un groupe d\'eploy\'e de type $A_2$, soit $\bpi=(\bs\bt)^3$ et
soit $y\in\uB^+$ apparaissant dans la
table ci-dessous. Alors, pour tout $n$,
on a $H(y\bpi^n)=(h^8 t^3)^n f(y)$ o\`u $f(y)$ est la valeur
de $H(y)$ donn\'ee par la table.

\halign{$#$\hfill&\quad$#$\hfill\cr
\strut y&H(y)\cr
\noalign{\hrule}\cr
1&2\cr
\bs,\bt&h^2t + h\cr
\us\,,\ut\,& h^2t + 1\cr
\bs\bs,\bt\bt&h^2+h^4t^2\cr
\us\ut,\ut\us&h^2t\cr
\us\bt,\ut\bs&h\cr
\bs\bt,\bt\bs&h^3t\cr
\bs\bt\bs, \bt\bs\bt, \bs\bt\bt, \bs\bs\bt&0\cr
\bs\ut\bs, \bt\us\bt, \us\bt\bt&h^4t^2 + h^2\cr
\bs\bt\bs\bt, \bs\bs\bs\bt&h^5t^2\cr
\bt\bs\bs\bt,\bs\bs\bt\bt&h^6t^3+h^4t\cr
\bt\bs\bt\bs\bt=\bs\bt\bs\bs\bt&h^7t^3+h^6t^2\cr
}
\end{theoreme}
\begin{pf*}{Preuve}
Nous proc\'edons par r\'ecurrence sur $n$. Nous traitons en m\^eme temps
le point  de d\'epart de la  r\'ecurrence ($n=0$) et le  cas g\'en\'eral
(o\`u  l'on suppose  le th\'eor\`eme  prouv\'e pour  $n-1$) pour  ne pas
avoir \`a dupliquer certains raisonnements identiques dans les deux cas.
Nous supposons donc $H(y\bpi^{n-1})$ connu pour tout $y$ dans la table,
et nous  allons d\'eterminer  $H(y\bpi^n)$, en raisonnant  par longueur
croissante de $y$ dans la table du th\'eor\`eme \ref{A2}.
Notons qu'un r\'esultat sur $H(x)$ est \'equivalent au r\'esultat
correspondant pour $H(\sigma(x))$ d'apr\`es le lemme \ref{rappelA2} (iv).

\casde{\bpi^n},  $\bs\bpi^n$ et $\bs^2\bpi^n$.
Si   $n=0$   la   valeur   se  d\'eduit   de la proposition \ref{H_s}.
Sinon,  nous   appliquons la proposition \ref{s^mb}
pour  $m=0,1,2,3,4$  et  $\bb=\bt\bs\bs\bt\bpi^{n-1}$.  Les  valeurs  de
$H(\bs^m\bb)$  sont  dans  la  table  (connues  par  r\'ecurrence)  pour
$m=0,1$ et permettent de d\'eterminer  avec les notations de la proposition \ref{s^mb}:
$H_i=h^4t(h^8t^3)^{n-1}$  et   $H_s=h^6t^3(h^8t^3)^{n-1}$,  d'o\`u,  pour
tout $m$,  on a  $H(\bs^m\bb)=(h^8t^3)^{n-1}(h^{m+6}t^3+h^4t(h^2t)^m)$ ce
qui  donne  les valeurs  recherch\'ees  pour  $m=2,3,4$ en  tenant  compte
de  $\bs^2\bb=\bpi^n$. 

\casde{\us\bpi^n}.
Si $n=0$ la  valeur se d\'eduit du lemme \ref{rappelA2}(v).
Sinon  on  a
\begin{multline*}
$$H(\us\bpi^n)=H(\us\bs\bs\bt\bs\bs\bt\bpi^{n-1})=  h^4t^2
H(\us\bt\bs\bs\bt\bpi^{n-1})=h^5t^2H(\us\ut\bs\bt\bpi^{n-1})=\\
=h^6t^2H(\us\,\ut\,\us\bpi^{n-1})=h^8t^3H(\us\bpi^{n-1})
\end{multline*}
d'apr\`es les lemmes  
\ref{rappel}(iii) et  \ref{rappelA2}(iii) puis encore (iii) et (i).

\casde{\us\ut\bpi^n}. Si $n=0$ on utilise le lemme \ref{rappelA2}(v) et sinon
on a
\begin{multline*}
H(\us\,\ut\bpi^n)= H(\us\,\ut\bt\bt\bs\bt\bt\bs\bpi^{n-1})=
h^4t^2H(\us\,\ut\bs\bt\bt\bs\bpi^{n-1})=
h^5t^2H(\us\,\ut\,\us\bt\bs\bpi^{n-1})=\\
=h^7t^3H(\us\bt\bs\bpi^{n-1})= h^8t^3H(\us\,\ut\bpi^{n-1})
\end{multline*}
d'apr\`es les lemmes \ref{rappel} (iii), \ref{rappelA2} (iii), (i) et (iii).

Jusqu'\`a la  fin de  la preuve,  tous les  \'el\'ements de  $\uB^+$ que
nous allons  consid\'erer sont  de la  forme $y\bpi^n$.  Pour simplifier
les notations,  nous poserons $\uH(y)=(h^8t^3)^{-n}H(y\bpi^n)$ et
aussi  $\uH^i_c(\bX(y))=H^{i+8n}_c(\bX(y\bpi^n))(3n)$.  Notons
que  les lemmes
\ref{rappel}  et  \ref{rappelA2}  restent  vrais  pour  $\uH$  sauf
\ref{rappel}  (i) qui  toutefois  reste  vrai pour  $x\in  B^+$, car  on
a  alors  $H(xy\bpi^n)=H(y\bpi^n x)=H(y  x\bpi^n)$  cette
derni\`ere \'egalit\'e car $\bpi^n$ est central dans $B^+$.

\casde{\us\bt} et $\bs\ut$.  
On  a   $\uH(\us\bt)=h^{-2}t\inv \uH(\bs\us\bt)  =h^{-2}
t\inv \uH(\us\bt\bs)= h\inv t\inv \uH(\us\ut)$ par \ref{rappel} (iii), (i)
et le lemme \ref{rappelA2} (iii) ce qui ram\`ene au cas pr\'ec\'edent.
Par les lemmes \ref{rappelA2} (iv) et \ref{rappel} (i), on obtient $\uH(\bs\ut)$.

\casde{\bs\ut\bs},       $\bt\us\bt$      et       $\us\bt\bt$.      Par
$\seof{\bs\ut\bs}{\bs\ut\us}{\bs\ut}$   (o\`u   $\bs\ut\us$  est   connu
par  $\uH(\bs\ut\us)=\uH(\ut\us\bs)=h^2t\uH(\ut\us)$ en  utilisant
\ref{rappel}(i)  et  \ref{rappel}(iii))  on  trouve  la  valeur de
$\uH(\bs\ut\bs)$;  par
le lemme \ref{rappelA2}(iv) on en d\'eduit  $\bt\us\bt$ et par \ref{rappel}(i) on
en d\'eduit $\us\bt\bt$.

\casde{\bs\bt}, $\bs\bt\bs$, $\bs\bs\bt$, $\bt\bs\bt$ et $\bs\bt\bt$.
Par  les   lemmes  \ref{rappelA2}   (iv)  et   \ref{rappel}  (i)   on  a
$\uH(\bs\bt\bs)=\uH(\bs\bs\bt)= \uH(\bt\bs\bt)=\uH(\bs\bt\bt)$.
Nous  allons \'etudier  simultan\'ement les  deux \'el\'ements  $\bs\bt$
et   $\bs\bt\bs$.   Par   $\seof{\bs\bt}{\us\bt}\bt$  o\`u   la   valeur
$\uH(\bt)=\uH(\bs)$ est connue, on trouve qu'il existe
$\varepsilon\in\{0,1\}$ tel que $\uH(\bs\bt)=\varepsilon(h+h^2)+h^3t$.
Par $\seof{\bs\bt\bs}{\bs\ut\bs}{\bs\bs}$ on trouve
qu'il existe $\theta,\gamma\in\{0,1\}$ tels que
$\uH(\bs\bt\bs)=(h^2+h^3)\theta+(h^4+h^5)t^2\gamma$.
En utilisant maintenant la suite exacte longue qui r\'esulte de la proposition \ref{sws}:
\begin{align*}
\cdots&\to \uH^{i-3}_c(\bX(\bt))(-1)\to \uH^{i-2}_c(\bX(\bs\bt))(-1)\oplus
\uH^{i-1}_c(\bX(\bs\bt))\to \uH^i_c(\bX(\bs\bt\bs))\cr&\to
\uH^{i-2}_c(\bX(\bt))(-1)\to\cdots,\cr
\end{align*}
on trouve $\varepsilon=\theta=0$.
Enfin, par $\seof{\bs\bt\bs}{\us\bt\bs}{\bt\bs}$ (o\`u $\us\bt\bs$ est connu
par $\uH(\us\bt\bs)=h\uH( \us\ut)$ par le lemme \ref{rappelA2} (iii)) on trouve
$\gamma=0$.

\casde{\bt\bs\bs\bt}. On a $\uH(\bt\bs\bs\bt)=\uH(\bs\bs\bt\bt)$
par \ref{rappel} (i)
et on trouve la valeur de $\bs\bs\bt\bt$ par $\seof{\bs\bs\bt\bt}{\bs\us\bt\bt}
{\bs\bt\bt}$ o\`u le terme du milieu est connu par
$\uH(\bs\us\bt\bt)=h^2t\uH(\us\bt\bt)$
(\cf\ \ref{rappel}(iii)).

\casde{\bs\bt\bs\bt}.                        On                        a
$\uH(\bs\bt\bs\bt=\bs\bs\bt\bs)=\uH(\bs^3\bt)$  par  \ref{rappel}(i)  et
$\seof{\bs^3\bt}{\us\bs^2\bt}{\bs^2\bt}$   donne   le  r\'esultat   o\`u
$\uH(\us\bs^2\bt)=h^4t^2\uH(\us\bt)$ par \ref{rappel}(iii).

\casde{\bt\bs\bt\bs\bt}. On a le r\'esultat par $\seof{\bt\bs\bt\bs\bt}
{\ut\bs\bt\bs\bt}{\bs\bt\bs\bt}$ o\`u le terme du milieu est connu
par  $\uH(\ut\bs\bt\bs\bt=\ut\bt\bs\bt\bt)=\uH(\bt\bt\ut\bt\bs)=
h^6t^3 \uH(\ut\bs)$ par \ref{rappel}(iii).
\end{pf*}

Avant de formuler une conjecture qui g\'en\'eralise le r\'esultat ci-dessus,
faisons quelques remarques sur la conjugaison dans le groupe de tresses
de type $A_2$.
Nous appelons \og conjugaison par permutation circulaire\fg\ la relation
d'\'equi\-va\-lence sur $B^+$ cl\^oture transitive de la relation
$\bx\by\sim\by\bx$.

\begin{proposition}\label{classes  de  B(A_2)} Soient  $B$  (resp.
$B^+$)  le groupe  (resp. le  mono{\"\i}de) des  tresses de  type $A_2$.
Soit $\bw\in B^+$. Alors, il existe un unique $n\ge 0$ et
un unique $\by$ dans la liste ci-dessous tel que
$\bw$ est conjugu\'e dans $B$ \`a $\bpi^n \by$.
\begin{enumerate}
\item $\bs^a$ avec $a\ge 0$.
\item $\bs\bt$.
\item $\bs^{a_1}\bt^{a_2}\bs^{a_3}\cdots\bs^{a_{2k-1}}\bt^{a_{2k}}$,
avec $k\ge 1$, $a_i\ge 2$, la suite $(a_i)$ est plus grande pour l'ordre
lexicographique que $(a_{i+d\pmod {2k}})$ pour tout $d$.
\item $\bw_0\bs^a$ avec $a\in\{0,1\}$.
\item $\bw_0\bs^{a_1}\bt^{a_2}\bs^{a_3}\cdots\bt^{a_{2k-2}}
\bs^{a_{2k-1}}$ avec $k\ge 1$, $a_i\ge 2$
et la suite $(a_i)$ est plus grande pour
l'ordre lexicographique que $(a_{i+d\pmod {2k-1}})$ pour tout $d$.
\end{enumerate}

De plus, tout \'el\'ement de $B^+$ est conjugu\'e  par permutation
circulaire \`a un des \'el\'ements (i)--(v) ou \`a un de leurs
conjugu\'es par $\bw_0$.
\end{proposition}

\begin{pf*}{Preuve}
La conjugaison  par  $\bw_0$ r\'ealise  l'automorphisme
de   $B^+$   donn\'e   par  $\bs\mapsto\bt,   \bt\mapsto\bs$.   Il   en
r\'esulte  que tout  \'el\'ement  de $B^+$  de  la forme  $\bx\bw_0\by$
(avec $\bx,\by\in  B^+$)  est  divisible par  $\bw_0$  \`a  gauche. 

Notons pour commencer qu'un \'el\'ement de la forme
$\bt^{a_1}\bs^{a_2}\cdots\bs^{a_{2k}}$ est conjugu\'e par $\bw_0$ \`a
$\bs^{a_1}\bt^{a_2}\cdots\bt^{a_{2k}}$. Par cons\'equent,
tout \'el\'ement
$\bs^{a_1}\bt^{a_2}\cdots\bt^{a_{2k}}$ avec $a_i\ge 2$
est conjugu\'e par permutation circulaire \`a un \'el\'ement du type
(iii) ou \`a un conjugu\'e par $\bw_0$ d'un tel \'el\'ement.
De m\^eme, tout \'el\'ement de la forme
$\bw_0\bs^{a_1}\bt^{a_2}\bs^{a_3}\cdots\bt^{a_{2k-2}}
\bs^{a_{2k-1}}$ avec $k\ge 1$, $a_i\ge 2$,
est conjugu\'e par permutation circulaire \`a un \'el\'ement du type
(v) ou \`a un conjugu\'e par $\bw_0$ d'un tel \'el\'ement.

Il  est  clair  que  tout  \'el\'ement  de  $B^+$  est  conjugu\'e  par
permutation   circulaire  \`a   un  \'el\'ement   $\bw=\bw_0^n  \bx$,
o\`u aucun   des  permut\'es  circulaire de $\bw$  n'est  divisible  par
$\bw_0^{n+1}$.
Fixons un tel $\bw$.
\`A conjugaison  pr\`es   par  $\bw_0$,  l'\'el\'ement $\bx$  est  de
la  forme  $\bs^{a_1}\bt^{a_2}\bs^{a_3}\cdots\bs^{a_{2k-1}}\bt^{a_{2k}}
\bs^{a_{2k+1}}$,  o\`u   $k\ge  0$,  $a_1,  \ldots,   a_{2k}\ge  1$  et
$a_{2k+1}\ge 0$. Puisque $\bx$ n'est pas divisible par $\bw_0$,
on a trois possibilit\'es: $\bx=\bs^{a_1}$ ou $\bx=\bs^{a_1}\bt$
ou
$a_2,a_3,\ldots a_{2k-2},a_{2k-1}\ge 2$, et  $a_{2k}\ge 2$ si
$a_{2k+1}\ne 0$.

Supposons $n$  pair.
Si $\bx=\bs^{a_1}$, on est dans le cas (i). Si $\bx=\bs^{a_1}\bt$, alors
$a_1=1$, car $\bx$ est conjugu\'e circulairement \`a $\bs\bt\bs^{a_1-1}$ qui
n'est pas divisible par $\bw_0$ par hypoth\`ese. On est alors dans le cas (ii).
Sinon, quitte \`a remplacer $\bx$ par un permut\'e circulaire, on peut
supposer que $a_{2k+1}=0$. Puisque $x$ est conjugu\'e circulairement
\`a $\bs^{a_1-1}\bt^{a_2}\cdots\bs^{a_{2k-1}-1}\bs\bt^{a_{2k}}\bs$, ce
dernier n'est pas divisible \`a droite par $\bw_0$, donc
$a_{2k}\ge 2$ et on est alors dans le cas (iii).

Supposons maintenant $n$ impair.
Si $\bx=\bs^{a_1}$, on est dans le cas (v) pour $a_1\ge 2$ et dans le
cas (iv) sinon. Si $\bx=\bs^{a_1}\bt$, alors
$\bw_0\bx=\bs\bt\bs^{a_1+1}\bt$ est conjugu\'e circulairement \`a
$\bw_0\bs^{a_1+1}$ et on est dans le cas (v).
Sinon, quitte \`a remplacer $\bx$ par un permut\'e circulaire, on peut
supposer que $a_{2k+1}\not=0$. Puisque
$\bw_0\bx$ est conjugu\'e circulairement \`a
$\bw_0\bs^{a_1-1}\bt^{a_2}\cdots\bt^{a_{2k}-1}\bt\bs^{a_{2k+1}}\bt$,
alors $a_{2k+1}\not=1$ et on est dans le
cas (v), au changement $k\mapsto k+1$ pr\`es.

Il reste  \`a v\'erifier  que la  liste  donn\'ee ne  contient pas  deux
\'el\'ements  conjugu\'es. Comme  toute  conjugaison  est compos\'ee  de
\og conjugaisons \'el\'ementaires\fg\ (\cf\ \cite[Corollary 4.5]{michel})
de la forme $\bw\mapsto \bv\inv\bw\bv$
o\`u $\bv\in\bW, \bv\inv\bw\bv\in B^+$,  il suffit d'\'etudier l'effet
de la  conjugaison par  chaque \'el\'ement de  $\bW$ sur  chaque forme
donn\'ee. Il suffit en outre de consid\'erer le cas o\`u
$\bv\not\preccurlyeq \bw$, car sinon il s'agit d'une conjugaison par
permutation circulaire.
Les seuls $\bw$ \`a consid\'erer sont donc ceux des types (i), (ii) et
(iii).
On montre alors facilement que si $\bv\in\bW$, 
$\bv\not\preccurlyeq \bw$, $\bv\preccurlyeq \bw\bv$, alors
$\bv$ ou $\bw_0\bv^{-1}$ centralise $\bw$.
\end{pf*}

\begin{corollaire}\label{H   fclasse} Si  $\bG$   est  de   type
$A_2$  d\'eploy\'e,  alors  pour   tout  $i$  l'application  $\bt\mapsto
H^i_c(\bX(\bt))$   d\'efinit    une   fonction   de    classe:   $B^+\to
\CR(\bG^F\times\genby{F})$.
\end{corollaire}
\begin{pf*}{Preuve}
C'est  une cons\'equence  imm\'ediate  de la proposition
\ref{xy=yFx},  qui montre  que
si  $\bx$  et $\by$  sont  conjugu\'es  par permutation  circulaire,  on
a  $H^i_c(\bX(\bx))\simeq H^i_c(\bX(\by))$,  du lemme  \ref{rappelA2} (iv)  qui
montre  que  si  $\bx$  et  $\by$ sont  conjugu\'es  par  $\bw_0$, alors  on  a
$H^i_c(\bX(\bx))\simeq H^i_c(\bX(\by))$, et de la derni\`ere remarque de
la proposition \ref{classes de B(A_2)} qui montre que toute conjugaison dans
$B^+$ peut
\^etre  r\'ealis\'ee  par  une conjugaison  par  permutation  circulaire
suivie d'une conjugaison par $\bw_0$.
\end{pf*}

La  liste  de la proposition \ref{classes  de B(A_2)} 
forme  un  syst\`eme  de
repr\'esentants des classes de conjugaison  dans $B^+$. On d\'efinit une
fonction de classe $\varphi$ sur $B^+$  en lui attribuant la valeur $0$
dans les cas (i) et (ii), la valeur $a$ dans le cas (iv) et
la valeur $k$ dans les cas (iii) et (v)
et  en
demandant que $\varphi(\bpi^n \bb)=n+\varphi(\bb)$. Nous conjecturons le
r\'esultat suivant,  qui outre le th\'eor\`eme \ref{A2},
est \'etay\'e par  de nombreux autres calculs:

\begin{conjecture}\label{conj A2}
Supposons $\bG$ de type $A_2$ d\'eploy\'e et soit $\bb\in B^+$.
Alors, on a $H(\bb)=(-h)^{l(\bb)-\varphi(\bb)}
\Trace(T_\bb\mid R_{-ht})$.
\end{conjecture}

\sub{}
Nous allons maintenant consid\'erer un groupe tel que $(\bG,F)$ soit de
type $\lexp  2A_2$. Par cons\'equent, on a $\delta=2$,
$\rho$ est  un caract\`ere
unipotent cuspidal  et les valeurs  propres de $F^\delta$ qui  lui sont
attach\'ees  sont  dans $-q^{\delta/2}q^{\bbN}$.

\begin{theoreme}\label{2A2}
Supposons   $(\bG,F)$  de   type  $\lexp   2A_2$  et   soit  $y\in\uB^+$
apparaissant  dans  la table  ci-dessous.  Alors,  pour  tout $n$,  on  a
$H(y\bpi^n)=(h^8  t^3)^n  f(y)$ o\`u  $f(y)$  est  la valeur  de  $H(y)$
donn\'ee par la table.
\halign{$#$\hfill&\quad$#$\hfill\cr
\strut y&H(y)\cr
\noalign{\hrule}\cr
1&0\cr
\bs,\bt,\us,\ut&ht^{1/2}\cr
\us\,\ut,\ut\,\us& h^3t^{3/2} + ht^{1/2}\cr
\bs\bt,\bt\bs,\bt\bt,\bs\bs&h^3 t^{3/2}+h^2 t^{1/2}\cr
\bs\bt\bs&2 h^4t^{3/2}\cr
\us\,\ut\bs&h^2t^{1/2}\cr
\us\,\ut\bs\bs&h^3t^{1/2}+h^5t^{5/2}\cr
\bs\bs\bt\bs,\bs\bt\bs\bs,\bt\bs\bs\bt&h^5t^{3/2}+h^6t^{5/2}\cr
\bs\bs\bs\bt\bs&h^8t^{7/2}+h^6t^{3/2}\cr
\bs\bs\bt\bs\bs, \bs\bt\bs\bs\bt&h^7t^{5/2}\cr}
\end{theoreme}
\begin{pf*}{Preuve}
Comme  pour le th\'eor\`eme \ref{A2}, nous  proc\'edons  par  r\'ecurrence sur  $n$  en
traitant  en m\^eme  temps le  cas $n=0$  et le  cas g\'en\'eral,  et en
raisonnant par longueur croissante de  $y$. Mais, \`a la diff\'erence du th\'eor\`eme
\ref{A2}, nous supposons \`a  l'\'etape $n$ le th\'eor\`eme d\'emontr\'e
pour  $\bpi^n$  (si  $n=0$  cela  r\'esulte  du lemme \ref{rappelA2}(v)). La
r\'ecurrence s'ach\`evera par la  preuve  pour $\bpi^{n+1}$.

Tous  les   \'el\'ements  de   $\uB^+$  que  nous   allons  consid\'erer
sont   de    la   forme   $y\bpi^n$.    Comme   dans   le    cas   $A_2$
d\'eploy\'e,    nous   poserons    $\uH(y)=(h^8t^3)^{-n}H(y\bpi^n)$   et
$\uH^i_c(\bX(y))=H^{i+8n}_c(\bX(y\bpi^n))(3n)$.

\casde{\us\bpi^n}  et  $\us\,\ut\bpi^n$. La  preuve  dans  le cas  $A_2$
d\'eploy\'e est encore valable ici.

\casde{\bs}  et  $\bt$.  On  trouve la  valeur  de $\uH(\bs)$ par $\seof\bs\us1$,  et
$\uH(\bt)$ a m\^eme valeur par \ref{rappel}(i).

\casde{\bs\bt},     $\bs\bs$,     $\bt\bs$      et     $\bt\bt$.     Par
$\seof{\bs\bt}{\us\bt}\bt$  o\`u  le terme  du  milieu  est donn\'e  par
$\uH(\us\bt)=\uH(\bs\us)= h^2t  \uH(\us)$ par \ref{rappel}(i)  et (iii),
on trouve  la valeur de $\uH(\bs\bt)$. Par \ref{rappel}(i)  on a $\uH(\bs\bs)=\uH(\bt\bs)
=\uH(\bt\bt)=\uH(\bs\bt)$.

\casde{\us\,\ut\bs}.  On a
$h^2t\uH(\us\,\ut\bs)=\uH(\bs\us\,\ut\bs)=
\uH(\us\,\ut\bs\bt)=h\uH(\us\,\ut\us)=h^3t\uH(\us)$, respectivement par
\ref{rappel} (iii), \ref{rappel} (i), \ref{rappelA2} (iii) et \ref{rappelA2} (i).

\casde{\bs\bt\bs}, $\bs\bs\bt\bs$, $\bs\bs\bs\bt\bs$, $\bt\bs\bs\bt$ et
$\bs\bt\bs\bs$.     Par
$\seof{\bs\bt\bs}{\us\bt\bs}{\bt\bs}$  o\`u  le   terme  du  milieu  est
donn\'e par  $\uH(\us\bt\bs)=h\uH(\us\,\ut)$ par le lemme \ref{rappelA2}(iii) on
trouve  qu'il  existe  $\varepsilon\in\{0,1\}$ tel  que  $\uH(\bs\bt\bs)
=2h^4t^{3/2}+\varepsilon(h^2+h^3)t^{1/2}$.

Nous appliquons maintenant la proposition \ref{s^mb} avec  $\bb=\bt\bs$,
pour    $m=0,1,2,3$.   Avec    les    notations    de la proposition \ref{s^mb}    on
a    $\uH(\bs\bt\bs)-h   \uH(\bt\bs)=    (th^2-h)H_i$   et    avec   les
valeurs    calcul\'ees    ci-dessus    on    trouve    $\uH(\bs\bt\bs)-h
\uH(\bt\bs)=(th^2-h)h^2t^{1/2}+\varepsilon(h^2+h^3)t^{1/2}$,  qui  n'est
divisible par $th^2-h$  que si $\varepsilon=0$, ce qui est  donc le cas.
On a donc $H_i=h^2t^{1/2}$ et $H_s=h^3t^{3/2}$, ce qui donne les valeurs
pour $\bs\bs\bt\bs$ ($m=2$) et $\bs\bs\bs\bt\bs$ ($m=3$).
Par \ref{rappel}(i) on d\'eduit $\bt\bs\bs\bt$ et $\bs\bt\bs\bs$.

\casde{\us\,\ut\bs\bs}.    On   trouve    la    valeur   en    utilisant
$\seof{\us\bt\bs\bs}{\us\,\ut\bs\bs}{\us\bs\bs}$ o\`u le terme de gauche
est  donn\'e par  $\uH(\us\bt\bs\bs  )=h \uH(\us\,\ut\bs)$
(lemme \ref{rappelA2} (iii)) et celui  de
droite par $\uH(\us\bs\bs)=h^4t^2\uH(\us)$ (\ref{rappel} (iii)).

\casde{\bs\bs\bt\bs\bs}, $\bs\bt\bs\bs\bt$     et    $\bpi$.         Par
$\seof{\bs\bs\bt\bs\bs}{\bs\us\bt\bs\bs}{\bs\bt\bs\bs}$      o\`u     le
terme   du    milieu   est   donn\'e    par   $\uH(\bs\us\bt\bs\bs)=h^2t
\uH(\us\bt\bs\bs)=h^3t\uH(\us\,\ut\bs)$  par les lemmes \ref{rappel} (iii) et
\ref{rappelA2} (iii)  on      trouve       qu'il
existe    $\varepsilon'\in\{0,1\}$    tel   que    $\uH(\bs\bs\bt\bs\bs)=
h^7t^{5/2}+\varepsilon'(h^5+h^6)t^{3/2}$ et c'est aussi la valeur de
$\uH(\bs\bt\bs\bs\bt)$ par \ref{rappel} (i).

Par $\seof{\bs\bs\bt\bs\bs\bt}{\bs\bs\bt\bs\bs\ut}{\bs\bs\bt\bs\bs}$
o\`u le terme du milieu est donn\'e par $$\uH(\bs\bs\bt\bs\bs\ut)=h^4t^2
\uH(\bt\bs\bs\ut)=h^4t^2 \uH(\ut\bs\bt\bt)=h^4t^2 \uH(\us\bt\bs\bs)=h^5t^2
\uH(\us\,\ut\bs),$$ 
la premi\`ere \'egalit\'e par \ref{rappel} (i) et (iii), les autres
respectivement par les lemmes
\ref{rappel} (i), \ref{rappelA2} (iv) et \ref{rappelA2} (iii),
on trouve qu'il existe $\varepsilon''$ tel que
$\uH(\bpi)=\uH(\bs\bs\bt\bs\bs\bt)=
\varepsilon'' t^{5/2}(h^7+h^8)+\varepsilon' t^{3/2}(h^6+h^7)$.

Par la proposition \ref{sws}, on a une suite exacte longue
\begin{multline*}\cdots\to \uH^{i-3}_c(\bX(\bt\bs\bs\bt))(-1)\to
\uH^{i-2}_c(\bX(\bs\bt\bs\bs\bt))(-1)\oplus
\uH^{i-1}_c(\bX(\bs\bt\bs\bs\bt))\to\\ \uH^i_c(\bX(\bs\bs\bt\bs\bs\bt))
\to\uH^{i-2}_c(\bX(\bt\bs\bs\bt))(-1)\to\cdots\\
\end{multline*}
et on d\'eduit que $\varepsilon'\le\varepsilon''$.

Enfin, par $\seof{\bs\bs\bt\bs\bs\bs}{\us\bs\bt\bs\bs\bs}{\bs\bt\bs\bs\bs}$
o\`u le terme de gauche est donn\'e par
$\uH(\bs\bs\bt\bs\bs\bs)=\uH(\bpi)$ (\cf\ \ref{rappel}(i)),
le terme du milieu est donn\'e par
$\uH(\us\bs\bt\bs\bs\bs)=h^2t \uH(\us\bt\bs\bs\bs)=h^3t \uH(\us\,\ut\bs\bs)$
(\cf\ lemmes \ref{rappel}(iii) et \ref{rappelA2}(iii))
et le terme de droite est donn\'e par
$\uH(\bs\bt\bs\bs\bs)=\uH(\bt\bt\bt\bs\bt)=\uH(\bs\bs\bs\bt\bs)$
par les lemmes \ref{rappel} (i) et \ref{rappelA2}(iv)
respectivement, on trouve $\varepsilon''=0$.
\end{pf*}

\subsection{Type $B_2$}
\sub{}
Nous  allons maintenant  consid\'erer  un groupe  $\bG$  de type  $B_2$.
Si  $\bG$ est  d\'eploy\'e  nous noterons  $\sigma$,  $\tau$, $\rho$  et
$\theta$  les  4  caract\`eres   unipotents  de  $\GF$  diff\'erents  de
$\St$  et  $\Id$,  o\`u  $\sigma$  (resp.  $\tau$)  est  le  caract\`ere
de  la  s\'erie  principale   correspondant  au  caract\`ere  $s\mapsto
1$,  $t\mapsto  -1$   (resp.  $s\mapsto  -1$,  $t\mapsto   1$)  de  $W$,
$\rho$  est le  caract\`ere de  la s\'erie  principale correspondant  au
caract\`ere de  dimension 2  de $W$,  et $\theta$  est le  caract\`ere
unipotent  cuspidal.  Nous  \'ecrirons  donc  $H(\by)$  pour  la  partie
$(\sigma,\tau,\rho,\theta)$-isotypique de  $\sum_i h^i H_c^i(\bX(\by))$,
repr\'esent\'ee  suivant   notre  convention  \ref{convention}   par  un
\'el\'ement de $\bbZ[t^{1/2},h][\sigma,\tau,\rho,\theta]$.

Si $\bG$  est de type  $\lexp 2B_2$, alors  $F$ est une  isog\'enie dont
le  carr\'e  est  le  Frobenius  pour  une  $\bbF_{q^2}$-structure  o\`u
$q^2$  est une  puissance  impaire  de $2$.  Il  y  a deux  caract\`eres
unipotents  diff\'erents  de $\St$  et  $\Id$.  Ils sont  cuspidaux,  de
dimension  $\frac q{\sqrt  2}(q^2-1)$,  associ\'es  aux  valeurs  propres
$\lambda_\rho=\zeta_8^3=(-1+i)/\sqrt   2$   et  $\omega_\rho=q$   (resp.
$\lambda_\rho=\zeta_8^5=(-1-i)/\sqrt  2$ et  $\omega_\rho=q$) de  $F^2$,
o\`u  l'on  a pos\'e  $\zeta_8=e^{2i\pi/8}$  (\cf\  corollaire
\ref{valeurs propres}  et
\cite[\S  7.3]{LuCox}). Nous les  noterons $\rho^+$ et  $\rho^-$. Ils ne
sont  pas rationnels, mais  d\'efinis sur $\bbQ(i)$  et \'echang\'es par
l'\'el\'ement non trivial de $\Gal(\bbQ(i)/\bbQ)$. D'apr\`es la
proposition \ref{indep conj},
qui   est  applicable
car   $\frac  1{\sqrt 2}(\rho^++\rho^-)$    est   un   $R_{\tilde\chi}$
(\cf\ Proposition \ref{Rchi   indepl}),    ils   ont   m\^eme
multiplicit\'e  dans  chaque  cohomologie  et  ces  multiplicit\'es  sont
ind\'ependantes    de   $\ell$.   Nous    repr\'esenterons   la   partie
$(\rho^+,\rho^-)$-isotypique  $H(\by)$ de  $\sum_i h^i H_c^i(\bX(\by))$,
par     $P\in\bbZ[t^{1/2},h]$     repr\'esentant    $(\rho^++\rho^-)P\in
\bbZ[t^{1/2},h][\rho^+,\rho^-]$.

Puisque $\bw_0=\bs\bt\bs\bt$  est central dans $B^+$,  nous allons
trouver  une  p\'eriodicit\'e  pour  la  translation  par  $\bw_0$,  \`a
condition pour  $B_2$ d\'eploy\'e d'introduire l'involution  $E$ sur les
caract\`eres  unipotents  qui \'echange  $\sigma$  et  $\tau$ ainsi  que
$\rho$ et $\theta$ ($E$ pour \og Ennola\fg, car cette involution correspond
\`a l'\'echange de $q$ et de $-q$).

Nous commen\c cons  par \'enoncer des cons\'equences  des r\'esultats
du \S \ref{sectioncohomologie} pour $B_2$.

\begin{lemme}\label{rappelB2}
Soit $y\in\uB^+$. Alors
\begin{enumerate}
\item $H(\us\bt\bs\bt y)=hH(\underline{sts}y)$.
\item $H(\ut\,\underline{sts}\,y)=h^2tH(\ut\us y)$.
\item $H(\us\,\ut\,\us y)=h^2tH(\us y)+H(\underline{sts}y)$.

Les assertions obtenues en \'echangeant $s$ et $t$ dans les assertions
pr\'ec\'e\-dentes sont aussi vraies.
\item
Consid\'erons  les repr\'esentations  de $\CH(h^2t)(W)\sdp  F$ donn\'ees
par les valeurs suivantes sur les g\'en\'erateurs:
$$\displaylines{\sigma_{h^2t}(T_s)=h^2t\hfill\sigma_{h^2t}(T_t)=-1\hfill
\sigma(F)=1\hfill\cr
\tau_{h^2t}(T_s)=-1\hfill\tau_{h^2t}(T_t)=h^2t\hfill\tau(F)=1\hfill\cr
\rho_{h^2t}(T_s)=\begin{pmatrix}-1&0\cr h\sqrt{2t}&h^2t\cr\end{pmatrix}\hfill
\rho_{h^2t}(T_t)=\begin{pmatrix}h^2t&h\sqrt{2t}\cr0&-1\cr\end{pmatrix}\hfill
\rho_{h^2t}(F)=\begin{pmatrix}1&0\cr0&1\cr\end{pmatrix}\hfill\cr
\text{si $\bG$ est d\'eploy\'e et }\rho_{h^2t}(F)=\begin{pmatrix}0&1\cr1&0\cr
\end{pmatrix}\text{sinon}\hfill\cr}$$
Alors si  $y$ est  dans le  sous-mono{\"\i}de engendr\'e  par $\underline  W$,
on a $$\displaylines{H(y)=\frac{1}{2}\Bigl(
\sigma\cdot\Trace(T_y\mid\sigma_{h^2t}-\tau_{h^2t}+\rho_{h^2t})+
\tau\cdot\Trace(T_y\mid-\sigma_{h^2t}+\tau_{h^2t}+\rho_{h^2t})
+\hfill\cr \hfill
\rho\cdot\Trace(T_y\mid\sigma_{h^2t}+\tau_{h^2t}+\rho_{h^2t})
+\theta\cdot\Trace(T_y\mid-\sigma_{h^2t}-\tau_{h^2t}+\rho_{h^2t})
\Bigr)\cr}$$
si $\bG$ est d\'eploy\'e et 
$H(y)=(1/\sqrt 2)\Trace(T_y F\mid\rho_{h^2t})$
sinon.
\item
Si  $\bG$ est  non d\'eploy\'e  et si  $y$ est produit d'\'el\'ements de
$C_{B^+}(F)$  et d'\'el\'ements  de $\uW^F$,  alors $H(y)$ a tous ses
coefficients pairs.
\end{enumerate}
\end{lemme}
\begin{pf*}{Preuve}   (i)    est   obtenu    par   \ref{rappel}(v)    et   par
$\seof{\us\bt\bs\bt  y}{\underline{stst} y}  {\us\ut\us y}$.  (ii) vient
de la proposition \ref{fibrationws}   appliqu\'e  avec   $w=sts$  et  $y=ts$,   et  de
\ref{rappel}(v).  (iii) vient  du lemme  \ref{lemme T}  (iii) appliqu\'e
avec $w=st$ et $y=s$, et de la proposition \ref{fibrationws}. On obtient
(iv) comme cons\'equence  imm\'ediate du corollaire \ref{rationnellement
lisse} (ii)  et  du  fait  que  dans  un groupe  de  type  $B_2$  toutes  les
vari\'et\'es $\bX(\uw)$ pour $w\in  W$ sont rationnellement lisses (\cf\
proposition \ref{X  rationnellement lisse}), en utilisant  \cite[theorem 4.23 et
\S 4.15]{Lubook} si $\bG$ est d\'eploy\'e et \cite[p. 373]{Lubook} si $\bG$
n'est pas d\'eploy\'e.

D\'emontrons  le (v). La proposition \ref{indep  conj} montre que
  pour tout $y$ les
multiplicit\'es de $\rho^+$ et $\rho^-$ dans $H^j_c(\bX(y),\bbQ_\ell)_k$ sont
ind\'ependantes   de   $\ell$   pour   tout   $j$   et  tout  $k$,  o\`u
$H^j_c(\bX(y),\bbQ_\ell)_k$  est la  partie de  $H^j_c(\bX(y),\bbQ_\ell)$
o\`u les
valeurs  propres de $F^\delta$  sont de module  $q^{\delta k/2}$. Le (v)
dit    que   la   multiplicit\'e   de    $\rho^+$   et  de $\rho^-$   dans
$H^j_c(\bX(y),\bbQ_\ell)_k$ est paire. Comme $H^j_c(\bX(y),\bbQ_\ell)_k$
est  un  sous-$(\GF\times\genby  F)$-module  de  $H^j_c(\bX(y),\bbQ_\ell)$
stable  par  $\Gal(\overline\bbQ_\ell/\bbQ_\ell)$,  on  a $\Trace(gF\mid
H^j_c(\bX(y),\bbQ_\ell)_k)\in\bbQ_\ell$.             La            somme
$|\GF|\inv\sum_{g\in\GF}\Trace(gF\mid
H^j_c(\bX(y),\bbQ_\ell)_k)\rho^+(g\inv)$  appartient donc  \`a la m\^eme
extension    de   $\bbQ_\ell$    que   $\rho^+$,    c'est-\`a-dire   \`a
$\bbQ_\ell(i)$.  Comme  $y$  est  stable  par  $F$, ce dernier induit un
automorphisme   sur  la   cohomologie  dont   les  valeurs  propres  sur
$H^j_c(\bX(y),\bbQ_\ell)_k$  sont de la forme $\pm\zeta_{16}^3 q^{k/2}$,
o\`u      $\zeta_{16}^3$     est     une      racine     carr\'ee     de
$\lambda_{\rho^+}=\zeta_8^3$. La somme ci-dessus vaut donc $\zeta_{16}^3
q^{k/2}(n^+-n^-)$,  o\`u $n^+$  (resp. $n^-$)  est la  multiplicit\'e de
${\rho^+}$  dans la partie  de $H^j_c(\bX(y),\bbQ_\ell)_k$ correspondant
\`a  l'espace  propre  g\'en\'eralis\'e  de  $F$  pour  la valeur propre
$\zeta_{16}^3  q^{k/2}$ (resp. $-\zeta_{16}^3 q^{k/2}$). Comme il existe
une infinit\'e de  $\ell$  tels  que  $\zeta_{16}^3  q^{k/2}$  n'appartienne  pas  \`a
$\bbQ_\ell(i)$ (th\'eor\`eme \ref{Chebotarev}),
         cela        force         $n^+=n^-$,        donc
$\scal{H^j_c(\bX(y),\bbQ_\ell)_k}{\rho^+}$   est   pair.   Comme   cette
multiplicit\'e est ind\'ependante de $\ell$ on obtient le r\'esultat.
Le m\^eme raisonnement s'applique \`a $\rho^-$.
\end{pf*}

\sub{}
Nous d\'emontrons maintenant un th\'eor\`eme analogue \`a \ref{A2}.
\begin{theoreme}\label{B2}
Supposons  $(\bG,F)$  de  type  $B_2$ d\'eploy\'e  et  soit  $y\in\uB^+$
apparaissant  dans  la table  ci-dessous.  Alors  pour  tout $n\in\bbN$,  on  a
$H(y\bw_0^n)=(h^5 t^2)^n E^n(f(y))$ o\`u $f(y)$  est la valeur de $H(y)$
donn\'ee par la table ci-dessous (o\`u
nous avons \'etendu $E$ par lin\'earit\'e \`a
$\bbZ[t^{1/2},h][\sigma,\tau,\rho,\theta]$).
\halign{$#$\hfill&\quad$#$\hfill\cr
\strut y&H(y)\cr
\noalign{\hrule}\cr
1&\sigma+\tau+2\rho\cr
\bs&h(\tau+\rho)+h^2t(\sigma+\rho)\cr
\bt&h(\sigma+\rho)+h^2t(\tau+\rho)\cr
\us& (h^2t + 1)(\sigma+\rho)\cr
\ut& (h^2t + 1)(\tau+\rho)\cr
\bs\bs&h^2(\tau+\rho)+h^4t^2(\sigma+\rho)\cr
\bt\bt&h^2(\sigma+\rho)+h^4t^2(\tau+\rho)\cr
\us\ut,\ut\us&h^2t(\sigma+\tau+\rho+\theta)\cr
\us\bt,\bt\us&h(\rho+\sigma)+h^2t(\tau+\theta)\cr
\ut\bs,\bs\ut&h(\rho+\tau)+h^2t(\sigma+\theta)\cr
\bs\bt,\bt\bs&h^2t\theta+h^3t\rho\cr
\us\,\ut\,\us&(h^2t+h^4t^2)(\sigma+\tau+\rho+\theta)\cr
\underline{sts}&(h^2t+h^4t^2)(\tau+\theta)\cr
\underline{tst}&(h^2t+h^4t^2)(\sigma+\theta)\cr
\bs\ut\bs&h^2(\rho+\tau)+h^3t(\sigma+\theta)+h^4t^2(\sigma+\rho+\tau+\theta)\cr
\bs\bt\bs&h^3t(\theta+\sigma)+h^4t^2(\tau+\theta)\cr
\bt\bs\bt&h^3t(\theta+\tau)+h^4t^2(\sigma+\theta)\cr}
\end{theoreme}
\begin{pf*}{Preuve}
Comme  dans le th\'eor\`eme \ref{A2}, nous  proc\'edons par  r\'ecurrence
croissante sur $n$  en traitant en m\^eme  temps le cas $n=0$  et le cas
g\'en\'eral,  et en  raisonnant  par longueur  croissante  de $y$.  Nous
commen\c cons par le

\casde{\us\bw_0^n} et $\ut\bw_0^n$. Si $n=0$ la valeur se d\'eduit du lemme  
\ref{rappelB2}(iv). Sinon, on a $H(\us\bw_0^n)=th^2H(\us\bt\bs\bt\bw_0^{n-1})=
th^3H(\underline{sts}\bw_0^{n-1})=h^5t^2E(H(\us\bw_0^{n-1}))$ 
respectivement  par les lemmes \ref{rappel}(iii), \ref{rappelB2}(i)  et  les
valeurs donn\'ees dans la table. Le cas de $\ut\bw_0^n$ est analogue.

\casde{\bw_0^n}, $\bs\bw_0^n$, $\bs\bs\bw_0^n$, $\bt\bw_0^n$ et $\bt\bt\bw_0^n$.
Si $n=0$ la valeur se d\'eduit par exemple de la proposition \ref{H_s}. Sinon,
nous appliquons la proposition \ref{s^mb} avec $\bb=\bt\bs\bt\bw_0^{n-1}$ pour $m=0,1,2,3$.
Par l'hypoth\`ese de r\'ecurrence nous avons $$H(\bb)=(h^5t^2)^{n-1}E^{n-1}
(h^3t(\theta+\tau)+h^4t^2(\sigma+\theta)).$$
Soit $\bY$ le ferm\'e compl\'ementaire de $\bX(\bs\bb)=\bX(\bw_0^n)$
dans $\bX(\uw_0\bw_0^{n-1})$. Par
$\seof{\bs\bt\bs\bw_0^{n-1}}\bY{\ut\us\ut\bw_0^{n-1}}$ et 
$\seof{\bt\bs\bt\bw_0^{n-1}}\bY{\us\ut\us\bw_0^{n-1}}$
on trouve qu'il existe $\varepsilon\in\{0,1\}$ tel que
$H(\bY)=(h^5t^2)^{n-1}\allowbreak E^{n-1}(h^4t^2(\tau+\sigma+2\theta)
+\varepsilon(h^2+h^3)t\theta)$.
Par  $\seof{\bw_0^n}{\uw_0\bw_0^{n-1}}{\bY}$  et \ref{rappel}(v)  on  en
d\'eduit  $H(\bw_0^n)=(h^5t^2)^{n-1}E^{n-1}(h^5t^2(\tau+\sigma+2\theta)+
\varepsilon(h^3+h^4)t\theta)$.   Avec   les  notations   de la proposition \ref{s^mb}
on     doit     avoir      $H(\bs\bb)-hH(\bb)=(ht^2-h)H_i$  d'o\`u
$H(\bs\bb)-hH(\bb)=(h^5t^2)^{n-1}E^{n-1}((ht^2-h)h^3t(\theta+\tau)
+\varepsilon(h^3+h^4)t   \theta)$,   ce    qui   n'est   divisible   par
$ht^2-h$   que   si   $\varepsilon=0$.  Donc   $\varepsilon=0$   et   on
trouve   les   valeurs  $H_i=h^{-2}t\inv(h^5t^2)^nE^n(\rho+\sigma)$   et
$H_s=h\inv(h^5t^2)^nE^n(\rho+\tau)$.   En  appliquant la proposition   \ref{s^mb}  avec
$m=2,3$ on en d\'eduit les valeurs pour $\bs\bw_0^n$ et $\bs\bs\bw_0^n$.
On   proc\`ede   de  fa\c   con   sym\'etrique   pour  $\bt\bw_0^n$   et
$\bt\bt\bw_0^n$, en \'echangeant les r\^oles de $\bs$ et $\bt$.

\casde{\us\ut}. Par  les lemmes \ref{rappel} (iii), \ref{rappelB2}(i) et
(ii), on a $H(\us\ut\bw_0^n)=
H(\us\ut\bt\bs\bt\bs\bw_0^{n-1})=
h^2tH(\us\ut\bs\bt\bs\bw_0^{n-1})=
h^3tH(\us\,\ut\us\ut\bw_0^{n-1})=
h^5t^2H(\us\ut\bw_0^{n-1})$,
ce qui  est  le  r\'esultat cherch\'e  car  la  table montre  que
$E(H(\us\ut\bw_0^{n-1}))=H(\us\ut\bw_0^{n-1})$.

\casde{\us\,\ut\,\us}. On proc\`ede exactement comme le cas pr\'ec\'edent.

\casde{\us\ut\us} et $\ut\us\ut$. On a
$H(\us\ut\us\bw_0^n)=
H(\us\,\ut\,\us\bw_0^n)-h^2tH(\us\bw_0^n)=\linebreak
h^5t^2E(H(\us\,\ut\,\us\bw_0^{n-1}))-h^7t^3E(H(\us\bw_0^{n-1}))=
h^5t^2E(H(\us\ut\us\bw_0^{n-1}))$   par le lemme \ref{rappelB2}(iii),   les  cas
pr\'ec\'edents,   et  \`a   nouveau le lemme \ref{rappelB2}(iii).   Le  cas   de
$\ut\us\ut$ est analogue.

Jusqu'\`a la fin de la preuve, tous les \'el\'ements de $\uB^+$ que nous
allons  consid\'erer  sont  de  la  forme  $y\bw_0^n$.  Pour  simplifier
les  notations,  nous posons  $\uH(y)=(h^5t^2)^{-n}E^n(H(y\bw_0^n))$  et
$\uH^i_c(\bX(y))=E^n(H^{i+5n}_c(\bX(y\bw_0^n))(2n))$.

\casde{\us\bt} et $\bs\ut$.
En utilisant $\uH(\us\bt\bs)=\uH(\bs\us\bt)=h^2t\uH(\us\bt)$ et
$\seof{\us\bt\bs}{\us\ut\us}{\us\ut}$ on trouve qu'il existe
$\varepsilon,\varepsilon'\in\{0,1\}$ tels que
$\uH(\us\bt)=h(\rho+\sigma)+h^2t(\tau+\theta)+(1+h)(\varepsilon\tau
+\varepsilon'\theta)$. En reportant ces valeurs dans la suite exacte
longue d\'eduite de $\seof{\us\bt}{\us\ut}{\us}$
on trouve $\varepsilon=\varepsilon'=0$.

On proc\`ede sym\'etriquement pour $\bs\ut$.

\casde{\bs\ut\bs}. On trouve le r\'esultat par
$\seof{\bs\ut\bs}{\bs\ut\us}{\bs\ut}$ o\`u le terme du milieu est donn\'e par
$\uH(\bs\ut\us)=\uH(\ut\us\bs)=h^2t\uH(\ut\us)$.

\casde{\bs\bt}, $\bt\bs$, $\bs\bt\bs$ et $\bs\bs\bt$.
Par $\seof{\bs\bt}{\us\bt}{\bt}$ on trouve qu'il existe
$\varepsilon,\varepsilon',\varepsilon''\in\{0,1\}$ tels que
$\uH(\bs\bt)=h^2t\theta+h^3t\rho+\varepsilon(h+h^2)\rho+
\varepsilon'(h+h^2)\sigma+\varepsilon''(h^2+h^3)t\tau$.
Sym\'etriquement, par $\seof{\bs\bt}{\bs\ut}{\bs}$ on trouve qu'il existe
$\eta,\eta',\eta''\in\{0,1\}$ tels que
$\uH(\bs\bt)=h^2t\theta+h^3t\rho+\eta(h+h^2)\rho+
\eta'(h+h^2)\tau+\eta''(h^2+h^3)t\sigma$. En comparant,
on trouve $\varepsilon=\eta$ et $\varepsilon'=\varepsilon''=\eta'=\eta''=0$.

\`A ce stade on a donc (la premi\`ere \'egalit\'e par \ref{rappel} (i))
$\uH(\bs\bt)=\uH(\bt\bs)
=h^2t\theta+h^3t\rho+\varepsilon(h+h^2)\rho$.
Par $\seof{\bs\bt\bs}{\us\bt\bs}{\bt\bs}$ et $\uH(\us\bt\bs)=h^2t\uH(\us\bt)$
(\ref{rappel}(iii)), on trouve qu'il existe
$\delta\in\{0,1\}$ tel que
$\uH(\bs\bt\bs)=\varepsilon(h^2+h^3)\rho+\delta(h^3+h^4)t \rho+
h^3t(\sigma+\theta)+h^4t^2(\tau+\theta)$. Par
$\seof{\bs\bt\bs}{\bs\ut\bs}{\bs\bs}$ on voit que $\delta=0$.
Enfin, en utilisant la suite exacte longue qui vient de la proposition \ref{sws}
\begin{align*}\cdots&\to \uH^0_c(\bX(\bt))(-1)\to
\uH^1_c(\bX(\bs\bt))(-1)\oplus
\uH^2_c(\bX(\bs\bt))\to \uH^3_c(\bX(\bs\bs\bt))\cr&\to
\uH^1_c(\bX(\bt))(-1)\to\cdots\cr
\end{align*}
(o\`u $\uH(\bs\bs\bt)=\uH(\bs\bt\bs)$ par \ref{rappel} (i))
on trouve $\varepsilon=0$.
On proc\`ede sym\'etri\-quement pour $\bt\bs\bt$.
\end{pf*}

\begin{theoreme}\label{2B2}
Supposons $(\bG,F)$ de type $\lexp 2B_2$ et soit $y\in\uB^+$
apparaissant dans la table ci-dessous. Alors pour
tout $n$, on a $H(y\bw_0^n)=(h^5 t^2)^nf(y)$ o\`u $f(y)$ est la valeur
de $H(y)$ donn\'ee par la table.
\halign{$#$\hfill&\quad$#$\hfill\cr
\strut y&H(y)\cr
\noalign{\hrule}\cr
1&0\cr
\bs,\bt,\us,\ut&h t^{1/2}\cr
\bt\bs,\bs\bt,\bs\bs,\bt\bt&h^2t^{1/2}+h^3t^{3/2}\cr
\us\ut\us& h^3t^{3/2}\cr
\us\bt\bs& h^2t^{1/2}\cr
\bs\bt\bs, \bt\bs\bt&h^4t^{3/2}\cr
\us\,\ut\,\us&2h^3t^{3/2}\cr}
\end{theoreme}
\begin{pf*}{Preuve}
La d\'emonstration proc\`ede comme dans le th\'eor\`eme \ref{2A2}.
\casde{\us\bw_0^n}, $\us\,\ut\,\us\bw_0^n$, $\us\ut\us\bw_0^n$.
On proc\`ede exactement comme
pour les m\^emes cas dans le th\'eor\`eme \ref{B2}.

\casde{\us\bt\bs\bw_0^n}.
On a $hH(\us\ut\us\bw_0^n)=H(\us\bt\bs\bt\bw_0^n)=H(\us\bt\bs\bw_0^n\bt)
=H(\bs\us\bt\bs\bw_0^n)=h^2t  H(\us\bt\bs\bw_0^n)$  o\`u  la  premi\`ere
\'egalit\'e est  par le lemme \ref{rappelB2} (i),  la deuxi\`eme car  $\bw_0$ est
central,  la  troisi\`eme par  \ref{rappel}  (i)  et la  derni\`ere  par
\ref{rappel}(iii).

Dans  la  suite,  les  \'el\'ements de  $\uB^+$  consid\'er\'es  sont  de
la  forme  $y\bw_0^n$.  Pour   simplifier  les  notations,  nous  posons
$\uH(y)=(h^5t^2)^{-n}H(y\bw_0^n)$ et
$\uH^i_c(\bX(y))=H^{i+5n}_c(\bX(y\bw_0^n))_{\rho^+}(2n)
=H^{i+5n}_c(\bX(y\bw_0^n))_{\rho^-}(2n)$.

\casde{\bs}. On conclut par $\seof\bs\us1$.

\casde{\bs\bt},  $\bs\bs$, $\bt\bs$,  $\bt\bt$. On  trouve $\bs\bs$  par
$\seof{\bs\bs}{\bs\us}\bs$  o\`u  le terme  du  milieu  est donn\'e  par
\ref{rappel}(iii).  Les  autres  \'el\'ements sont  $F$-conjugu\'es  par
permutation circulaire (\cf\ \ref{rappel} (i)).

\casde{\bs\bt\bs} et $\bw_0$. Par $\seof{\bs\bt\bs}{\us\bt\bs}{\bt\bs}$
on trouve qu'il existe $\varepsilon\in\{0,1\}$ tel que
$\uH(\bs\bt\bs)=h^4t^{3/2}+\varepsilon(h^2+h^3)t^{1/2}$. Par
$\seof{\bw_0}{\us\bt\bs\bt}{\bt\bs\bt}$, o\`u le terme du milieu est
\'egal \`a $h^2t\uH(\us\bt\bs)$ par \ref{rappel} (iii), on trouve qu'il existe
$\eta\in\{0,1\}$ tel que $\uH(\bw_0)=\varepsilon(h^3+h^4)t^{1/2}+
\eta(h^4+h^5)t^{3/2}$.  On conclut en utilisant le lemme \ref{rappelB2}(v) qui montre
qu'on doit avoir $\varepsilon=\eta=0$.
\end{pf*}

\subsection{Type $G_2$}
\sub{}
Nous consid\'erons maintenant un groupe $\bG$  de type $G_2$. 

Si $\bG$ est d\'eploy\'e, nous notons $\sigma$ (resp. $\tau$) le caract\`ere
lin\'eaire  de  $W$  donn\'e  par $s\mapsto  -1$,  $t\mapsto  1$  (resp.
$s\mapsto  1$,  $t\mapsto  -1$)  et  nous  notons  $A$  (resp.  $B$)  la
repr\'esentation  irr\'eductible de  degr\'e  2 o\`u $w_0$ agit par  $-1$
(resp. 1). Par abus  de notation, nous noterons  de la m\^eme
fa\c con  les caract\`eres  unipotents de  la s\'erie  principale autres
que  $\Id$ et  $\St$.  Il  y a  4  autres  caract\`eres unipotents,  qui
sont  cuspidaux  et  associ\'es  aux  valeurs  propres  $\lambda_\rho=1,
-1,  j,   j^2$  respectivement  o\`u  l'on   a  pos\'e  $j=e^{2i\pi/3}$.
Nous   les  noterons   $\rho,\rho_{-1},\rho_j$   et  $\rho_{j^2}$.
Nous ne  saurons  d\'emontrer  un  th\'eor\`eme  que  pour  les  caract\`eres
diff\'erents  de  $B$  et   $\rho_{-1}$.  Nous  \'ecrirons  donc  $H(y)$
pour  la  partie $\sigma,\tau,A,\rho,  \rho_j,\rho_{j^2}$-isotypique  de
$\sum_i  h^i  H^i_c(\bX(y))$,  repr\'esent\'ee  par  un  \'el\'ement  de
$\bbZ[t,h][\sigma,\tau,A,\rho,  \rho_j,\rho_{j^2}]$.

Si  $(\bG,F)$   est  de   type  $\lexp{2}\bG_2$,   alors  $F$   est  une
isog\'enie  dont  le  carr\'e  est  l'endomorphisme  de  Frobenius  pour
une  $\bbF_{q^2}$-structure  o\`u $q^2$  est  une  puissance impaire  de
$3$. Le  groupe $\GF$  poss\`ede 6 caract\`eres  unipotents diff\'erents
de  $\St$  et  $\Id$,  tous   cuspidaux.  Ils  ne  sont  pas  rationnels
mais  d\'efinis  sur  $\bbQ(\zeta_3)$  o\`u  $\zeta_3=e^{2i\pi/3}$.  Ils
forment  3   orbites  sous  $\Gal(\bbQ(\zeta_3)/\bbQ)$.   La  premi\`ere
contient  deux caract\`eres  de dimension  $\frac q{ 2\sqrt  3}(q^2+\sqrt
3q+1)(q^2-1)$,  associ\'es  aux   valeurs  propres  $\lambda_\rho=i$  et
$\omega_\rho=q$  de $F^2$  (resp. $\lambda_\rho=-i$  et $\omega_\rho=q$)
(\cf\ corollaire \ref{valeurs propres} et \cite[4.7]{GM}); nous noterons
$\rho_i$  et  $\rho_{-i}$ ces  caract\`eres.  La  deuxi\`eme orbite  est
identique sauf  que la dimension  des caract\`eres qu'elle  contient est
$\frac q{  2\sqrt  3}(q^2-\sqrt  3q+1)(q^2-1)$; nous  noterons  $\rho'_i$
et   $\rho'_{-i}$  les   caract\`eres  correspondants.   La  troisi\`eme
orbite   contient   deux   caract\`eres   de   dimension   $\frac q{\sqrt
3}(q^4-1)$, associ\'es  aux valeurs  propres $\lambda_\rho=\zeta_{12}^5$
et   $\omega_\rho=q$   de   $F^2$   (resp.   $\lambda_\rho=\zeta_{12}^7$
et  $\omega_\rho=q$)   o\`u  $\zeta_{12}=e^{2i\pi/12}$.   Nous  noterons
$\rho_{\zeta_{12}^5}$  et $\rho_{\zeta_{12}^7}$  ces deux  caract\`eres.
Nous noterons $H(\by)$  la partie $(\rho_i,\rho_{-i},\rho'_i,\rho'_{-i},
\rho_{\zeta_{12}^5},\rho_{\zeta_{12}^7})$-isotypique   de  $\sum_i   h^i
H_c^i     (\bX(\by))$.     Nous     repr\'esenterons     $H(\by)$    par
$P\in\bbZ[t^{1/2},h]    [\rho_i,    \rho_{-i},    \rho'_i,   \rho'_{-i},
\rho_{\zeta_{12}^5}, \rho_{\zeta_{12}^7}]$.

Ici  encore $\bw_0=\bs\bt\bs\bt\bs\bt$  est central  dans $B^+$, et nous
allons  trouver une p\'eriodicit\'e pour la translation par $\bw_0$, \`a
\og l'involution  d'Ennola\fg\ $E$ pr\`es (qui correspond \`a l'\'echange de
$q$ et de $-q$) qui dans le cas d\'eploy\'e \'echange $B$ et $\rho_{-1}$
d'une  part, et $A$ et $\rho$ d'autre part. Dans le cas non d\'eploy\'e,
$E$  \'echange $\rho_i$ et  $\rho'_i$, et $\rho_{-i}$  et $\rho'_{-i}$.

Nous commen\c cons  par rappeler des cons\'equences  des r\'esultats des
deux premiers chapitres pour un groupe de type $G_2$.

\begin{lemme}\label{rappelG2} Soit $y\in\uB^+$. Alors
\begin{enumerate}
\item $H(\us\bt\bs\bt\bs\bt y)=h H(\underline{ststs}y)$.
\item $H(\us\ \underline{ts}y)=H(\underline{sts}y)+h^2t H(\us y)$.
\item $H(\ut\ \underline{sts}y)=H(\underline{tsts}y)+h^2t
                                                    H(\underline{ts}y)$.
\item $H(\us\ \underline{tsts}y)=H(\underline{ststs}y)+h^2t
                                                    H(\underline{sts}y)$.
\item $H(\ut\ \underline{ststs}y)=h^2tH(\underline{tsts}y)$.
\par
Les sym\'etriques des assertions ci-dessus obtenues en \'echangeant $s$ et $t$
sont aussi vraies.
\item Consid\'erons les repr\'esentations $A_{h^2t}$
et $B_{h^2t}$ de $\CH(h^2t)\sdp F$ d\'efinies
sur les g\'en\'erateurs par:
$$\displaylines{\hfill
A_{h^2t}(T_s)=\begin{pmatrix}-1&0\cr h\sqrt{3t}&h^2t\cr\end{pmatrix}\qquad
A_{h^2t}(T_t)=\begin{pmatrix}h^2t&h\sqrt{3t}\cr0&-1\cr\end{pmatrix}
\hfill\cr
\hfill
B_{h^2t}(T_s)=\begin{pmatrix}-1&0\cr h\sqrt{t}&h^2t\cr\end{pmatrix}\qquad
B_{h^2t}(T_t)=\begin{pmatrix}h^2t&h\sqrt{t}\cr0&-1\cr\end{pmatrix}
\hfill\cr}$$
et dans les deux cas
$A_{h^2t}(F)=B_{h^2t}(F)=\begin{pmatrix}1&0\cr0&1\cr\end{pmatrix}$
si $\bG$ est d\'eploy\'e  et
$\begin{pmatrix}0&1\cr1&0\cr\end{pmatrix}$ sinon.
Alors si $y$ est dans le sous-mono{\"\i}de engendr\'e par $\uW$, on a
si $\bG$ est d\'eploy\'e:
\begin{multline*} 6H(y)=
A\Trace(T_y\mid A_{h^2t}+3B_{h^2t}+2\sigma_{h^2t}+2\tau_{h^2t})\\
+\sigma\Trace(T_y\mid 2A_{h^2t}+4\sigma_{h^2t}-2\tau_{h^2t})
+\tau\Trace(T_y\mid 2A_{h^2t}-2\sigma_{h^2t}+4\tau_{h^2t}) \\
+\rho\Trace(T_y\mid A_{h^2t}-B_{h^2t}+2\sigma_{h^2t}+2\tau_{h^2t})
+(\rho_j+\rho_{j^2})\Trace(T_y\mid 2A_{h^2t}-2\sigma_{h^2t}-2\tau_{h^2t})\\
\end{multline*}
et sinon
\begin{multline*}
H(y)=\frac{1}{2}\bigl((\rho_i+\rho_{-i})\Trace(T_yF\mid A_{h^2t}/\sqrt 3
+B_{h^2t})\\
   +(\rho'_i+\rho'_{-i})\Trace(T_yF\mid A_{h^2t}/\sqrt 3-B_{h^2t})\bigr)
   +(\rho_{\zeta_{12}^5}+\rho_{\zeta_{12}^7})\Trace(T_yF\mid A_{h^2t}/\sqrt 3).
\end{multline*}
\item
Supposons $\bG$ non d\'eploy\'e et $y$ produit d'\'el\'ements de
$C_{B^+}(F)$  et  d'\'el\'ements  de  $\uW^F$. 
Soit  $H^j_c(y,\bbQ_\ell)_k$ la  partie de $H^j_c(y,\bbQ_\ell)$ o\`u
les valeurs propres de $F^\delta$ sont de module $q^{\delta k/2}$.
Si  les
multiplicit\'es  de  tous  les  caract\`eres  unipotents  cuspidaux dans
$H^j_c(y,\bbQ_\ell)_k$  sont  ind\'ependantes  de  $\ell$ pour tout $j$,
alors tous les coefficients de $H(y)$ sont pairs.
\end{enumerate}
\end{lemme}
\begin{pf*}{Preuve}   (i)    est   obtenu    par   \ref{rappel}(v)    et   par
$\seof{\us\bt\bs\bt\bs\bt   y}{\underline{ststst}   y}{\underline{ststs}
y}$. (ii), (iii) et  (iv) viennent du lemme  \ref{lemme T} (iii) et
de  la  proposition  \ref{fibrationws}. Pour   (v), on obtient \`a partir du
lemme \ref{lemme T} (iii) et de  la  proposition  \ref{fibrationws}
$$H(\ut\ \underline{ststs}\ y)=H(\uw_0 y)+h^2tH(\underline{tsts} y).$$
En  utilisant  de  plus le lemme \ref{rappel}(v) on obtient le r\'esultat.
On  obtient  (vi) imm\'ediatement  par  le  corollaire
\ref{rationnellement  lisse} (ii) en  utilisant  \cite[p. 376]{Lubook}  (voir
aussi \cite[4.7]{GM}) gr\^ace  au fait que dans un groupe  de type $G_2$
toutes les  vari\'et\'es $\bX(\uw)$  pour $w\in W$  sont rationnellement
lisses (\cf\ proposition \ref{X rationnellement lisse}).

Le  (vii)  dit  que  la  multiplicit\'e  de  tout  caract\`ere unipotent
cuspidal   $\rho$  dans  $H^j_c(\bX(y),\bbQ_\ell)_k$  est  paire.  Comme
$H^j_c(\bX(y),\bbQ_\ell)_k$ est un sous-$(\GF\times\genby F)$-module de 
$H^j_c(\bX(y),\bbQ_\ell)$ stable par $\Gal(\overline\bbQ_\ell/\bbQ_\ell)$,
on a $\Trace(gF\mid H^j_c(\bX(y),\bbQ_\ell)_k)\in\bbQ_\ell$. La somme
$|\GF|\inv\sum_{g\in\GF}\Trace(gF\mid H^j_c(\bX(y),\bbQ_\ell)_k)\rho(g\inv)$
appartient donc \`a  la   m\^eme  extension   de  $\bbQ_\ell$  que  $\rho$,
c'est-\`a-dire \`a $\bbQ_\ell(\zeta_3)$. Comme $y$ est stable par $F$,
ce  dernier induit un automorphisme sur  la cohomologie dont les valeurs
propres  sur $H^j_c(\bX(y),\bbQ_\ell)_k$  sont de  la forme $\pm\mu_\rho
q^{k/2}$,  o\`u $\mu_\rho$ est une racine carr\'ee de $\lambda_\rho$. La
somme ci-dessus vaut donc $\mu_\rho q^{k/2}(n^+-n^-)$, o\`u $n^+$ (resp.
$n^-$)   est   la   multiplicit\'e   de   $\rho$   dans   la  partie  de
$H^j_c(\bX(y),\bbQ_\ell)_k$    correspondant   \`a   l'espace   propre
g\'en\'eralis\'e  de $F$ pour la valeur propre $\mu_\rho q^{k/2}$ (resp.
$-\mu_\rho  q^{k/2}$).  D'apr\`es le th\'eor\`eme \ref{Chebotarev},
il  existe  des $\ell$ tels que $\mu_\rho
q^{k/2}$   n'appartienne  pas   \`a  $\bbQ_\ell(\zeta_3)$ (car $\mu_\rho$
est une racine primitive 8-i\`eme ou 24-i\`eme de l'unit\'e), donc
$n^+=n^-$  et $\scal{H^j_c(\bX(y),\bbQ_\ell)_k}\rho$  est pair. Comme
cette   multiplicit\'e  est  ind\'ependante  de  $\ell$  on  obtient  le
r\'esultat.
\end{pf*}

On  notera  que  dans  la   table  ci-dessous,  l'\'echange  de  $s$  et
$t$  dans  une   vari\'et\'e  \'echange  $\sigma$  et   $\tau$  dans  la
cohomologie  correspondante.  En   caract\'eristique  3,  ce  r\'esultat
peut  se d\'emontrer  en  utilisant l'isog\'enie  de  $\bG$ induite  par
l'automorphisme du diagramme; dans les autres caract\'eristiques nous ne
savons que le constater.
\begin{theoreme}\label{G2}
Supposons  $(\bG,F)$  de  type  $G_2$ d\'eploy\'e  et  soit  $y\in\uB^+$
apparaissant  dans  la table  ci-dessous.  Alors  pour  tout $n\in \bbN$,  on  a
$H(y\bw_0^n)=(h^7 t^3)^n E^n(f(y))$ o\`u $f(y)$  est la valeur de $H(y)$
donn\'ee par la table (dans la table
nous avons pos\'e $J=\rho_j+\rho_{j^2}$ et nous avons
\'etendu $E$ par lin\'earit\'e \`a $\bbZ[t^{1/2},h][\sigma,\tau,A,\rho,J]$).
\halign{$#$\hfill&\quad$#$\hfill\cr
\strut y&H(y)\cr
\noalign{\hrule}\cr
1&\sigma+\tau+2A\cr
\bs&h(\sigma+A)+h^2t(\tau+A)\cr
\bt&h(\tau+A)+h^2t(\sigma+A)\cr
\us&(h^2t+1)(\tau+A)\cr
\ut&(h^2t+1)(\sigma+A)\cr
\bs\bt, \bt\bs&h^3t A+h^2tJ\cr
\us\bt,\bt\us&h(\tau+A)+h^2t(\sigma+J)\cr
\ut\bs,\bs\ut&h(\sigma+A)+h^2t(\tau+J)\cr
\us\ut,\ut\us&h^2t(\sigma+\tau+A+J)\cr
\bs\bs&h^2(\sigma+A)+h^4t^2(\tau+A)\cr
\bt\bt&h^2(\tau+A)+h^4t^2(\sigma+A)\cr
\bs\bt\bs&h^3t(\tau+J)+h^4t^2(\sigma+J)\cr
\bt\bs\bt&h^3t(\sigma+J)+h^4t^2(\tau+J)\cr
\bs\ut\bs&h^4t^2(\sigma+\tau+A+J)+h^3t(\tau+J)+h^2(\sigma+A)\cr
\underline{sts}&(h^2t+h^4t^2)(\sigma+J)\cr
\us\bt\bs&h^3t(\tau+A)+h^4t^2(\sigma+J)\cr
\bs\bt\bs\bt&h^4t^2\rho+h^5t^2J\cr
\underline{tsts},\underline{stst}&h^4t^2(\tau+\sigma+\rho+J)\cr
\us\bt\bs\bt&h^3t(\sigma+J)+h^4t^2(\tau+\rho)\cr
\ut\bs\bt\bs&h^3t(\tau+J)+h^4t^2(\sigma+\rho)\cr
\underline{ststs}&(h^6t^3+h^4t^2)(\tau+\rho)\cr
\bt\bs\bt\bs\bt&h^6t^3(\sigma+\rho)+h^5t^2(\tau+\rho)\cr
\ut\bs\bt\bs\bt&h^5t^2(\tau+J)+h^6t^3(\sigma+\rho)\cr
\us\bt\bs\bt\bs&h^5t^2(\sigma+J)+h^6t^3(\tau+\rho)\cr}
\end{theoreme}
\begin{pf*}{Preuve}
Nous proc\'edons comme pour le th\'eor\`eme \ref{2A2}, en d\'emontrant le th\'eor\`eme par
r\'ecurrence sur  $n$; \`a  l'\'etape $n$  nous le  supposons d\'emontr\'e
pour  $y\bw_0^{n-1}$,  o\`u  $y$  est  dans la  table,  ainsi  que  pour
$\bw_0^n$;ainsi  nous terminons  la d\'emonstration  en le  d\'emontrant
pour $\bw_0^{n+1}$.

\casde{\us\bw_0^n} et $\ut\bw_0^n$. 
Si $n=0$ la valeur se d\'eduit du lemme \ref{rappelG2}(vi). Sinon, on a
$H(\us\bw_0^n)=H(\us\bs\bt\bs\bt\bs\bt\bw_0^{n-1})=
h^2t H(\us\bt\bs\bt\bs\bt\bw_0^{n-1})=
h^3t H(\underline{ststs}\bw_0^{n-1})=h^7t^3
E(H(\us\bw_0^{n-1}))$, en utilisant les lemmes \ref{rappel}(iii),
\ref{rappelG2}(i) et la r\'ecurrence.
Le cas de $\ut\bw_0^n$ est analogue.
\casde{\us \ut\bw_0^n} et de $\ut\us\bw_0^n$.
Si $n=0$ la valeur se d\'eduit du lemme \ref{rappelG2}(vi). Sinon, on a
\begin{multline*}
H(\ut\us\bw_0^n)=H(\ut\us\bs\bt\bs\bt\bs\bt\bw_0^{n-1})=
h^2t H(\ut\us\bt\bs\bt\bs\bt\bw_0^{n-1})=
h^3t H(\ut\ \underline{ststs}\bw_0^{n-1})=\\
=h^5t^2H(\underline{tsts}\bw_0^{n-1})
=h^7t^3 E(H(\ut\us\bw_0^{n-1}))
\end{multline*}
par les lemmes \ref{rappel}(iii), \ref{rappelG2}(i), \ref{rappelG2}(v) et la
r\'ecurrence.
Le cas de $\us\ut\bw_0^n$ est analogue.
\casde{\underline{sts}\bw_0^n}.
Si $n=0$ la valeur se d\'eduit du lemme \ref{rappelG2}(vi). Sinon, on a
$H(\underline{sts}\bw_0^n)=$(par le lemme \ref{rappelG2}(ii))
$H(\underline{st}\ \us\bw_0^n)-h^2tH(\us\bw_0^n)=$
(voir la preuve du cas $\us\bw_0^n$)
$h^3t H(\underline{st}\ \underline{ststs}\bw_0^{n-1})
-h^5t^2H(\underline{ststs}\bw_0^{n-1})=$(par le lemme \ref{rappelG2}(v))
$h^5t^2 H(\us\ \underline{tsts}\bw_0^{n-1})
-h^5t^2 H(\underline{ststs}\bw_0^{n-1})=$(par le lemme \ref{rappelG2}(iv))
$h^7t^3 H(\underline{sts}\bw_0^{n-1})=$(par la table)
$h^7t^3 E(H(\underline{sts}\bw_0^{n-1}))$. Le cas de $\underline{tst}\bw_0^n$
est analogue
\casde{\underline{tsts}\bw_0^n}.
Si $n=0$ la valeur se d\'eduit du lemme \ref{rappelG2}(vi). Sinon, on a
\begin{equation*}
\begin{split}
H(\underline{tsts}\bw_0^n)&=
H(\underline{tst}\ \us\bw_0^n)-h^2tH(\ut\us\bw_0^n)
\text{ (par le lemme \ref{rappelG2}(iii)) }\\
&=h^3tH(\underline{tst}\ \underline{ststs}\bw_0^{n-1})
-h^7t^3H(\underline{tsts}\bw_0^{n-1})\\&\phantom{=}
\text{ (voir les preuves des cas $\us\bw_0^n$ et $\ut\us\bw_0^n$) }\\
&=h^3tH(\underline{ts}\ \ut\ \underline{ststs}\bw_0^{n-1})
-h^5t^2H(\ut\ \underline{ststs}\bw_0^{n-1})
-h^7t^3H(\underline{tsts}\bw_0^{n-1})\\&\phantom{=}
\text{ (par le lemme \ref{rappelG2}(ii)) }\\
&=h^5t^2H(\underline{ts}\ \underline{tsts}\bw_0^{n-1})
-2h^7t^3H(\underline{tsts}\bw_0^{n-1})\\&\phantom{=}
\text{ (par le lemme \ref{rappelG2}(v)) }\\
&=h^5t^2H(\ut\ \underline{ststs}\bw_0^{n-1})
+h^7t^3H(\ut\ \underline{sts}\bw_0^{n-1})
-2h^7t^3H(\underline{tsts}\bw_0^{n-1})\\&\phantom{=}
\text{ (par le lemme \ref{rappelG2}(iv)) }\\
&=h^7t^3H(\ut\ \underline{sts}\bw_0^{n-1})
-h^7t^3H(\underline{tsts}\bw_0^{n-1})
\text{ (par le lemme \ref{rappelG2}(v)) }\\
&=h^9t^4H(\underline{ts}\bw_0^{n-1})
\text{ (par le lemme \ref{rappelG2}(iii)) }\\
&=h^7t^3E(H(\underline{tsts}\bw_0^{n-1}))\text{ (par la table) }.
\end{split}
\end{equation*}
Le cas de $\underline{stst}\bw_0^n$ est analogue.
\casde{\underline{ststs}\bw_0^n}
Si $n=0$ la valeur se d\'eduit du lemme \ref{rappelG2}(vi). Sinon, on a
\begin{equation*}
\begin{split}
H(\underline{ststs}\bw_0^n)&=
H(\us\ \underline{tsts}\bw_0^n)-h^2tH(\underline{sts}\bw_0^n)
\text{ (par le lemme \ref{rappelG2}(iv)) }\\
&=h^9t^4H(\us\ \underline{ts}\bw_0^{n-1})-h^9t^4H(\underline{sts}\bw_0^{n-1})\\
&\phantom{=}\text{ (voir les cas de $\underline{tsts}\bw_0^n$ et
de $\underline{sts}\bw_0^n$) }\\
&=h^{11}t^5H(\us\bw_0^{n-1})
\text{ (par le lemme \ref{rappelG2}(ii)) }\\
&=h^7t^3 E(H(\underline{ststs}\bw_0^{n-1}))
\text{ (par la table) }.
\end{split}
\end{equation*}

\casde{\bs\bw_0^n}, $\bs\bs\bw_0^n$, $\bt\bw_0^n$ et $\bt\bt\bw_0^n$.
Si $n=0$ la valeur se d\'eduit par exemple de la proposition \ref{H_s}. Sinon,
nous appliquons la proposition \ref{s^mb} avec $\bb=\bt\bs\bt\bs\bt\bw_0^{n-1}$
pour $m=0,1,2,3$.
Par l'hypoth\`ese de r\'ecurrence nous avons $H(\bb)=(h^7t^3)^{n-1}E^{n-1}
(h^5t^2(\tau+\rho)+h^6t^3(\sigma+\rho))=(h^7t^3)^nE^n
(h^{-2}t\inv(\tau+A)+h\inv(\sigma+A))$, et $H(\bs\bb)=H(\bw_0^n)=
(h^7t^3)^n(\tau+\sigma+2A)$. Ces deux \'egalit\'es permettent dans la proposition 
\ref{s^mb} de d\'eterminer $H_i=(h^7t^3)^nE^n(h^{-2}t\inv(\tau+A))$ et 
$H_s=(h^7t^3)^nE^n(h\inv(\sigma+A))$
et on  en d\'eduit les  valeurs pour $m=2,3$.  On proc\`ede de  fa\c con
analogue pour $\bt\bw_0^n$ et $\bt\bt\bw_0^n$.

Jusqu'\`a la fin de la preuve, tous les \'el\'ements de $\uB^+$ que nous
allons  consid\'erer  sont  de  la  forme  $y\bw_0^n$.  Pour  simplifier
les  notations,  nous posons  $\uH(y)=(h^7t^3)^{-n}E^n(H(y\bw_0^n))$  et
$\uH^i_c(\bX(y))=E^n(H^{i+7n}_c(\bX(y\bw_0^n))(3n))$.

\casde{\us\bt} et $\us\bt\bs$, $\ut\bs$ et $\ut\bs\bt$.
Par $\seof{\us\bt}{\us\ut}\us$ on trouve qu'il existe
$\varepsilon_A,\varepsilon_\tau\in\{0,1\}$ tels que
$\uH(\us\bt)=h(\tau+A)+h^2t(\sigma+J)+(h^2+h^3)t(\varepsilon_A
A+\varepsilon_\tau\tau)$. Par $\seof{\us\bt\bs}{\underline{sts}}{\us\ut}$
on trouve qu'il existe $\varepsilon_\sigma,\varepsilon_j,\varepsilon_{j^2}
\in\{0,1\}$ tels que
$\uH(\us\bt\bs)=h^4t^2(\sigma+J)+h^3t(\tau+A)+(h^2+h^3)t(\varepsilon_\sigma\sigma
+\varepsilon_j \rho_j+\varepsilon_{j^2}\rho_{j^2})$. Mais par \ref{rappel} (i) et (iii) on a
$\uH(\us\bt\bs)=h^2t \uH(\us\bt)$,  donc en comparant les  valeurs on trouve
que  $\varepsilon_A=\varepsilon_\tau=\varepsilon_\sigma=\varepsilon_j=\varepsilon_{j^2}=0$
d'o\`u  le  r\'esultat  pour   les  deux  premi\`eres  vari\'et\'es.  On
proc\`ede de fa\c con analogue pour $\ut\bs$ et $\ut\bs\bt$.

\casde{\bs\ut\bs}.   Par   $\seof{\bs\ut\bs}{\bs\ut\us}{\bs\ut}$,   o\`u
$\uH(\bs\ut\us)=h^2t \uH(\ut\us)$  par \ref{rappel} (i) et  (iii), on trouve
la valeur.

\casde{\us\bt\bs\bt} et de $\us\bt\bs\bt\bs$, 
$\ut\bs\bt\bs$, $\ut\bs\bt\bs\bt$.
Par $\seof{\us\bt\bs\bt\bs}{\underline{ststs}}{\underline{stst}}$
on trouve qu'il existe $\varepsilon_\tau,\varepsilon_\rho\in\{0,1\}$ tels que
$\uH(\us\bt\bs\bt\bs)=h^6t^3(\tau+\rho)+h^5t^2(\sigma+J)+(h^4+h^5)t^2
(\varepsilon_\tau\tau+\varepsilon_\rho\rho)$. Par
$\seof{\us\bt\bs\bt}{\underline{stst}}{\underline{sts}}$ on trouve qu'il existe
$\varepsilon_\sigma,\varepsilon_j,\varepsilon_{j^2}\in\{0,1\}$ tels que
$\uH(\us\bt\bs\bt)=h^3t(\sigma+J)+h^4t^2(\tau+\rho)+(h^4+h^5)t^2
(\varepsilon_\sigma\sigma+\varepsilon_j\rho_j+\varepsilon_{j^2}\rho_{j^2})$. Mais on a
$\uH(\us\bt\bs\bt\bs)=\uH(\bt\bs\bt\bs\us)=h^2t \uH(\bt\bs\bt\us)=h^2t
\uH(\us\bt\bs\bt)$. En comparant les valeurs  on trouve le r\'esultat. Les
vari\'et\'es  $\ut\bs\bt\bs$ et  $\ut\bs\bt\bs\bt$ se  traitent par  des
arguments sym\'etriques.

\casde{\bs\bt}     et      de     $\bs\bt\bs$,      $\bt\bs\bt$.     Par
$\seof{\bs\bt}{\us\bt}\bt$  on   trouve  qu'il   existe  $\varepsilon_A,
\varepsilon_\tau,     \varepsilon_\sigma     \in\{0,1\}$    tels     que
$\uH(\bs\bt)=h^3tA+h^2tJ+  (h^2+h)  (\varepsilon_\tau\tau +  \varepsilon_A
A)+   (h^2+h^3)t   \varepsilon_\sigma   \sigma$.   Par   la   sym\'etrie
entre    $\sigma$     et    $\tau$,    (qu'on    peut     obtenir    ici
en     consid\'erant      $\seof{\bs\bt}{\bs\ut}\bs$),     on     trouve
$\varepsilon_\sigma=\varepsilon_\tau=0$.

  Par $\seof{\bs\bt\bs}{\us\bt\bs}{\bt\bs}$ on trouve qu'il existe
$\varepsilon\in\{0,1\}$ tel que $\uH(\bs\bt\bs)=h^3t(\tau+J)+h^4t^2(\sigma+J)
+(h^2+h^3)\varepsilon_A A+(h^3+h^4)t\varepsilon A$. En comparant
avec $\seof{\bs\bt\bs}{\bs\ut\bs}{\bs\bs}$ on trouve que $\varepsilon=0$.
Enfin, en comparant avec la suite exacte longue qui vient de la proposition \ref{sws}
\begin{align*}\cdots&\to \uH^0_c(\bX(\bt))(-1)\to
\uH^1_c(\bX(\bs\bt))(-1)\oplus
\uH^2_c(\bX(\bs\bt))\to \uH^3_c(\bX(\bs\bs\bt))\\
&\to\uH^1_c(\bX(\bt))(-1)\to\cdots\\
\end{align*}
on trouve que $\varepsilon_A=0$.
On obtient $\uH(\bt\bs\bt)$ par un argument sym\'etrique.

\casde{\bs\bt\bs\bt}. Par $\seof{\bs\bt\bs\bt}{\us\bt\bs\bt}{\bt\bs\bt}$
on trouve qu'il existe $\varepsilon_\tau,\varepsilon_\sigma,
\varepsilon_j,\varepsilon_{j^2}\in\{0,1\}$ tels que $\uH(\bs\bt\bs\bt)=
h^4t^2\rho+h^5t^2J+(h^3+h^4)t(\varepsilon_\sigma\sigma+\varepsilon_j\rho_j+\varepsilon_{j^2}\rho_{j^2})+
(h^4+h^5)t^2\varepsilon_\tau\tau$. La sym\'etrie entre $\sigma$ et $\tau$
(qu'on peut obtenir ici en consid\'erant
$\seof{\bs\bt\bs\bt}{\bs\bt\bs\ut}{\bs\bt\bs}$) montre que
$\varepsilon_\tau=\varepsilon_\sigma=0$.

\casde{\bt\bs\bt\bs\bt}. Par
$\seof{\bt\bs\bt\bs\bt}{\ut\bs\bt\bs\bt}{\bs\bt\bs\bt}$
on trouve qu'il existe $\alpha_j,\alpha_{j^2}\in\{0,1\}$ tels que
$\uH(\bt\bs\bt\bs\bt)=h^6t^3(\sigma+\rho)+(h^5t^2)(\tau+\rho)+
(h^5+h^6)t^2(\alpha_j\rho_j+\alpha_{j^2}\rho_{j^2})+(h^4+h^5)t(\varepsilon_j\rho_j+\varepsilon_{j^2}\rho_{j^2})$.

\casde{\bw_0}. Par $\seof{\bw_0}{\us\bt\bs\bt\bs\bt}{\bt\bs\bt\bs\bt}$ o\`u
$\uH(\us\bt\bs\bt\bs\bt)=h\uH(\underline{ststs})$ (\cf\ lemme \ref{rappelG2}(i))
on trouve qu'il
existe $\varepsilon_\tau,\varepsilon_\rho\in\{0,1\}$ tels que $\uH(\bw_0)=
h^7t^3(\sigma+\tau+2\rho)+(h^5+h^6)t(\varepsilon_j\rho_j+\varepsilon_{j^2}\rho_{j^2})
+(h^6+h^7)t^2(\alpha_j\rho_j+\alpha_{j^2}\rho_{j^2})
+(h^5+h^6)t^2(\varepsilon_\tau\tau+\varepsilon_\rho\rho)$.

Nous    appliquons   maintenant    la   proposition    \ref{s^mb}   avec
$\bb=\bt\bs\bt\bs\bt$,  sous   la  forme   (en  suivant   les  notations
de  la  proposition  \ref{s^mb}):  $\uH(\bs\bb)-h  \uH(\bb)=(h^2t-h)H_i$,  en
particulier  on  voit  que  $\uH(\bs\bb)-h  \uH(\bb)$  doit  \^etre  divisible
par   $h^2t-h$.   Avec   les   valeurs  obtenues,   cela   implique   que
$\varepsilon_\tau=\varepsilon_\rho=0$.

Nous finissons en utilisant le corollaire
\ref{valeurs  propres dans H(w0)} qui montre
que les caract\`eres $\rho_j$ et $\rho_{j^2}$ ne peuvent intervenir dans
la cohomologie de $X(\bpi^n\bw_0)$, donc que $\varepsilon_j=\varepsilon_{j^2}=\alpha_j
=\alpha_{j^2}=0$
(notons que le corollaire \ref{valeurs  propres dans H(w0)} d\'epend du
th\'eor\`eme \ref{2G2} mais pas du th\'eor\`eme \ref{G2}).
\end{pf*}
\begin{theoreme}\label{2G2}
Supposons  $(\bG,F)$  de  type  $\lexp 2G_2$  et  soit  $y\in\uB^+$
apparaissant  dans  la table  ci-dessous.  Alors  pour  tout $n$,  on  a
$H(y\bw_0^n)=(h^7 t^3)^n E^n(f(y))$ o\`u $f(y)$  est la valeur de $H(y)$
donn\'ee par la table; dans cette table nous avons pos\'e
$A=\rho_i+\rho_{-i}+\rho_{\zeta_{12}^5}+\rho_{\zeta_{12}^7}$ et
$B=\rho'_i+\rho'_{-i}+\rho_{\zeta_{12}^5}+\rho_{\zeta_{12}^7}$; avec ces
notations $E$ \'echange $A$ et $B$.
\halign{$#$\hfill&\quad$#$\hfill\cr
\strut y&H(y)\cr
\noalign{\hrule}\cr
1&0\cr
\bs,\bt,\us,\ut&ht^{1/2}A\cr
\us\ut,\ut\us&(ht^{1/2}+h^3t^{3/2})A\cr
\bs\bs,\bs\bt,\bt\bs,\bt\bt&(h^2t^{1/2}+h^3t^{3/2})A\cr
\underline{sts},\underline{tst}&h^3t^{3/2}(A+B)\cr
\us\bt\bs&h^2t^{1/2}A+h^3t^{3/2}B\cr
\bs\bt\bs,\bt\bs\bt&h^4t^{3/2}A+h^3t^{3/2}B+\varepsilon_n(h^3+h^4)t^{3/2}
(\rho_{\zeta_{12}^5}+\rho_{\zeta_{12}^7})\cr
\underline{stst},\underline{tsts}&(h^3t^{3/2}+h^5t^{5/2})B\cr
\us\bt\bs\bt,\ut\bs\bt\bs& h^4t^{3/2}A+h^5t^{5/2}B\cr
\bs\bt\bs\bt,\bt\bs\bt\bs&(h^4t^{3/2}+h^5t^{5/2})B\cr
\underline{ststs}& h^5t^{5/2}B\cr
\us\bt\bs\bt\bs& h^4t^{3/2}B\cr
\bs\bt\bs\bt\bs,\bt\bs\bt\bs\bt,\us\bt\bs\bt\bs\bt& h^6t^{5/2}B\cr
}
Ci-dessus $\varepsilon_n\in\{0,-1\}$ est une constante ind\'ependante
de $\ell$.
\end{theoreme}
\begin{pf*}{Preuve}
La preuve proc\`ede comme pour le th\'eor\`eme \ref{G2}, sauf que nous allons
utiliser l'ingr\'edient suppl\'ementaire (les notations sont celles de
\S \ref{indep}):

\begin{proposition}\label{indepG2}
Nous fixons $i,j\in\bbN, \bb\in\uB^+$. Soit
$M_\ell=H^i_c(\bX(\bb),\bbQ_\ell)_j$; supposons que
$\scal {M_\ell}{M_\ell^*}$ est ind\'ependant de $\ell$, et que
$\scal{\rho'_i}{M_\ell}=\scal{\rho'_{-i}}{M_\ell}=0$. Alors
$\scal{\rho_i}{M_\ell}=\scal{\rho_{-i}}{M_\ell}$, $\scal
{\rho_{\zeta_{12}^5}}{M_\ell}=\scal{\rho_{\zeta_{12}^7}}{M_\ell}$ et
$M_\ell$ est ind\'ependant de $\ell$.
\end{proposition}
\begin{pf*}{Preuve}
Soient      $a_i,a_{{-i}},a_{{\zeta_{12}^5}},a_{{\zeta_{12}^7}}$     les
multiplicit\'es                      respectives                      de
$\rho_i,\rho_{-i},\rho_{\zeta_{12}^5},\rho_{\zeta_{12}^7}$          dans
$M_\ell$.                                                          Comme
$\rho_i+\rho_{-i}+\rho_{\zeta_{12}^5}+\rho_{\zeta_{12}^7}$            et
$\rho'_i+\rho'_{-i}+\rho_{\zeta_{12}^5}+\rho_{\zeta_{12}^7}$        sont
combinaisons  lin\'eaires de $R_{\tilde\chi}$,  on d\'eduit de
la proposition \ref{Rchi indepl}    et   de    l'hypoth\`ese   que
$\alpha:=a_i+a_{{-i}}$   et
$\beta:=a_{{\zeta_{12}^5}}+a_{{\zeta_{12}^7}}$  sont  ind\'ependants  de
$\ell$.          L'hypoth\`ese          nous          donne          que
$\gamma:=a_ia_{{-i}}+a_{{\zeta_{12}^5}}
a_{{\zeta_{12}^7}}=a_i(\alpha-a_i)+
a_{{\zeta_{12}^5}}(\beta-a_{{\zeta_{12}^5}})$   est   ind\'ependant   de
$\ell$.  Il  existe  une  infinit\'e  de  $\ell$ tels que ni $\rho_i$ ni
$\rho_{\zeta_{12}^5}$  ne soient  \`a valeurs  dans $\bbQ_\ell$
(th\'eor\`eme \ref{Chebotarev},  ce qui
pour  ces $\ell$ force $a_i=\alpha/2$ et $a_{{\zeta_{12}^5}}=\beta/2$ ce
qui  donne  $\gamma=(\alpha^2+\beta^2)/4$.  Donc  pour  tout $\ell$ on a
$\gamma=(\alpha^2+\beta^2)/4$.  Mais  on  a toujours $a_i(\alpha-a_i)\le
\alpha^2/4$   avec   \'egalit\'e   uniquement   si   $a_i=\alpha/2$  (et
similairement   pour  $a_{{\zeta_{12}^5}}$   et  $\beta$).   On  a  donc
n\'ecessairement                $a_i=\alpha/2=a_{{-i}}$               et
$a_{{\zeta_{12}^5}}=\beta/2=a_{{\zeta_{12}^7}}$. \end{pf*}

Nous avons \'evidemment la proposition analogue en \'echangeant $\rho_i$
et $\rho'_i$ et $\rho_{-i}$ et $\rho'_{-i}$. En g\'en\'eral nous saurons
que    l'hypoth\`ese   de    la   proposition    a   lieu   car   $\scal
{M_\ell}{M^*_\ell}$   sera   le   seul   terme   inconnu   de  la  somme
$S_{2i,2j,\bb,\bb}$.

Nous commen\c cons par
\casde{\us\bw_0^n}, $\ut\bw_0^n$, $\us\ut\bw_0^n$, $\ut\us\bw_0^n$,
$\underline{sts}\bw_0^n$, $\underline{tsts}\bw_0^n$,
$\underline{stst}\bw_0^n$,
$\underline{ststs}\bw_0^n$.  La preuve  proc\`ede exactement  comme pour
les m\^emes \'el\'ements dans le cas de $G_2$.

Jusqu'\`a la fin de la preuve, tous les \'el\'ements de $\uB^+$ que nous
allons  consid\'erer sont  de la  forme $y\bw_0^n$.  Pour simplifier les
notations,  nous  posons  $\uH(y)=(h^7t^3)^{-n}E^n(H(y\bw_0^n))$
et nous noterons $\uH^i(y)_j$ la partie
$\rho_i,\rho_{-i},\rho'_i,\rho'_{-i},\rho_{\zeta_{12}^5},
\rho_{\zeta_{12}^7}$-isotypique de
$E^n(H^{i+7n}_c(\bX(y\bw_0^n),\bbQ_\ell)_{j+6n})$.  
Comme
expliqu\'e  dans  la  preuve  du  th\'eor\`eme  \ref{A2},  les  items de
\ref{rappel}  sont encore  vrais pour  $\uH$ sauf  le (i) qui est encore
vrai si $x\in B^+$.

\casde{\bs}   et  de   $\bt$.  L'\'el\'ement   $\bs$  est   donn\'e  par
$\seof\bs\us1$. On en d\'eduit $\bt$ par \ref{rappel}(i).

\casde{\bs\bs},  $\bt\bs$, $\bs\bt$ et  $\bt\bt$. L'\'el\'ement $\bs\bs$
est donn\'e par $\seof{\bs\bs}{\bs\us}\bs$ o\`u $\bs\us$ est donn\'e par
\ref{rappel}(iii).    On   d\'eduit   les    autres   \'el\'ements   par
\ref{rappel}(i).

\casde{\us\bt\bs},  de $\us\bt\bs\bt$  et de  $\ut\bs\bt\bs$. Remarquons
tout  d'abord qu'on  a $\uH(\us\bt\bs\bt)=  \uH(\bs\us\bt\bs)\text{ (par
\ref{rappel}(i)) }=h^2t\uH(\us\bt\bs)\text{ (par \ref{rappel}(iii))}$.

Par $\seof{\us\bt\bs}{\underline{sts}}{\underline{ts}}$ on trouve que
$\uH(\us\bt\bs)=h^2t^{1/2}A+h^3t^{3/2}B+(h^3+h^4)t^{3/2}{}_\le A$, o\`u
${}_\le A$ d\'enote un sous-module de $A$. Alors par
$\seof{\us\bt\bs\bt}{\underline{stst}}{\underline{tst}}$ et la remarque
pr\'eliminaire on trouve que ${}_\le A=0$.

On obtient $\ut\bs\bt\bs$ par un raisonnement sym\'etrique.

\casde{\us\bt\bs\bt\bs\bt} et de $\us\bt\bs\bt\bs$. Par un raisonnement
analogue au cas pr\'ec\'edent, on a $\uH(\us\bt\bs\bt\bs\bt)=
h^2t\uH(\us\bt\bs\bt\bs)$, et on obtient $\uH(\us\bt\bs\bt\bs\bt)$
par le lemme \ref{rappelG2}(i), d'o\`u les deux valeurs.

\casde{\bs\bt\bs}, $\bt\bs\bt\bs$, $\bs\bt\bs\bt\bs$ et $\bw_0$,
ainsi que $\bt\bs\bt$, $\bs\bt\bs\bt$ et $\bt\bs\bt\bs\bt$.
Par $\seof{\bs\bt\bs}{\us\bt\bs}{\bt\bs}$ on trouve qu'il existe
$\beta,\beta^*\in\{0,-1\}$ et un sous-module $A_1$ de $A$ tels que
$\uH(\bs\bt\bs)=h^3t^{3/2}B+h^4t^{3/2}A+(h^3+h^4)t^{3/2}
(\beta \rho_{\zeta_{12}^5}+\beta^*\rho_{\zeta_{12}^7})+
(h^2+h^3)t^{1/2}A_1$. On peut appliquer \ref{indepG2} avec
$M_\ell=H^2(\bs\bt\bs)_1$ car $\scal{M_\ell}{M_\ell^*}$ est le
seul terme non nul de $S_{4+14n,2+12n,\bs\bt\bs,\bs\bt\bs}$, ce qui prouve que 
$A_1$
est rationnel et ind\'ependant de $\ell$. On voit de m\^eme que $\beta$ est
ind\'ependant de $\ell$ et que $\beta=\beta^*$ en appliquant \ref{indepG2}
avec $M_\ell=H^4(\bs\bt\bs)_3$.

Ensuite,  par  $\seof{\bt\bs\bt\bs}{\ut\bs\bt\bs}{\bs\bt\bs}$  on trouve
qu'il    existe    un    sous-module    $A_2$    de    $A$    tel    que
$\uH(\bt\bs\bt\bs)=(h^4t^{3/2}+h^5t^{5/2})B+(h^3+h^4)t^{1/2}A_1+
(h^4+h^5)t^{3/2}A_2$. En appliquant \ref{indepG2} avec
$M_\ell=H^5(\bt\bs\bt\bs)_3$ on voit que $A_2$ est ind\'ependant de $\ell$.

Puis, par $\seof{\bs\bt\bs\bt\bs}{\us\bt\bs\bt\bs}{\bt\bs\bt\bs}$ on trouve
qu'il existe un sous-module $B_1$ de $B$ tel que
$\uH(\bs\bt\bs\bt\bs)=h^6t^{5/2}B+
(h^4+h^5)t^{1/2}A_1+ (h^4+h^5)t^{3/2}B_1+ (h^5+h^6)t^{3/2}A_2$.
En appliquant \ref{indepG2} avec
$M_\ell=H^4(\bs\bt\bs\bt\bs)_3$ on voit que $B_1$ est ind\'ependant de $\ell$.

Enfin,   par   $\seof{\bw_0}{\us\bt\bs\bt\bs\bt}{\bs\bt\bs\bt\bs}$  o\`u
$\us\bt\bs\bt\bs\bt$  est  donn\'e  par  le  lemme \ref{rappelG2}(i), on
trouve qu'il existe un sous-module $B_2$ de $B$ tel que
$\uH(\bw_0)= (h^5+h^6)t^{1/2}A_1+ (h^5+h^6)t^{3/2}B_1+
(h^6+h^7)t^{3/2}A_2+ (h^6+h^7)t^{5/2}B_2$.
En appliquant \ref{indepG2} avec
$M_\ell=H^7(\bw_0)_5$ on voit que $B_2$ est ind\'ependant de $\ell$.

En  appliquant alors le lemme  \ref{rappelG2}(vii)  on  trouve
que $A_1=A_2=B_1=B_2=0$ (des inconnues  que nous avons introduites,
seule  $\beta$  n'a  pas  \'et\'e   prouv\'ee  nulle;  nous  l'avons  not\'ee
$\varepsilon_n$ dans la table).

Les \'el\'ements $\bt\bs\bt$, $\bt\bs\bt\bs\bt$ donnent la m\^eme valeur
que respectivement $\bs\bt\bs$ et $\bs\bt\bs\bt\bs$ par \ref{rappel}(i),
et $\bs\bt\bs\bt$  s'obtient par  un raisonnement sym\'etrique  de celui
utilis\'e pour $\bt\bs\bt\bs$.
\end{pf*}

\section{Endomorphismes des vari\'et\'es $\bX(\bw)$}
\label{sectionendo}

Pour  $\bt\in\BW$, suivant  \cite{Sydney},  nous nous  int\'eresserons
dans cette  partie \`a  certains morphismes  commutant \`a  l'action de
$\GF$ entre vari\'et\'es $\bX(\bt)$ qui nous permettront de construire
l'action de $C_\bB(\bw F)$ mentionn\'ee dans l'introduction. 

\subsection{Combinatoire}
Nous consid\'erons la cat\'egorie $\CB$ d'objets les \'el\'ements de $B$,
avec
$$\Hom_{\CB}(\bw,\bw')=\{\by\in B|\bw'=\by\inv\bw F(\by)\}$$
 et la composition est
donn\'ee par le produit. Soit $\CB^+$ la sous-cat\'egorie pleine d'objets
les \'el\'ements de $B^+$.

Nous introduisons deux sous-cat\'egories de $\CB^+$, dont l'ensemble
des objets co{\"\i}ncide avec $B^+$.

La cat\'egorie $\cD^+$ est la plus petite sous-cat\'egorie telle que
$$\{\by\in B^+\mid\by\preccurlyeq\bw, \by\inv\bw F(\by)=\bw'\}
\subset\Hom_{\cD^+}(\bw,\bw')$$
et $\cD$ est la plus petite sous-cat\'egorie contenant
$\cD^+$ et dont toutes les fl\`eches sont inversibles.

\begin{remarque}\label{B+w->w'}
\noindent
\begin{enumerate}
\item  $\Hom_{\cD^+}(\bw,\bw')$ est l'ensemble  des $\by\in
B^+$   qui  admettent   une  d\'ecompo\-sition   $\by=\by_1\cdots\by_n$ 
avec $\by_i\in B^+$, telle
qu'il  existe  une   suite  $\bw_1,\ldots,\bw_{n+1}$  d'\'el\'ements  de
$B^+$   avec   $\bw_1=\bw$,   $\bw_{n+1}=\bw'$,   $\by_i\preccurlyeq\bw_i$   et
$\bw_{i+1}=\by_i\inv\bw_i\ F(\by_i)$.
\item  De   m\^eme,  $\Hom_\cD(\bw,\bw')$  est l'ensemble
des  $\by\in   B$  qui  admettent   une  d\'ecomposition 
$\by=\by_1\cdots\by_n$  avec $\by_i\in B^+$ ou $\by\inv_i\in B^+$,
telle qu'il  existe une  suite $\bw_1,\ldots,\bw_{n+1}$  d'\'el\'ements
de  $B^+$ avec  $\bw_1=\bw$, $\bw_{n+1}=\bw'$,  $\by_i\preccurlyeq\bw_i$ si
$\by_i\in B^+$ et
$\by_i\inv\preccurlyeq\bw_{i+1}$ si $\by_i\inv\in B^+$,
 et $\bw_{i+1}=\by_i\inv\bw_i F(\by_i)$.
\end{enumerate}
\end{remarque}

Soit $\tilde{\cD}^+$ la cat\'egorie d'ensemble d'objets $B^+$ et dont
les fl\`eches sont engendr\'ees par les
$\gamma_\by=\gamma_\by^{\bw,\bw'}:\bw\to \bw'$ pour $\by\in B^+$ v\'erifiant
$\by\preccurlyeq\bw$ et $\by\inv\bw F(\by)=\bw'$, avec les relations
$\gamma^{\bw,\bw''}_\by=\gamma^{\bw',\bw''}_{\by_2}\gamma^{\bw,\bw'}_{\by_1}$ lorsque
$\by=\by_1\by_2\preccurlyeq\bw$.

On a un foncteur canonique $\tilde{\cD}^+\to\CB$ qui est l'identit\'e sur
les objets et envoie $\gamma_\by^{\bw,\bw'}$ sur $\by$.

\begin{proposition}
\label{fidele+}
Le foncteur canonique $\tilde{\cD}^+\to\CB$ induit un isomorphisme
$\tilde{\cD}^+\iso \cD^+$.
\end{proposition}

\begin{pf*}{Preuve}
Il s'agit de montrer que $\tilde{\cD}^+\to\CB$ est pleinement fid\`ele.

Soit $\CE(\by,\bw,\bw')$ l'ensemble  des suites  $(\by_1,\ldots,\by_n)$
d'\'el\'ements de $B^+$ de produit $\by$ et telles qu'il
existe  une   suite  $\bw_1,\ldots,\bw_{n+1}$  d'\'el\'ements  de
$B^+$   avec   $\bw_1=\bw$,   $\bw_{n+1}=\bw'$,   $\by_i\preccurlyeq\bw_i$   et
$\bw_{i+1}=\by_i\inv\bw_i\ F(\by_i)$.
L'image du foncteur s'identifie \`a la sous-cat\'egorie de $\CB$ d'objets
$B^+$ et d'ensemble de fl\`eches de $\bw$ \`a $\bw'$ les
$\by\in\Hom_\CB(\bw,\bw')$ tels que $\CE(\by,\bw,\bw')\not=\emptyset$.

Nous allons d\'emontrer par r\'ecurrence sur $l(\by)$ que:
\begin{enumerate}
\item Si $\bv\in\bW$,  $\bv\preccurlyeq\bw$  et  $\bv\preccurlyeq\by$,
alors   il    existe   $(\bu_1,\ldots,\bu_m)\in\CE(\by,\bw,\bw')$   avec
$\bu_1=\bv$.
\item Si  $(\by_1, \ldots,\by_n),  (\bz_1,\ldots,\bz_k)
\in\CE(\by,\bw,\bw')$, alors  on   a  $\gamma_{\by_n}\cdots
\gamma_{\by_1}= \gamma_{\bz_k} \cdots \gamma_{\bz_1}$.
\end{enumerate}
La proposition r\'esultera de l'assertion (ii).

Consid\'erons  $\bv$  comme en  (i).
Posons  $\bv=\bs\bv'$ avec  $\bs\in\bS$. Soit
$(\by_1,\ldots,\by_n)\in \CE(\by,\bw,\bw')$.

Supposons tout d'abord  $\bs\preccurlyeq\by_1$.
Posons $\bw_1=\bs\inv\bw F(\bs)$   et   $\by_1=\bs\by'_1$.
Alors,   $\bv'\preccurlyeq\bw_1$ et    $\bv'\preccurlyeq\by'_1   \by_2\cdots\by_n$,   donc
par   r\'ecurrence  il   existe 
$(\bv',\bu_2,  \ldots,\bu_m)\in\CE(\by'_1\by_2
\ldots\by_n,\bw_1,\bw')$
et  $(\bv,\bu_2\ldots,\bu_m)$  satisfait l'assertion (i).
Par r\'ecurrence, on a
$\gamma_{\bu_m}\cdots \gamma_{\bu_2} \gamma_{\bv'}=
\gamma_{\by_n}\cdots \gamma_{\by_2} \gamma_{\by'_1}$, donc
$\gamma_{\bu_m}\cdots \gamma_{\bu_2} \gamma_{\bv}=
\gamma_{\by_n}\cdots \gamma_{\by_1}$.

Supposons maintenant $\bs\not\preccurlyeq\by_1$.
Soit $\bt\in\bS$ tel que $\by_1=\bt\by'_1$
Soit $\bd$ le p.p.c.m. de $\bs$ et  $\bt$, \ie, le plus court \'el\'ement
de $B^+$ divisible \`a gauche par $\bs$ et $\bt$.
On a 
$\bs\preccurlyeq\bw$,  $\bs\preccurlyeq\by$,   $\bt\preccurlyeq\bw$ et $\bt\preccurlyeq\by$, donc
$\bd\preccurlyeq\bw$ et $\bd\preccurlyeq\by$ (lemme \ref{alpha}).
Soit $\bd'\in B^+$ tel que $\bd=\bt\bd'$ et soit $\bw_1=\bt\inv\bw F(\bt)$.
Alors, $\bd'\preccurlyeq\bw_1$ et $\bd'\preccurlyeq\by'_1\by_2\cdots\by_n$.
Par r\'ecurrence, on d\'eduit de (i)
qu'il existe $(\bd',\by'_2,\ldots,\by'_k)\in 
\CE(\by'_1\by_2\cdots\by_n,\bw_1,\bw')$, donc
$(\bd,\by'_2,\ldots,\by'_k)\in\CE(\by,\bw,\bw')$.
Puisque $\bs\preccurlyeq\bd$, on d\'eduit de l'\'etude pr\'ec\'edente
qu'il existe
$(\bv,\bu_2,  \ldots,\bu_m)\in\CE(\by,\bw,\bw')$, d'o\`u (i).
L'\'etude pr\'ec\'edente montre aussi que
$\gamma_{\by'_k}\cdots \gamma_{\by'_2} \gamma_{\bd}=
\gamma_{\by_n}\cdots \gamma_{\by_1}$, donc,
en utilisant \`a nouveau le cas pr\'ec\'edent,
$\gamma_{\bu_m}\cdots \gamma_{\bu_2} \gamma_{\bv}=
\gamma_{\by_n}\cdots\gamma_{\by_1}$.

Soient    maintenant   $(\bz_1,    \ldots,\bz_k) \in\CE(\by,\bw,\bw')$.
Soit  $\bu\in\bW$ tel  que $\bz_1=\bu\bz'_1$.
Nous avons montr\'e qu'il existe
$(\bu,\bu_2,\ldots,\bu_m)\in\CE(\by,\bw,\bw')$ tel que
$\gamma_{\bu_m}\cdots \gamma_{\bu_2} \gamma_{\bu}=
\gamma_{\by_n}\cdots \gamma_{\by_1}$.
Par  r\'ecurrence, (ii) montre que 
$\gamma_{\bu_m}\cdots \gamma_{\bu_2}=
\gamma_{\bz_k}\cdots \gamma_{\bz_2}\gamma_{\bz'_1}$, donc
$\gamma_{\bu_m}\cdots \gamma_{\bu_2}\gamma_{\bu}=
\gamma_{\bz_k}\cdots \gamma_{\bz_2}\gamma_{\bz_1}$,
d'o\`u le (ii).
\end{pf*}

\begin{remarque}
On en d\'eduit que le sous-mono{\"\i}de $B^+_\bw\subset   C_{B^+}(\bw   F)$
consi\-d\'er\'e dans \cite[5.6(ii)]{Sydney} 
s'identifie   \`a   $\End_{\cD^+}(\bw)$.
\end{remarque}

Soit $\tilde{\cD}$ la cat\'egorie d\'eduite de $\tilde{\cD}^+$ en rendant toutes
les fl\`eches inversibles.
Puisque $\CB$ est un groupo{\"\i}de, 
le foncteur canonique $\tilde{\cD}^+\to\CB$ s'\'etend de mani\`ere unique
en un foncteur $\tilde{\cD}\to\CB$.

\begin{proposition}
\label{fidele}
Le foncteur canonique $\tilde{\cD}\to\CB$ induit un isomorphisme
$\tilde{\cD}\iso \cD$.
\end{proposition}

Commen\c{c}ons par un lemme.

\begin{lemme}\label{prolongement}
Soit $X$ un ensemble fini de fl\`eches de $\tilde\cD^+$ de m\^eme source
v\'erifiant:
\begin{enumerate}
\item  Tout facteur initial d'une fl\`eche de  $X$ est dans $X$, \ie\ si
$\gamma'\gamma''\in X$ alors $\gamma''\in X$.
\item Si, pour $\bs,\bt\in\bS$ on a $\gamma_\bs\gamma\in X$ et
$\gamma_\bt\gamma\in X$, alors $\gamma_{\bw_0^{\{\bs,\bt\}}}\gamma\in X$.
\end{enumerate}
Alors $X$ est l'ensemble des facteurs initiaux d'une fl\`eche
$\gamma\in\tilde\cD^+$.
\end{lemme}
\begin{pf*}{Preuve}
Remarquons   d'abord  que   la  condition   (ii)  a   un  sens: si
$\gamma_\bs^{\bw,\bs\inv\bw   F(\bs)}$   et  $\gamma_\bt^{\bw,\bt\inv\bw
F(\bt)}$  sont  deux  fl\`eches  de  m\^eme  source,  alors  la fl\`eche
$\gamma_{\bw_0^{\{\bs,\bt\}}}^{\bw,{\bw_0^{\{\bs,\bt\}}}\inv\bw
F(\bw_0^{\{\bs,\bt\}})}$    existe,    car    $\bs\preccurlyeq\bw$    et
$\bt\preccurlyeq\bw$   implique  ${\bw_0^{\{\bs,\bt\}}}\preccurlyeq\bw$.
Soit  $\gamma$ une fl\`eche  de longueur maximale  dans $X$. Si $\gamma$
n'est pas multiple de toutes les fl\`eches de $X$, alors il existe $f\in
X$  qui  est  un  facteur  initial  de  $\gamma$ et $\bs\in\bS$ tels que
$\gamma_\bs  f\in X$ et  $\gamma_\bs f$ n'est  pas un facteur initial de
$\gamma$.  Fixons un tel $f$ de  longueur maximale. Comme la longueur de
$f$  est inf\'erieure \`a  celle de $\gamma$,  il existe $\bt\in\bS$ tel
que $\gamma_\bt f$ soit un facteur de $\gamma$. Par (i) on a $\gamma_\bt
f\in  X$ et (ii) implique que $\gamma_1=\gamma_{\bw_0^{\{\bs,\bt\}}}f\in
X$.  Comme  $\gamma_\bs  f$  est  un  facteur  initial de $\gamma_1$, la
fl\`eche  $\gamma_1$ n'est pas  un facteur initial  de $\gamma$. Donc il
existe  un facteur initial $f'$ de  $\gamma_1$, dont $\gamma_\bt f$ soit
un  facteur et qui soit facteur  de $\gamma$, et un $\bs'\in\{\bs,\bt\}$
tel  que $\gamma_{\bs'}f'$ ne  soit pas facteur  de $\gamma$. Comme $f'$
est plus long que $f$, on a une contradiction.
\end{pf*}

\begin{corollaire}\label{ppcm}
L'ensemble   des  fl\`eches  $\gamma_\bx$  de  $\tilde\cD^+$  de  source
donn\'ee  telles que $\bx\in\bW$ est form\'e des facteurs initiaux d'une
unique fl\`eche $\gamma_\bu$ avec $\bu\in\bW$.
\end{corollaire}
\begin{pf*}{Preuve}
L'ensemble   $X$  des  fl\`eches de l'\'enonc\'e v\'erifie  les  conditions  du
lemme \ref{prolongement}   car   $\bx\bs\in\bW$   et  $\bx\bt\in\bW$  implique
$\bx\bw_0^{\{\bs,\bt\}}\in\bW$.
\end{pf*}

\begin{pf*}{Preuve de la Proposition \ref{fidele}}
Il s'agit de montrer que $\tilde{\cD}\to\CB$ est pleinement fid\`ele.
Le lemme \ref{x-1y=x'y'-1} ci-dessous ram\`ene la preuve de la proposition
au   cas  de  fl\`eches  positives   (\ie,  de  $\tilde{\cD}^+$)  et  la
proposition \ref{fidele+} l'a d\'ej\`a r\'esolu.
\end{pf*}

\begin{lemme}\label{x-1y=x'y'-1}
Soit $\gamma$ une fl\`eche de $\tilde{\cD}$ et $\bx,\by\in B^+$ tels que
$\gamma=\gamma_\bx\gamma_\by\inv$.  Alors,  il  existe $\bx',\by'\in\BW$
tels   que   $\gamma=\gamma_{\bx'}\inv   \gamma_{\by'}$  et  $\by\inv\bx
=\by'\bx^{\prime-1}$.
\end{lemme}
\begin{pf*}{Preuve}
Soient     $\bx=\bx_1\ldots\bx_m$    et    $\by=\by_1\ldots\by_n$    des
d\'ecompositions  en \'el\'ements  de $\bW$.  Nous d\'emontrons le lemme
par  r\'ecurrence  sur  $m+n$;  nous  rajoutons  dans  l'hypoth\`ese  de
r\'ecurrence   que  $\bx'$   et  $\by'$   ont  des  d\'ecompositions  en
\'el\'ements  de $\bW$ de longueurs respectives  $n'$ et $m'$ telles que
$n'+m'\le m+n$.

Par \ref{ppcm}, $\gamma_{\bx_1}$ et $\gamma_{\by_1}$ sont facteurs d'une
fl\`eche  $\gamma_\bu$. \'Ecrivons  $\bu=\bx_1\bx'_1=\by_1\by'_1$. Alors
$\gamma_{\bx_1}\gamma_{\by_1}\inv=\gamma_{\bx'_1}\inv\gamma_{\by'_1}$.

Notons  que comme la propri\'et\'e est claire si $m=0$ ou $n=0$, \ie\ si
$\bx$  ou $\by$ est l'identit\'e, on l'a d\'emontr\'ee si $m+n\le 2$; on
peut aussi supposer $n,m\ge 1$.

On a
$\gamma_\bx\gamma_\by\inv=\gamma_{\bx_2\ldots\bx_m}\gamma_{\bx_1}
\gamma_{\by_1}\inv\gamma_{\by_2\ldots\by_n}\inv=
\gamma_{\bx_2\ldots\bx_m}\gamma_{\bx'_1}\inv             \gamma_{\by'_1}
\gamma_{\by_2\ldots\by_n}\inv$.   Comme  $m<n+m$,  par  r\'ecurrence  il
existe   $\bx'',\by''\in   B^+$   tels   que  $\gamma_{\bx_2\ldots\bx_m}
\gamma_{\bx'_1}\inv=   \gamma_{\bx''}\inv\gamma_{\by''}$.   Notons   que
$\by''$  est produit d'au plus $m-1$  termes dans $\bW$; si $\bx''\ne 1$
ceci   est   assur\'e   par   l'hypoth\`ese   de   r\'ecurrence.   Sinon
$\by''=\bx_1^{\prime-1}\bx_2\ldots\bx_m$ est dans $B^+$, et donc produit
de  $m-1$ termes de $\bW$ au  plus. On peut donc appliquer l'hypoth\`ese
de  r\'ecurrence qui affirme qu'il  existe $\bx''',\by'\in  B^+$ tels que
$\gamma_{\by'_1\by''}                     \gamma_{\by_2\ldots\by_n}\inv=
\gamma_{\bx'''}\inv\gamma_{\by'}$.  On  en  d\'eduit  le lemme en posant
$\bx'=\bx''\bx'''$.
\end{pf*}

\subsection{Traces}
\sub{}
\label{defDb}
Soit $\CC^+$ la cat\'egorie des vari\'et\'es quasi-projectives sur $\Fqbar$
munies des morphismes propres.
On a un foncteur $\tilde{\cD}^+\to\CC^+$ qui envoie
$\bw$ sur $\bX(\bw)$ et $\by\in\Hom_{\tilde{\cD}^+}(\bw,\bw')$ sur $D_\by$ (\cf\ proposition \ref{xy=yFx}).

Soit $\CC$ la localisation de $\CC^+$ en les morphismes qui induisent
des \'equivalences de sites \'etales. Le foncteur pr\'ec\'edent s'\'etend en
un foncteur $\tilde{\cD}\to \CC$.

En composant avec le foncteur $H_c^*$, on obtient un foncteur
de $\CC$ vers la cat\'egorie des $\Qlbar\GF$-modules gradu\'es.

Par   abus    de   notation,    nous   noterons   encore    $$D_\bb   \in
\Hom_{\Qlbar\GF}(H^*_c(\bX(\bw)),  H^*_c(\bX(\bw')))$$ l'image par
le foncteur ci-dessus de $D_\bb\in\Hom_\CC(\bX(\bw),\bX(\bw'))$.
Nous nous int\'eresserons \`a calculer la
repr\'esentation de $\End_\cD(\bw)$ ainsi obtenue dans
$\End_\GF(H^*_c(\bX(\bw)))$, qui sera dans certains cas
une \og alg\`ebre de Hecke cyclotomique\fg\ associ\'ee au groupe $\End_\cD(\bw)$.

Soit $\Sh$ l'op\'erateur de torsion sur les classes de conjugaison
de $\GF$ d\'efini comme suit. Soit $g\in\GF$ et $h\in\bG$ tel que
$g=h\inv F(h)$. Alors, $\Sh((g))=(F(h) h\inv)$. On note encore
$\Sh$ l'automorphisme de l'espace des fonctions centrales sur
$\GF$ qui s'en d\'eduit.

Pour deux endomorphismes de Frobenius $F_1$ et $F_2$ de
$\bG$,
soit $\Sh_{F_1/F_2}$ la  descente de Shintani de  $\bG^{F_1}$ \`a
$\bG^{F_2}$ (\cf\ \cite[I, 7]{DM}). Pour $g\in\bG^{F_1}$, soit $h\in\bG$ tel que
$g=h\inv F_1(h)$. Alors, $\Sh_{F_1/F_2}((g))=(h F_2(h)\inv)$.
On note encore
$\Sh_{F_1/F_2}$ l'application lin\'eaire des fonctions centrales sur
$\bG^{F_1}$ vers celles sur $\bG^{F_2}$ qui s'en d\'eduit.

On note $T_w$ l'endomorphisme de $\Qlbar[\CB^F]=
\Ind_{\bB^F}^{\bG^F}\Id$ donn\'e
par $\bB\mapsto\sum_{\bB'}\bB'$, o\`u $\bB'$ d\'ecrit les \'el\'ements de
$\CB^F$ tels que $\bB'\xrightarrow{w} \bB$.

\begin{proposition}\label{lefschetz pour racines d-iemes}
Soit $\bw\in(\BW)^F$, $g\in\GF$ et soit $d$ un entier strictement
positif. Alors,
\begin{enumerate}
\item 
L'endomorphisme $g D_\bw F^n$ de $\bX(\bw^d)$
v\'erifie la formule des traces pour tout $n\ge0$.

\item
Si $\bw\in\bW^F$, on a
\begin{multline*}
g\mapsto|\bX(\bw^d)^{gD_\bw F^n}|)=\\
\begin{cases}\Sh^d(
g\mapsto\Trace(g  T_w   \mid  \Ind^\GF_{\bB^F}\Id))&\text{si $n=0$}\\
\Sh_{F^n/F}\circ\Sh_{F^{nd+1}/F^n}  (g\mapsto\Trace(g  T_w  F^n
\mid\Ind_{\bB^{F^{nd+1}}}^{\bG^{F^{nd+1}}}\Id))&\text{si $n>0$}.\\
\end{cases}\end{multline*}
\end{enumerate}
\end{proposition}

\begin{pf*}{Preuve}
L'\'el\'ement $\bw\in (\BW)^F$  peut s'\'ecrire sous la forme
$\bw_0^{I_1}\cdots\bw_0^{I_k}$ o\`u $I_1,\ldots,I_k$ sont des
sous-ensembles $F$-stables de $S$ (proposition \ref{monoidepointsfixes}).
Soit      
$$\bY=\bX(\uw_0^{I_1},\ldots,\uw_0^{I_k},\ldots,\uw_0^{I_1},
\ldots,\uw_0^{I_k})$$
  (la   suite  $\uw_0^{I_1},\ldots,\uw_0^{I_k}$  est
r\'ep\'et\'ee   $d$   fois),    une   compactification lisse  $F$-stable   de
$\bX(\bw^d)$ (propositions \ref{quasiaffine}, \ref{desingularisation} et
 corollaire \ref{O(W_I)lisse}).   

Soit $D'=D_{\uw_0^{I_1}\cdots\uw_0^{I_k}}$.
Le  morphisme $g D'  F^n$ est
fini et, comme il a une puissance  \'egale \`a une puissance de $F$, son
graphe est transverse \`a la diagonale.
Par cons\'equent, l'op\'erateur  $g   D'  F^n$   sur  $\bY$
v\'erifie la formule des traces (th\'eor\`eme \ref{Lefschetz}).
Notons qu'il prolonge l'endomorphisme $g D_\bw  F^n$ de $\bX(\bw^d)$.

Nous  consid\'erons  maintenant   la  stratification  de $\bY$ donn\'ee par
(\ref{lt ferme}). L'o\-p\'erateur $g D' F^n$  permute les pi\`eces
de cette stratification de la m\^eme mani\`ere que $D'$.
Nous allons montrer par r\'ecurrence
sur la  longueur de $\bw$ et  pour chaque $\bw$ par  r\'ecurrence sur la
dimension des  pi\`eces que $g D'  F^n$ v\'erifie la formule  des traces
sur toute union $D'$-stable de pi\`eces.

Les deux  membres de  la formule  des traces  pour un  endomorphisme $f$
d'une vari\'et\'e $\bX$ sont additifs par rapport \`a la d\'ecomposition
de $\bX$  comme union d'un  nombre fini de  sous-vari\'et\'es localement
ferm\'ees $f$-stables. Cela est clair pour le second membre et r\'esulte
pour le  premier membre de  la longue  suite exacte de  cohomologie pour
l'union  d'un ouvert  et  d'un  ferm\'e. Pour  prouver  que  $g D'  F^n$
v\'erifie la formule  des traces sur une union de  pi\`eces, on est donc
ramen\'e \`a le prouver pour une orbite de pi\`eces sous $D'$.

Sur une orbite qui a au moins deux  pi\`eces $g D' F^n$ n'a pas de point
fixe.  La  cohomologie d'une  telle  orbite  est  la somme  directe  des
cohomologies  des  pi\`eces  et  $g D'  F^n$  permute  cycliquement  les
facteurs de cette somme directe, donc  a une trace nulle et v\'erifie la
formule  des traces  dans  ce cas.  Si l'orbite  est  une seule  pi\`ece
$D'$-stable, cette pi\`ece est  de la forme $\bX((\bv_1\cdots\bv_k)^d)$,
o\`u $\bv_i\in W_{I_i}^F$, et $g  D' F^n$ induit $g D_{\bv_1\cdots\bv_k}
F^n$  sur cette  pi\`ece.  Si $\bv_1\cdots\bv_k$  n'est  pas \'egal  \`a
$\bw$, l'op\'erateur  $g D_{\bv_1\cdots\bv_k} F^n$ v\'erifie  la formule
des  traces  sur  cette  pi\`ece   par  r\'ecurrence  sur  $l(\bw)$.  Si
$\bv_1\cdots\bv_k=\bw$, par  r\'ecurrence sur la dimension  des pi\`eces
on  sait que  $g  D_\bw F^n$  v\'erifie  la formule  des  traces sur  la
r\'eunion de toutes les autres pi\`eces  de $\bY$. Comme il la v\'erifie
sur $\bY$, il la v\'erifie par additivit\'e sur $\bX(\bw^d)$.

D\'emontrons     maintenant    la deuxi\`eme assertion de la proposition.
Un    point    fixe    de    $gD_\bw    F^n$    sur    la    vari\'et\'e
$\bX(w,\ldots,w)$   (o\`u  $w$   est   r\'ep\'et\'e   $d$  fois)   est
une  suite   $(\bB_1,\ldots,\bB_d)$  telle   que  $$(\bB_1,\ldots,\bB_d)=
(\lexp{g}F^n(\bB_2),\ldots,\lexp{g}F^n(\bB_d),\lexp{g}F^{n+1}(\bB_1)).$$
Une telle
suite correspond \`a la     donn\'ee    de  $\bB_1$    tel     que
$\bB_1=\lexp{g^d}F^{dn+1}(\bB_1)$  et
$\lexp{g}F^n(\bB_1)\xrightarrow{w}\bB_1$.
Soit  $k$  tel  que $g^d=k\inv F^{nd+1}(k)$  et  posons 
$\bB'=\lexp k\bB_1$  et  $g'=kg F^n(k)\inv$. Alors les conditions deviennent
$\bB'\in\CB^{F^{nd+1}}$   et   $\lexp{g'}F^n(\bB')\xrightarrow{w}\bB'$.
On a 
$$\Trace(g'T_wF^n \mid\Ind_{\bB^{F^{nd+1}}}^{\bG^{F^{nd+1}}}\Id)=
\mid \{\bB'\in\CB^{F^{nd+1}}\mid \lexp{g'}F^n(\bB')\xrightarrow{w}\bB'\}.$$

Pour conclure, il faut donc montrer que 
si $n=0$  alors $(g')=\Sh^d((g))$  et   sinon
$(g')=\Sh_{F^n/F^{nd+1}}\circ\Sh_{F/F^n} (g)$.

Commen\c cons par  le cas  $n=0$. On a une bijection de l'ensemble des classes
de conjugaison rationnelles des conjugu\'es g\'eom\'etriques de $g$ vers
$H^1(F,C_\bG(g))$, donn\'ee par $\lexp  kg\mapsto k\inv F(k)$. Ici,
la classe  de $g'$  est param\'etr\'ee  par l'image  de $g^d$  dans
$H^1(F,C_ \bG(g))$. Or, d'apr\`es \cite[IV, proposition 1.1]{DM}, telle est l'image
de la classe de $g$ par $\Sh^d$.

Supposons  maintenant   $n  \ge  1$.   Soit  alors  $h\in\bG$   tel  que
$g=h F^n(h)\inv$.  Comme  $g\in\GF$,  on  a  $g^d=g F^n(g)\cdots
F^{(d-1)n}(g)=h F^{dn}(h)\inv$  et $g'=(kh)F^n(kh)\inv$.
Donc, en utilisant que  $\Sh_{F/F^n}(g) \subset\bG^{F^n}$, on a
\begin{multline*}
\Sh_{F/F^n}(g)= F^{nd}(\Sh_{F/F^n}(g))  =   F^{nd}(h\inv F(h))\\
=(h\inv g^d F^{nd+1}(h))   =  ( (kh)\inv F^{nd+1}(kh))
\end{multline*}
et $\Sh_{F^n/F^{nd+1}}$ de cette classe vaut bien $(g')$.
\end{pf*}

La  proposition  suivante nous  permettra  de  d\'emontrer que  certains
op\'e\-ra\-teurs  $D_\by$  ont  une  trace   nulle  sur  la  cohomologie  de
$\bX(\bw)$  (quand  nous  saurons  qu'ils  v\'erifient  la  formule  des
traces).
\begin{proposition}\label{divpi->trace nulle}
Soient $\bw\in\BW$ et $\by\in \End_{\cD^+}(\bw)$ tel que $\by\preccurlyeq\bpi$
et que l'image  $\beta(\by)\in W$ soit non  triviale. Alors l'endomorphisme
$D_\by$ de $\bX(\bw)$ n'a pas de points fixes.
\end{proposition}
\begin{pf*}{Preuve}
D'apr\`es   les   propri\'et\'es  des   mots   dans   $\BW$  (\cf\   par
exemple  \cite[3.20]{Sydney}),  le   fait  que  $\by$  divise  $\bpi$
est   \'equivalent    \`a  l'existence de $\bx,\bx'\in\bW$ tels que
$\by=\bx\bx'$.
Soit   $\by_1\cdots\by_n$  une   d\'ecomposition  de   $\by$  comme  dans la
remarque \ref{B+w->w'}(i)   o\`u   nous    supposons $\by_i\in\bW$.
Soit  $a\in\bX(\bw)$  et   $\bB=p'(a)$  (voir proposition \ref{transitif}). 
Alors, si  on   pose  $\bB_1=p'(D_{\by_1}(a))$,  $\bB_2=p'(D_{\by_1\by_2}(a))$,
\etc, un  point   fixe $a$ de   $D_\by$  donne   lieu   \`a  une   suite
$\bB\xrightarrow{y_1}\bB_1\xrightarrow{y_2}\bB_2\cdots\xrightarrow{y_n}\bB$, \ie\
\`a  un  \'el\'ement  $b\in\CO(\by)$ tel  que  $p'(b)=p''(b)$. 

Montrons
que  ceci   est  impossible.   Un tel \'el\'ement $b$ correspond, par  l'isomorphisme canonique
$\CO(\bx\bx')\iso\CO(x,x')$, \`a  une  suite
$\bB\xrightarrow x\bB'\xrightarrow{x'}\bB$.  Mais ceci  n'est possible  que si
$x'=x\inv$, ce que nous avons exclu en supposant que $\beta(\by)\ne 1$.
\end{pf*}

\sub{}

Nous  rappelons des  constructions  de \cite[\S 5.A et 6.D]{Sydney}.

\begin{proposition}\label{Db sur Xw}
Soit  $d$ un  entier  strictement positif.  
\begin{enumerate}
\item Soit   $\cD^+_{F^d}$  la  cat\'egorie   $\cD^+$  relative
au  Frobenius  $F^d$. On  a  $\End_{\cD^+_{F^d}}(\bpi)=C_\BW(F^d)$.  On
en  d\'eduit  une  action  \`a droite,  que  nous  noterons  $\bb\mapsto
D_\bb^{\bpi F^d}$, du mono{\"\i}de $C_\BW(F^d)$ sur $\bX(\bpi,F^d)$.
\item Soit  $\bw\in   B^+$   tel  que   $(\bw  F)^d=\bpi   F^d$.
L'application
$$\CO(\bw)  \xrightarrow i  \CO(\bw) \times  \CO(F(\bw))
\times \cdots \times\CO(F^{d-1}(\bw)),\  x\mapsto (x,F(x),\ldots,F^{d-1}(x))$$
 se restreint en  un plongement de $\bX(\bw)$
dans  $\bX(\bpi,F^d)$. Pour  tout $\bb\in  C_\BW(\bw F)$,  l'op\'erateur
$D^{\bpi  F^d}_\bb$  d\'efini  en  (i)  stabilise  $\bX(\bw)$,  et  nous
noterons $D_\bb$ sa restriction \`a $\bX(\bw)$. Ceci fournit une action
\`a droite
de $C_\BW(\bw F)$ sur $\bX(\bw)$ et une action \`a droite de
$C_B(\bw F)$ sur $H^*_c(\bX(\bw))$.
\end{enumerate}
\end{proposition}
\begin{pf*}{Preuve}  Tout  \'el\'ement  de   $\bW$  divisant  $\bpi$,  on  a
$C_\bW(\bpi   F^d)\subset\End_{\cD^+_{F^d}}(\bpi)\subset  C_{B^+}(\bpi
F^d)$, d'o\`u  (i) car  $C_\bW(\bpi F^d)=C_\bW(F^d)$  engendre $C_\BW(F^d)$
comme mono{\"\i}de (\cf\ proposition \ref{monoidepointsfixes}).

L'image      par     $i$      de      $\bX(\bw)$     est      clairement
dans     
$$\CO(\bw)     \times_\CB     \CO(F(\bw))     \times_\CB
\cdots   \times_\CB\CO(F^{d-1}(\bw))=\CO(\bw F(\bw)\cdots
F^{d-1}(\bw) =\bpi)$$
  d'o\`u l'assertion sur $i$  dans (ii). Enfin,
le  fait  que  $D^{\bpi  F^d}_\bb$ stabilise  $\bX(\bw)$  quand  $\bb\in
C_\BW(\bw F)$  n'est pas \'ecrit  dans \cite{Sydney}, mais on  peut le
voir  comme suit:  $\bX(\bw)$ est  caract\'eris\'e dans  $\bX(\bpi,F^d)$
comme  l'ensemble des  solutions de  l'\'equation $D_\bw  x=F(x)$.
Cette  \'equation est  clairement  pr\'eserv\'ee par  l'action d'un  tel
$D^{\bpi F^d}_\bb$.
\end{pf*}

Il est  clair que  dans le  cas o\`u  $\bb$ divise  $\bw$, l'op\'erateur
d\'efini en  (ii) ci-dessus  est bien  le m\^eme  que celui de la d\'efinition
\ref{Dt}, donc  quand $\bb\in  \End_{\cD^+}(\bw)$ c'est le  m\^eme qu'en
\S \ref{defDb}.

\sub{}
Nous  allons maintenant  introduire une  technique
qui  nous permettra
d'obtenir des r\'esultats pour les op\'erateurs
$D_\bx$ induits par les \'el\'ements $\bx\in\End_{\cD_I}(\bw)$ en nous
ramenant \`a un sous-groupe de Levi de type $I$.

\begin{definition}\label{Bw}
Soit $I$ une partie de $S$. Nous notons $\cD^+_I$ la
plus petite
sous-cat\'egorie de  $\cD^+$ avec  m\^eme ensemble d'objets et
telle que $\{\by\in B^+_I\mid\by\preccurlyeq\bw, \by\inv\bw F(\by)=\bw'\}
\subset\Hom_{\cD_I^+}(\bw,\bw')$.
De m\^eme, nous notons $\cD_I$ la plus petite sous-cat\'egorie de
$\cD$ contenant $\cD^+_I$ et dont toutes les fl\`eches sont inversibles.
\end{definition}

\begin{proposition}\label{Dy|pour y in B_I}
On se place dans le cadre des notations et des hypoth\`eses de la d\'efinition
\ref{Xtilde_n} et de la proposition \ref{X_w produit de varietes}.
Soit $\bx\in  B_I^+$ tel  que
$\bx\preccurlyeq\bw$.
Alors, on a un diagramme commutatif
$$\xymatrix{
\tilde\bX_{(I)}(\dz_1,\ldots,\dz_k)\times_{\bL_I^{\dz F}}
\bX_{\bL_I}(\by,\dz F)\ar[rr]^-\sim \ar[d]_{\Id\times D_\bx} &&
 \bX(\bw)\ar[d]^{D_{\bx}} \\
\tilde\bX_{(I)}(\dz_1,\ldots,\dz_k)\times_{\bL_I^{\dz F}}
\bX_{\bL_I}(\bx\inv\by\lexp{\bz} F(\bx),\dz F)\ar[rr]^-\sim &&
 \bX(\bx\inv\bw F(\bx)) 
}$$
\end{proposition}

\begin{pf*}{Preuve}
Notons tout d'abord que l'op\'erateur $D_\bx$ est bien d\'efini sur
$\bX_{\bL_I}(\by,\dz F)$ car $\bx\preccurlyeq\by$. Notons aussi que l'on a
$\bx\inv\by\lexp{z}F(\bx)=\alpha_I(\bx\inv\bw F(\bx))$ et
$\bz=\omega_I(\bx\inv\bw F(\bx))$ car
$\bx\inv\bw F(\bx)=(\bx\inv\by\lexp{\bz} F(\bx))\bz$.

Puisque $D_{\bx_1\bx_2}=D_{\bx_2}\circ D_{\bx_1}$ si
$\bx_1\bx_2\preccurlyeq\bw$, on peut supposer
$\bx\in\bW$. Quitte \`a changer de d\'ecomposition de $\by$, on peut en plus
supposer $\bx=\by_1$.
L'op\'erateur $\Id\times D_\bx$ envoie 
$$((g_1\bU_{I_1},\ldots,g_k\bU_{I_k}),(\bB_1,\ldots,\bB_h))$$
sur
$$((g_1\bU_{I_1},\ldots,g_k\bU_{I_k}),(\bB_2,\ldots,\bB_h,\lexp z F(\bB_1))).$$
Via les isomorphismes de la proposition \ref{X_w produit de varietes},
cela correspond \`a envoyer
\begin{multline*}
(\lexp{g_1}(\bB_1\bU_I),\lexp{g_1}(\bB_2\bU_I),\ldots,\lexp{g_1}(\bB_h\bU_I),
\lexp{g_1}(\lexp z F(\bB_1)\bU_{I_1}),
\lexp{g_2}(\lexp{z_2\cdots z_k} F(\bB_1)\bU_{I_2}),\\
\ldots,\lexp{g_k}(\lexp{z_k} F(\bB_1)\bU_{I_k}))
\end{multline*}
sur l'\'el\'ement $A$ suivant
de la vari\'et\'e
$\bX(y_2,\ldots,y_h,\lexp zF(y_1),z_1,\ldots,z_k)$:
\begin{multline*}
(\lexp{g_1}(\bB_2\bU_I),\ldots,\lexp{g_1}(\bB_h\bU_I),
\lexp{g_1}(\lexp z F(\bB_1)\bU_I),
\lexp{g_1}(\lexp z F(\bB_2)\bU_{I_1}),
\lexp{g_2}(\lexp{z_2\cdots z_k} F(\bB_2)\bU_{I_2}),\\
\ldots,\lexp{g_k}(\lexp{z_k} F(\bB_2)\bU_{I_k})).
\end{multline*}

Rappelons que 
$\lexp{\bz_i\cdots \bz_k}F(\by_1)\bz_i=\bz_i\lexp{\bz_{i+1}\cdots \bz_k}F(\by_1)
\in\bW$.
L'isomorphisme canonique
$$\CO(\lexp{z_i\cdots z_k}F(y_1),z_i)\iso
\CO(z_i,\lexp{z_{i+1}\cdots z_k}F(y_1))$$
envoie 
$$(\lexp{g_i}(\lexp{z_i\cdots z_k} F(\bB_1)\bU_{I_i}),
\lexp{g_i}(\lexp{ z_i\cdots z_k} F(\bB_2)\bU_{I_i}),
\lexp{g_{i+1}}(\lexp{ z_{i+1}\cdots z_k} F(\bB_2)\bU_{I_{i+1}}))$$
sur 
$$(\lexp{g_i}(\lexp{ z_i\cdots z_k} F(\bB_1)\bU_{I_i}),
\lexp{g_{i+1}}(\lexp{ z_{i+1}\cdots z_k} F(\bB_1)\bU_{I_{i+1}}),
\lexp{g_{i+1}}(\lexp{ z_{i+1}\cdots z_k} F(\bB_2)\bU_{I_{i+1}})).$$
L'isomorphisme canonique
$$\CO(\lexp{z_k}F(y_1),z_k)\iso
\CO(z_k,F(y_1))$$
envoie 
$$(\lexp{g_k}(\lexp{ z_k} F(\bB_1)\bU_{I_k}),
\lexp{g_k}(\lexp{ z_k} F(\bB_2)\bU_{I_k}),F(\lexp{g_1}(\bB_2\bU_I)))$$
sur
$$(\lexp{g_k}(\lexp{ z_k} F(\bB_1)\bU_{I_k}),F(\lexp{g_1}(\bB_1\bU_I)),
F(\lexp{g_1}(\bB_2\bU_I))).$$
Par composition des isomorphismes canoniques, on en d\'eduit que
$A$ s'envoie sur l'\'el\'ement
\begin{multline*}
(\lexp{g_1}(\bB_2\bU_I),\ldots,\lexp{g_1}(\bB_h\bU_I),
\lexp{g_1}(\lexp{ z} F(\bB_1)\bU_{I_1}),
\lexp{g_2}(\lexp{ z_2\cdots z_k} F(\bB_1)\bU_{I_1}),\\
\ldots,
\lexp{g_k}(\lexp{ z_k}
F(\bB_1)\bU_{I_k}),F(\lexp{g_1}(\bB_1\bU_I)))
\end{multline*}
de la vari\'et\'e
$\bX(y_2,\ldots,y_h,z_1,\ldots,z_k,F(y_1))$: c'est
bien l'image par $D_\bx$ de la suite initiale, d'o\`u le r\'esultat.
\end{pf*}

La proposition pr\'ec\'edente a pour cons\'equence:
\begin{corollaire}\label{Dx pour x in B+w,I}
Soit $\bx\in \End_{\cD_I}(\bw)$. Alors,
$D_\bx\in\End_\CC(\bX_\bw)$ correspond par
l'isomor\-phisme de la proposition \ref{X_w produit de varietes} \`a
$\Id\times D_\bx$ (o\`u $D_\bx$ est vu dans
$\End_\CC(\bX_{\bL_I}(\by,\dz F))$).
\end{corollaire}

Soit $\bL$  (resp. $\bU$)  un  compl\'ement de  Levi (resp.  le radical
unipotent) d'un sous-groupe parabolique de $\bG$ et soit $n\in\bG$ tel
que $\bL$ est $nF$-stable.
Soit $R^{\bG}_{\bL,\bU,n}:\CR(\bL^{nF})\to \CR(\GF)$ le morphisme donn\'e
par 
$$[V]\mapsto \sum_i (-1)^i [H^i_c(\tilde\bX^{\bL,\bU}(n),\Qlbar)
\otimes_{\Qlbar\bL^{nF}} V]$$
D'apr\`es \cite[Theorem 1.33]{BMM} la restriction aux
caract\`eres unipotents de $R^{\bG}_{\bL,\bU,n}$
ne d\'epend pas de $\bU$. Nous noterons alors 
$R^{\bG}_{\bL,n}$ cette restriction.

\begin{theoreme}\label{XhDyFn}
Soient $\bw\in B^+$, $I\subset S$ tels que $\lexp{\bw} F(B_I)=B_I$
et soit $\dz$ un repr\'esentant dans $N_\bG(\bT)$ de $\beta(\omega_I(\bw))$.
Soit $\bx\in \End_{\cD^+_I}(\bw)$.
Alors, la fonction centrale sur $\GF$ donn\'ee par $g\mapsto
\TrH{gD_\bx}{\bX(\bw)}$ est l'image par $R^{\bG}_{\bL_I,\dz}$
de la fonction
$$l\mapsto\TrH{lD_\bx}{\bX_{\bL_I}(\alpha_I(\bw),\dz F)}.$$
\end{theoreme}
\begin{pf*}{Preuve}
Le th\'eor\`eme r\'esulte des propositions \ref{X_w produit de varietes}
et \ref{Dy|pour y in B_I}, compte tenu du fait
qu'on peut donc,
si la forme normale de $\omega_I(\bw)$ est $\bz_1\cdots\bz_k$,
calculer $R^\GF_{\bL_I,\dz}$
en utilisant la vari\'et\'e $\tilde\bX_{(I)}(\dz_1,\ldots,\dz_k)$,
o\`u les $\dz_i$ sont des repr\'esentants dans $N_\bG(\bT)$ des $z_i$
tels que $\dz=\dz_1\cdots\dz_k$.
\end{pf*}

\subsection{Cas $\bw=\bpi^n$}
\sub{}

Soit $(W,S)$ un groupe de Coxeter cristallographique muni d'un automorphisme
$\sigma$.
Soit $\CH_x(W,\sigma)$ le quotient de l'al\-g\`ebre
du groupe $C_{B^+}(\sigma)$ sur l'anneau $\bbZ[x^{1/2},x^{-1/2}]$
par l'id\'eal engendr\'e par les
$(\bw_0^I+1)(\bw_0^I-x^{l(w_0^I)})$ pour ${I\in S/\sigma}$.
On note $T_I$ l'image de  $\bw_0^I$ dans $\CH_x(W,\sigma)$. Cette alg\`ebre
admet pour base $\{T_w\}_{w\in  W^\sigma}$ et se  sp\'ecialise en
l'alg\`ebre de groupe  $\Qlbar  W^\sigma$  par $x^{1/2}\mapsto  1$ (rappelons
que $C_B(\sigma)$ est un mono{\"\i}de de tresses, \cf\ proposition
\ref{monoidepointsfixes}).
Pour $\chi\in\Irr(W^\sigma)$,
nous  notons $\chi_x\in\Irr(\CH_x(W,\sigma))$ le caract\`ere  qui  se
sp\'ecialise en  $\chi$ par  $x^{1/2}\mapsto 1$. 

On suppose maintenant \`a nouveau que $W$ est le groupe de Weyl de $\bG$.
Soit $\CH_q(W,F)=\CH_x(W,F)\otimes_f\Qlbar$,
o\`u $f$ est la sp\'ecialisation $x^{1/2}\mapsto q^{1/2}$.

On prend comme convention que les endomorphismes d'un espace vectoriel
agissent \`a droite (donc, si $M$ est un $(A,B)$-bimodule, alors on
a un morphisme d'alg\`ebres $B\to\End_A(M)$).

\begin{proposition}\label{cH(q)}
On a 
$\CH_q(W,F)\simeq\End_{\Qlbar\GF}(\Ind_{\bB^F}^\GF\Id)$.
\end{proposition}
\begin{pf*}{Preuve}
Cette propri\'et\'e est bien connue, mais faute de r\'ef\'erence commode
nous en rappelons une d\'emonstration. On sait que
$\End_{\Qlbar\bG^F}(\Ind_{\bB^F}^{\bG^F}\Id)$
est une alg\`ebre de Hecke du groupe $W^F$, o\`u les valeurs
propres du g\'en\'erateur associ\'e \`a $w_0^I$ sont $-1$ et
$|(\bB/(\bB\cap\lexp{w_0^I}\bB))^F|$
(\cf\  \cite[IV, \S 1, exercice 24]{Bbki}). 
Pour calculer cette derni\`ere valeur propre on utilise le fait que
$(\bB/(\bB\cap\lexp{w_0^I}\bB))^F\simeq(\bU\cap\lexp{w_0^I}\bU^-)^F$, o\`u
$\bU^-$ est le radical unipotent du sous-groupe de Borel oppos\'e \`a $\bB$.
Comme le  nombre de points rationnels  sur $\Fq$ d'un espace  affine est
ind\'ependant  de la  $\Fq$-structure consid\'er\'ee  (voir par  exemple
\cite[Example 3.7]{DMb}), le nombre de points rationnels vaut
$q^{\dim(\bU\cap\lexp{w_0^I}\bU^-)}=q^{l(w_0^I)}$, ce
qui est bien la valeur annonc\'ee.
\end{pf*}

\sub{}
Dans le cas o\`u l'action de $F$ sur $W$ est triviale et $n=1$,
le th\'eor\`eme suivant est \cite[th\'eor\`eme 2.7]{Sydney}.

\begin{theoreme}\label{action sur Xpi}
Pour tout $n\ge 1$,
l'action  de   $C_{B^+}(F)$  sur  $\bX(\bpi^n)$  (\cf\   proposition
\ref{Db  sur  Xw})
induit  un  morphisme  
$$\CH_q(W,F)\to
\End_{\Qlbar\GF}(\bigoplus_iH^i_c(\bX(\bpi^n))),\ T_I\mapsto D_{\bw_0^I}.$$
\end{theoreme}
\begin{pf*}{Preuve}
Il  suffit  de voir  que  les  $D_{\bw_0^I}$ v\'erifient  les  relations
quadratiques   de  $\CH_q(W,F)$.
Nous  appliquons  la proposition \ref{Dy|pour
y  in  B_I}  avec $I$  une  partie  $F$-stable  de  $S$ (les
hypoth\`eses sont  v\'erifi\'ees pour tout $\bx\in  B_I^+$ tel
que $\bx\preccurlyeq\bpi^n$). Si $\bx\in
C_{B^+_I}(F)$ et $\bx\preccurlyeq\bpi^n$,  alors  $D_\bx$  est un
endomorphisme  de  $\bX(\bpi^n)$
qui  correspond d'apr\`es  la proposition \ref{Dy|pour  y in  B_I} \`a 
l'endomorphisme
$D_\bx$ de $\bX_{\bL_I}(\bpi^n_I)$,  o\`u $\bpi_I=(\bw_0^I)^2$ est
l'\'el\'ement
analogue  \`a $\bpi$  pour  $\bL_I$ (car  on a  $\alpha_I(\bpi^n)=\bpi_I^n$,
et  $\dz=1$). Via la  formule
de  K\"unneth,  il  suffit   donc  de  d\'emontrer  que  l'endomorphisme
$D_{\bw_0^I}$ de $\bX_{\bL_I}(\bpi_I^n)$ v\'erifie la relation quadratique
$(D_{\bw_0^I}+1)(D_{\bw_0^I}- q^{l(w_0^I)})=0$ quand  $I$ est une orbite
de $F$ dans $S$.

Soit $I$  une  orbite  de $F$ sur $S$. Le   groupe    adjoint
$(\bL_I)_{\ad}$   de   $\bL_I$   est   une   descente   des   scalaires,
c'est-\`a-dire  de   la  forme   $\bH^k$  o\`u   $\bH$  est   un  groupe
r\'eductif   quasi-simple,   et  $(\bL_I)_{\ad}^F\simeq\bH^{F^k}$.   Les
repr\'esentations unipotentes  se factorisant par le  groupe adjoint, la
cohomologie  de  $\bX_{\bL_I}(\bpi_I^n)$  est  isomorphe  \`a  celle  de
$\bX'=\bX_\bH(\bpi^{kn}_\bH,F^k)$  (o\`u  $\bpi_\bH$  est  l'\'el\'ement
analogue  \`a $\bpi$  pour  $\bH$) sur  laquelle $D_{\bw_0^I}$  induit
l'action de  $D_{\bw_0^\bH}$ (proposition \ref{produit}).

Le groupe  $\bH^{F^k}$ est d'un  des types
$A_1(q^k)$, $\lexp  2A_2(q^k)$, $\lexp 2B_2(q^k)$ ou  $\lexp 2G_2(q^k)$.
D'apr\`es les th\'eor\`emes \ref{2A2}, \ref{2B2} et \ref{2G2}, seules
les repr\'esentations $\Id$ et $\St$ apparaissent dans la cohomologie
de $\bX'$.
D'apr\`es les propositions \ref{Id} et \ref{Steinberg},
la  vari\'et\'e $\bX'$  n'a que  deux groupes  de
cohomologie non  nuls: en degr\'es $2knl(\bpi_\bH)$  et $knl(\bpi_\bH)$.
Sur $H_c^{2knl(\bpi_\bH)}(\bX')$,  le groupe  $\bH^{F^k}$ agit  par $\Id$
et la valeur propre de $F^k$ est $q^{k^2nl(\bpi_\bH)}$. Sur
$H_c^{knl(\bpi_\bH)}(\bX')$, il agit par $\St$ et la valeur propre de
$F^k$ est $1$.  Puisque les $(\Qlbar\bG^F)$-modules 
$H_c^{2knl(\bpi_\bH)}(\bX')$ et $H_c^{knl(\bpi_\bH)}(\bX')$ sont
irr\'eductibles, l'action de $D_{\bw_0^\bH}$ sur chacun d'eux est donn\'ee par
un scalaire.

Comme $D_{\bw_0^\bH}^{2kn}=D_{\bpi_\bH^{kn}}=F^k$ sur
$\bX'$, on en  d\'eduit que $D_{\bw_0^\bH}$ agit sur
$H_c^{knl(\bpi_\bH)}(\bX')$ (un  espace de dimension $q^{kl(w_0^\bH)}$) par
$\alpha$, une  racine  $2kn$-\`eme  de l'unit\'e.
En outre, $D_{\bw_0^\bH}$ agit sur $H_c^{2knl(\bpi_\bH)}(\bX')$ (un espace
de dimension $1$)
par $q^{kl(w_0^\bH)}\beta$, o\`u $\beta$ est une racine  $2kn$-\`eme de
l'unit\'e.

Les propositions \ref{lefschetz  pour  racines  d-iemes} (i)  et
\ref{divpi->trace nulle} donnent  $\TrH{D_{\bw_0^\bH}}   {\bX'}=0$,
donc $\alpha=-\beta$.
Soit $g$   est  un  \'el\'ement  unipotent non trivial
de  $\bH^{F^k}$. Alors,  $\Trace(g,\St)=0$   et $\Trace(g,\Id)=1$.
D'apr\`es la proposition \ref{lefschetz pour    racines     d-iemes} (i), 
on a $|(\bX')^{gD_{\bw_0^\bH}}|=
\TrH{gD_{\bw_0^\bH}}{\bX'}=\beta   q^{kl(\bw_0^\bH)}$, donc
$\beta=1$ puisque $|(\bX')^{gD_{\bw_0^\bH}}|$ est  un entier positif. On
a donc bien la relation quadratique annonc\'ee.
\end{pf*}

\subsection{Cas $\bw_0\bpi^n$}

Notons $\CH_{-q}(W,{w_0 F})$ la sp\'ecialisation de
$\CH_x(W,{w_0 F})$ via $x^{1/2}\mapsto i\sqrt q$.
dans cette partie, nous construisons une repr\'esentation de l'alg\`ebre
$\CH_{-q}(W,{w_0 F})$ dans $\End_\GF(\oplus_iH^i_c(\bX(w_0)))$.

Soit   $I\in   S/w_0   F$;   comme   $(\bw_0\bpi^n   F)^2=\bpi^{2n+1}F^2$
et  $\bw_0^I\in   C_{B^+}(\bw_0\bpi^n  F)$, la proposition
\ref{Db  sur  Xw}  (ii) fournit un  endomorphisme  $D_{\bw_0^I}$
de  $\bX(\bw_0\bpi^n)$   pour  tout   $n\in\bbN$, donnant naissance
\`a une action de $C_{B^+}(\bw_0 F)$ sur $\bX(\bw_0\bpi^n)$.

Les deux \'enonc\'es suivants
g\'en\'eralisent \cite[3.10 (b)]{LuMa} (cas o\`u $w_0F$ agit trivialement
sur $W$).

\begin{theoreme}\label{End(H(X_w0))} Pour tout $n\ge 0$, l'action de
$C_{B^+}(\bw_0 F)$ sur $\bX(\bw_0\bpi^n)$ induit un morphisme d'alg\`ebres
$$\CH_{-q}(W,{w_0 F})\to\End_{\Qlbar\GF}(\oplus_iH^i_c(\bX(\bw_0\bpi^n))),\
T_I\mapsto (-1)^{l(w_0^I)}D_{\bw_0^I}.$$
\end{theoreme}
\begin{pf*}{Preuve}
On proc\`ede comme dans la preuve du th\'eor\`eme \ref{action sur Xpi}.
Il s'agit de prouver que pour tout $I$, l'action de
$(-1)^{l(w_0^I)}D_{\bw_0^I}$
sur la cohomologie de $\bX(\bw_0\bpi^n)$
v\'erifie la m\^eme relation quadratique que $T_I$.
Gr\^ace au corollaire \ref{Dx pour x in B+w,I}, il suffit de
montrer que $(-1)^{l(w_0^I)}D_{\bw_0^I}$ v\'erifie la bonne
relation quadratique sur la vari\'et\'e
$\bX_{\bL_I}(\bw_0^I\bpi_I^n,w_0^Iw_0F)$.

Soit $\bH$ 
un groupe r\'eductif quasi-simple tel que le groupe adjoint
de $\bL_I$ soit isomorphe \`a $\bH^k$ et soit
$F'$ l'endomorphisme de Frobenius sur
$\bH$ correspondant par cet isomorphisme \`a $(w_0^Iw_0F)^k$.
On se ram\`ene \`a l'\'etude de
l'op\'erateur $D_{\bw_0^\bH}$, induit par $D_{\bw_0^I}$, sur la vari\'et\'e
$\bX'=\bX_\bH((\bw_0^\bH\bpi_\bH^n)^k,F')$.
Le groupe $\bH^{F'}$ est d'un des types $A_1(q^k)$, $\lexp2B_2(q^k)$,
$\lexp2G_2(q^k)$, $\lexp2A_2(q^k)$ (si $k$ est pair) ou $A_2(q^k)$
(si $k$ est impair).
D'apr\`es les th\'eor\`emes \ref{2B2}, \ref{2G2}, \ref{2A2} et \ref{A2},
seules les repr\'esentations $\Id$ et $\St$ apparaissent dans la cohomologie
de la vari\'et\'e $\bX'$.
D'apr\`es les propositions \ref{Id} et \ref{Steinberg},
la vari\'et\'e $\bX'$
n'a que deux groupes de cohomologie non nuls: en degr\'es
$kl(\bw_0^\bH\bpi_\bH^n)$ et
$2kl(\bw_0^\bH\bpi_\bH^n)$. Le groupe $\bH^{F'}$ agit sur
$H_c^{kl(\bw_0^\bH\bpi_\bH^n)}(\bX')$
par $\St$
la valeur propre de $F'$ est 1. Sur
$H_c^{2kl(\bw_0^\bH\bpi_\bH^n)}(\bX')$, le groupe $\bH^{F'}$ agit 
par $\Id$ et $F'$ a comme valeur
propre $q^{k^2l(\bw_0^\bH\bpi^n_\bH)}$.

L'op\'erateur $D_{\bw_0^\bH}$ est une racine
$(2n+1)k$-\`eme de $F'$ et commute \`a $\bH^{F'}$. Par cons\'equent,
il agit sur $H_c^{kl(\bw_0^\bH\bpi^n)}(\bX')$ (un espace de dimension
$q^{kl(w_0^\bH)}$) par $\alpha$ et sur $H_c^{2k{l(\bw_0^\bH\bpi^n)}}(\bX')$
(un espace de dimension $1$) par
$\beta q^{kl(w_0^\bH)}$,
o\`u $\alpha$ et $\beta$ sont des racines $(2n+1)k$-\`emes de 1,

Les propositions \ref{lefschetz  pour  racines  d-iemes} (i) et
\ref{divpi->trace nulle} donnent   $\TrH{D_{\bw_0^\bH}}   {\bX'}=0$,
donc $\alpha=(-1)^{{kl(\bw_0^\bH)}-1}\beta$.
Soit $g$   est  un  \'el\'ement  unipotent non trivial
de  $\bH^{F'}$.
D'apr\`es la proposition \ref{lefschetz pour    racines     d-iemes} (i), 
on a
\begin{multline*}
|(\bX')^{gD_{\bw_0^\bH}}|=\TrH{gD_{\bw_0^\bH}}{\bX'}\\
=-\beta \St(g)+\beta q^{kl(\bw_0^\bH)}=\beta q^{kl(\bw_0^\bH)},
\end{multline*}
donc $\beta=1$ puisque $|(\bX')^{gD_{\bw_0^\bH}}|$ est  un entier positif. On
a donc bien la relation quadratique annonc\'ee.
\end{pf*}

\begin{remarque}
On pourra pr\'ef\'erer consid\'erer le morphisme $T_I\mapsto D_{\bw_0^I}$ o\`u
l'alg\`ebre de Hecke de $W^{w_0F}$ a pour param\`etres
$\{((-1)^{1+l(w_0^I)},q^{l(w_0^I)})\}_{I\in S/w_0F}$.
\end{remarque}

L'\'enonc\'e  pr\'ec\'edent  a  la  cons\'equence   suivante.

\begin{corollaire}\label{valeurs   propres   dans   H(w0)}   Pour   tout
$n\in  \bbN$,  les  valeurs  propres  de  $F$  dans  la  cohomologie  de
$\bX(\bw_0\bpi^n)$ sont  soit de la  forme $\pm  q^m$, soit de  la forme
$\pm i q^{m+\frac12}$ o\`u $m\in\bbN$. On est toujours dans le premier cas
si $w_0$  est central dans  $W$ et $\bG$ n'a  pas de composante  de type
$E_7$ ou $E_8$.
\end{corollaire}
\begin{pf*}{Preuve}     L'action de $F$     sur     $\oplus_i
H^i_c(\bX(\bw_0\bpi^n))$ est \'egale \`a l'action de
l'\'el\'ement $(-1)^{l(w_0)} T_{w_0}^{2n+1}$
dans la  repr\'esentation de $\CH_{-q}(W,{w_0 F})$  sur cette cohomologie.
Pour  $\chi\in\Irr(W^{w_0 F})$,  nous  notons  $\chi_{-q}$ le  caract\`ere
correspondant  de  $\CH_{-q}(W,{w_0  F})$ et  $d_{\chi_{-q}}$  le  degr\'e
g\'en\'erique de $\chi_{-q}$. Enfin  nous noterons $A_\chi$ (resp. $a_\chi$)
le  degr\'e (resp.  la  valuation) de  $d_{\chi_{-q}}$.  La valeur  propre
de  $T^2_{w_0}$  dans la  repr\'esentation  de  caract\`ere $\chi_{-q}$  est
$(-q)^{2N-a_{\chi}-A_{\chi}}$  \cite[corollaire   4.20]{Sydney}  et  les
valeurs propres de  $T_{w_0}$ sont des racines carr\'ees  de cette valeur
(ici, $N$ est le nombre de racines de $W$).
Ceci d\'emontre  l'\'enonc\'e sauf  la derni\`ere phrase.  Si $w_0$  est central
alors pour  tout $\chi$ irr\'eductible on  a $\chi(w_0)\ne 0$ et  si $\bG$
n'a pas de composante de type  $E_7$ ou $E_8$ alors par \cite[corollaire
4.19]{Sydney}  on  en  d\'eduit   que  $a_{\chi}+A_{\chi}$  est  pair  (en
effet tout  caract\`ere est  alors \og g\'en\'eriquement rationnel\fg\  suivant la
terminologie de {\it loc. cit.}).
\end{pf*}

\end{document}